\documentclass[1p,number,times,sort&compress]{elsarticle}
\usepackage[dvipsnames,table,xcdraw]{xcolor}
\usepackage[normalem]{ulem} 
\usepackage{cancel}
\usepackage{graphicx} 
\usepackage{subfigure} 
\usepackage{chngcntr} 
\usepackage{amssymb} 
\usepackage{amsmath} 
\usepackage{mathrsfs}
\usepackage{epstopdf}
\usepackage{float} 
\usepackage{booktabs} 
\usepackage{microtype}
\usepackage[font=footnotesize]{caption}
\usepackage[colorlinks,breaklinks]{hyperref}
 
\newcommand{\Rn}{\mathbb{R}^3}
 
\makeatletter
\renewcommand{\p@subfigure}{\thefigure}
\makeatother

\begin{document} 
 

\title{Building Three-Dimensional Differentiable Manifolds Numerically}

\author[cass,ccom]{Lee Lindblom\corref{cor1}}
\ead{llindblom@ucsd.edu}
\author[HTWBerlin]{Oliver Rinne}
\ead{oliver.rinne@htw-berlin.de}
\author[Cornell]{Nicholas W.~Taylor}
\ead{nwt2@cornell.edu}

\cortext[cor1]{Corresponding author}

\address[cass]{Center
  for Astrophysics and Space Sciences, University of California at San
  Diego,\\ 9500 Gilman Drive, La Jolla, CA 92093, USA}
\address[ccom]{Center for Computational Mathematics,
  University of California at San Diego,
  \\ 9500 Gilman Drive, La Jolla, CA 92093, USA}
\address[HTWBerlin]{Faculty 4, HTW Berlin -- University of Applied Sciences,
Treskowallee 8, 10318 Berlin, Germany} 
\address[Cornell]{Department of Physics, Cornell University,
Ithaca, NY 14853, USA}

\date{\today}
 
\begin{abstract}
A method is developed here for building differentiable
three-dimensional manifolds on multicube structures. This method
constructs a sequence of reference metrics that determine
differentiable structures on the cubic regions that serve as
non-overlapping coordinate charts on these manifolds. It uses
solutions to the two- and three-dimensional biharmonic equations in a
sequence of steps that increase the differentiability of the reference
metrics across the interfaces between cubic regions.  This method is
algorithmic and has been implemented in a computer code that
automatically generates these reference metrics.  Examples of
three-manifolds constructed in this way are presented here, including
representatives from five of the eight Thurston geometrization
classes, plus the well-known Hantzsche-Wendt, the Poincar\'e
dodecahedral space, and the Seifert-Weber space.
\end{abstract}


\maketitle

\section{Introduction}
\label{s:Introduction}

Differentiable manifolds are the mathematical structures on which the
differential equations of the physical sciences are solved to provide
descriptions of the universe as we understand it.  This paper develops
methods that allow these equations to be solved numerically in a
convenient way on a much broader class of manifolds.

In the traditional literature, an $n$-dimensional differentiable
manifold is defined as a space that can be covered by a collection of
open sets, plus invertible maps that take each member of this
collection onto some open subset of $\mathbb{R}^n$.  In practical
terms, these open subsets in $\mathbb{R}^n$ are the coordinate charts
used to identify points in the manifold.  For points having images in
two coordinate patches, the inferred maps in the overlap regions
between the patches must be differentiable.  The differentiability
of these overlap maps defines the differentiable structure of the
manifold.  This structure is used to define what it means for global
tensor fields on the manifold to be continuous and differentiable.
The existence of differentiable global tensor fields is fundamental to
finding global solutions to the equations of the physical sciences on
manifolds.  Therefore having, or if necessary creating, a suitably
smooth differentiable structure on a manifold is essential.

The traditional description of a differentiable manifold is difficult
to implement numerically in a computer code for several reasons: Such
an implementation must keep track of the exact size and shape of each
coordinate patch in $\mathbb{R}^n$, plus the exact sizes and shapes of
the overlap regions containing points represented in two patches, plus
the maps between the coordinates in the overlap regions.  These
structures can of course be designed and implemented in a code for any
particular manifold. However, each case is unique and each case
requires a lot of work to design and implement properly.  It requires
a great deal of effort even to transform a numerical code designed for
use on one manifold into one that can be used on another.  In
addition, there does not exist in the literature (so far as we know) a
catalog containing the needed information (i.e. the needed collections
of coordinate regions, plus all the needed information about their
overlaps, plus the maps between the overlap regions) that would allow
these traditional methods to be implemented in a code in a
straightforward way for a broad collection of three-dimensional
manifolds.

An alternative description of a differentiable manifold was introduced
in Ref.~\cite{Lindblom2013} that is simpler in ways that make it more
suitable for use in a computer code.  In this multicube approach the
coordinate charts in $\mathbb{R}^n$ are standardized, requiring each
patch to be a cube of uniform coordinate size and orientation.  These
coordinate patches are chosen not to overlap in $\mathbb{R}^n$, except
for points on the boundaries of the cubes.  The global coordinates in
$\mathbb{R}^n$ can therefore be used to identify points globally in
these manifolds.  Since the coordinate patches have uniform sizes and
shapes in this approach, the maps that identify points on the
boundaries between neighboring patches are particularly simple,
consisting of a rigid translation that maps the center of a face into
the center of its neighbor's face, followed by a simple rotation
(and/or reflection) that aligns the two faces in the appropriate way.
In three dimensions, the case of primary interest in this paper, the
number of possible rotations/reflections is quite small (just 48), so
all the possible maps are easily included in a computer code.  It was
shown in Ref.~\cite{Lindblom2013} that this multicube structure is
sufficiently general to represent any two- or three-dimensional
manifold in this way.

The simplicity of the structures of the coordinate charts and their
overlap regions makes it much easier to implement the multicube
description of a manifold in a computer code.  In addition, describing
manifolds in this way makes it possible to access and easily make use
of published catalogs that contain thousands of three-dimensional
manifolds represented by their triangulations~\cite{Matveev2005,
  Martelli2001, Martelli2006, Burton2011}.  Some of these catalogs
include online access to the explicit triangulations for these
manifolds~\cite{Regina}.  Converting a triangulation into a multicube
structure is straightforward, see e.g.  Ref.~\cite{Lindblom2013}.  A
computer code that implements this procedure has been developed as
part of this project and is described in some detail in
\ref{s:TriangulationToMultiCubeCode}. Most of the manifolds
included in this study are based on triangulations given in
Ref.~\cite{Regina}, and then converted to multicube structures by
this new code.  The basic multicube structures constructed in this
way do not come with differentiable structures.  So the problem of
constructing those differentiable structures--the main focus of this
paper--remains.

Since the coordinate patches in a multicube representation do not
overlap, it is not possible to construct differentiable structures on
these manifolds in the traditional way.  Instead,
Ref.~\cite{Lindblom2013} showed how these structures could be
constructed using a reference metric.  Given a reference metric that
is continuous across each interface boundary in a multicube
structure, a simple analytical formula can be used to determine
special Jacobians at those boundaries.  Those Jacobians can then be
used to define what it means for vector and tensor fields to be
continuous across those boundaries.  A reference metric that is both
continuous and differentiable (in the appropriate sense) across the
interfaces can also be used to define a covariant derivative that
(together with the Jacobians) can be used to determine what it means
for vector and tensor fields to be differentiable across those
boundaries.  This approach was used to construct differentiable
structures on a few simple three-dimensional manifolds in
Ref.~\cite{Lindblom2013}.  An algorithmic method for constructing the
needed reference metrics numerically for arbitrary two-dimensional
manifolds was developed and tested in Ref.~\cite{Lindblom2015}.  This
paper focuses on the more difficult and complicated problem of
developing analogous algorithmic methods for constructing reference
metrics on arbitrary three-dimensional manifolds.

Most of the equations of the physical sciences require fixing some
combination of the values and normal derivatives of the fields at the
boundaries of computational domains.  This means that a differentiable
structure must be present on the manifold that is capable of defining
what it means for fields and their derivatives to be continuous across
those boundaries.  For a manifold constructed by the multicube method,
this means that a global $C^{\,1}$ metric is required.  The purpose of
this paper is to develop a step-by-step algorithm for constructing
global $C^{\,1}$ metrics on these manifolds.  These steps consist of
building a sequence of metrics $\hat g_{ab}$, $\bar g_{ab}$,
$\bar{\bar g}_{ab}$, and $\tilde g_{ab}$ described in detail in
Secs.~\ref{s:C0ReferenceMetrics} and
\ref{s:C1ReferenceMetrics}.\footnote{The notation $g_{ab}$ is often
used to represent the physical metric (as determined by solving
Einstein's equation for example). To avoid confusion, that notation is
not used here in the construction of the $C^1$ reference metric
$\tilde g_{ab}$ that is designed only to define the differential
structure of the manifold.} The first part of this procedure,
described in Sec.~\ref{s:C0ReferenceMetrics}, constructs a global
$C^{\,0}$ metric, $\hat g_{ab}$, whose intrinsic parts (i.e., the
components that define the intrinsic metric on a given face) are
continuous across the interface boundaries between the cubic regions,
and which is free from conical singularities at the vertices and along
the edges of those regions.  The first step, described in
Sec.~\ref{s:Step1}, re-organizes the multicube structure into a set of
overlapping star-shape domains that surround each of the vertices in
the multicube structure.  Singularity-free flat metrics are
constructed on these star-shaped domains in the second step, described
in Sec.~\ref{s:Step2}.  These flat metrics are combined together using
a special partition of unity to produce a global $C^{\,0}$ reference
metric, $\hat g_{ab}$, in the third step, described in
Sec.~\ref{s:Step3}.

In Sec.~\ref{s:C1ReferenceMetrics} the $C^{\,0}$ metric, $\hat
g_{ab}$, is transformed into a $C^{\,1}$ metric $\tilde g_{ab}$ in
three additional steps. These steps build two additional intermediate
metrics, $\bar g_{ab}$ and $\bar{\bar g}_{ab}$ in
Secs.~\ref{s:Step3.1} and ~\ref{s:Step3.2}. In the first of these,
Sec.~\ref{s:Step3.1}, a conformal transformation is applied to $\hat
g_{ab}$ that produces a new metric, $\bar g_{ab}$, that makes all the
edges of each cubic region into geodesics. This transformation also
makes one component of the associated extrinsic curvatures $\bar
K^{\{\alpha\}}_{ab}$ vanish along the edges.  The conformal factor
needed for this step is produced by solving two-dimensional biharmonic
equations on each cube face, with boundary conditions along the edges
that enforce the geodesic conditions. The pseudo-spectral numerical
methods used to solve those equations for this study are described in
\ref{s:BiharmonicMethods}.  In Sec.~\ref{s:Step3.2} gauge
transformations are performed on the metric $\bar g_{ab}$ at the
interfaces of the cubic regions.  The resulting metric $\bar{\bar
  g}_{ab}$ has the property that its intrinsic components on each cube
face are identical to those of $\bar g_{ab}$, but the gauge components
of the metric on those faces are deformed in a way that makes all the
components of the associated extrinsic curvatures $\bar{\bar
  K}^{\{\alpha\}}_{ab}$ vanish on all the edges of each cubic region.
In Sec.~\ref{s:Step3.3} the metric $\bar{\bar g}_{ab}$ is adjusted in
the interiors of each cubic region (keeping the boundary values fixed)
by solving three-dimensional biharmonic equations whose boundary
conditions are chosen to make the extrinsic curvatures $\tilde
K^{\{\alpha\}}_{ab}$ vanish on each cube face.  This $\tilde g_{ab}$
retains the continuity of its intrinsic components across each
interface boundary inherited from $\hat g_{ab}$ and $\bar g_{ab}$.
The continuity of the intrinsic metric together with the continuity of
the extrinsic curvature are the geometric conditions, often referred
to as the Israel junction conditions~\cite{Israel1966}, needed to
ensure that the metric $\tilde g_{ab}$ is $C^{\,1}$ across the
interface boundaries.
 
Section~\ref{s:NumericalExamples} describes a number of
three-dimensional manifolds on which $C^{\,1}$ differentiable
structures have been constructed for this study using the methods
described in Secs.~\ref{s:C0ReferenceMetrics} and
\ref{s:C1ReferenceMetrics}.  Numerical convergence of the Israel
junction conditions, the necessary and sufficient conditions that
$\tilde g_{ab}$ be $C^{\,1}$ across the interface boundaries, is
demonstrated for these examples. \ref{s:3DMulticubeManifolds} presents
detailed multicube structures for a variety of three-dimensional
manifolds, including examples from the Thurston geometrization
classes~\cite{Thurston1997,Scott1983} $E^3$, $S^3$, $S^2\times S^1$,
$H^2\times S^1$, and $H^3$.  The manifolds studied here include 29
that were constructed from triangulations given in Ref.~\cite{Regina}
using the code described in \ref{s:TriangulationToMultiCubeCode}.  In
addition a few multicube structures were constructed by hand for
several well known three-manifolds: including the Poincar\'e
dodecahedral space~\cite{ThrelfallSeifert1931}, Seifert-Weber
space~\cite{SeifertWeber1933}, and all six compact orientable
three-manifolds that admit flat metrics~\cite{riazuelo2004,
  Hitchman2018}, including the Hantzsche-Wendt
space~\cite{Hantzsche1934}.  Section~\ref{s:Discussion} gives a brief
summary of the basic methods developed in this paper and the ways they
have been tested numerically.  In addition a number of interesting
questions and possible extensions of the current results are outlined.


\section{Constructing $C^{\,0}$ Three-Dimensional Reference Metrics}
\label{s:C0ReferenceMetrics}

The procedure to create a continuous ($C^{\,0}$) three-dimensional
reference metric, $\hat g_{ij}$, on a multicube structure has three
basic steps: In the first, described in Sec.~\ref{s:Step1}, the
multicube structure is re-organized to create a collection of
overlapping star-shaped domains on the manifold.  In the second step,
described in Sec.~\ref{s:Step2}, flat metrics are constructed in each
of these overlapping domains.  In the third step, described in
Sec.~\ref{s:Step3}, a global $C^{\,0}$ reference metric, $\hat
g_{ab}$, is constructed using these flat metrics and a special
partition of unity.  Explicit analytic formulas are given in
Secs.~\ref{s:Step2} and \ref{s:Step3} for the $C^{\,0}$ metric, $\hat
g_{ab}$, along with the flat metrics and partition of unity functions
used to construct it.

All these steps can be, and have been, implemented in a computer code
that automatically generates these $C^{\,0}$ metrics using only the
multicube structures as input.  In the simplest version of this
procedure (the one described in most detail here, and the one
presently implemented in our code) all the dihedral angles between the
cube faces that meet along a particular edge are chosen to have the
same size.  While this simplifying assumption cannot be applied to
most multicube structures, it is general enough that compliant
structures have been constructed here on a diverse set of manifolds in
Sec.~\ref{s:NumericalExamples} to illustrate these methods.

\subsection{Step 1: Assembling Star-Shaped Domains.}
\label{s:Step1} 

\begin{figure*}[!h]  
\centering
\subfigure{
  \includegraphics[height=0.23\textwidth]{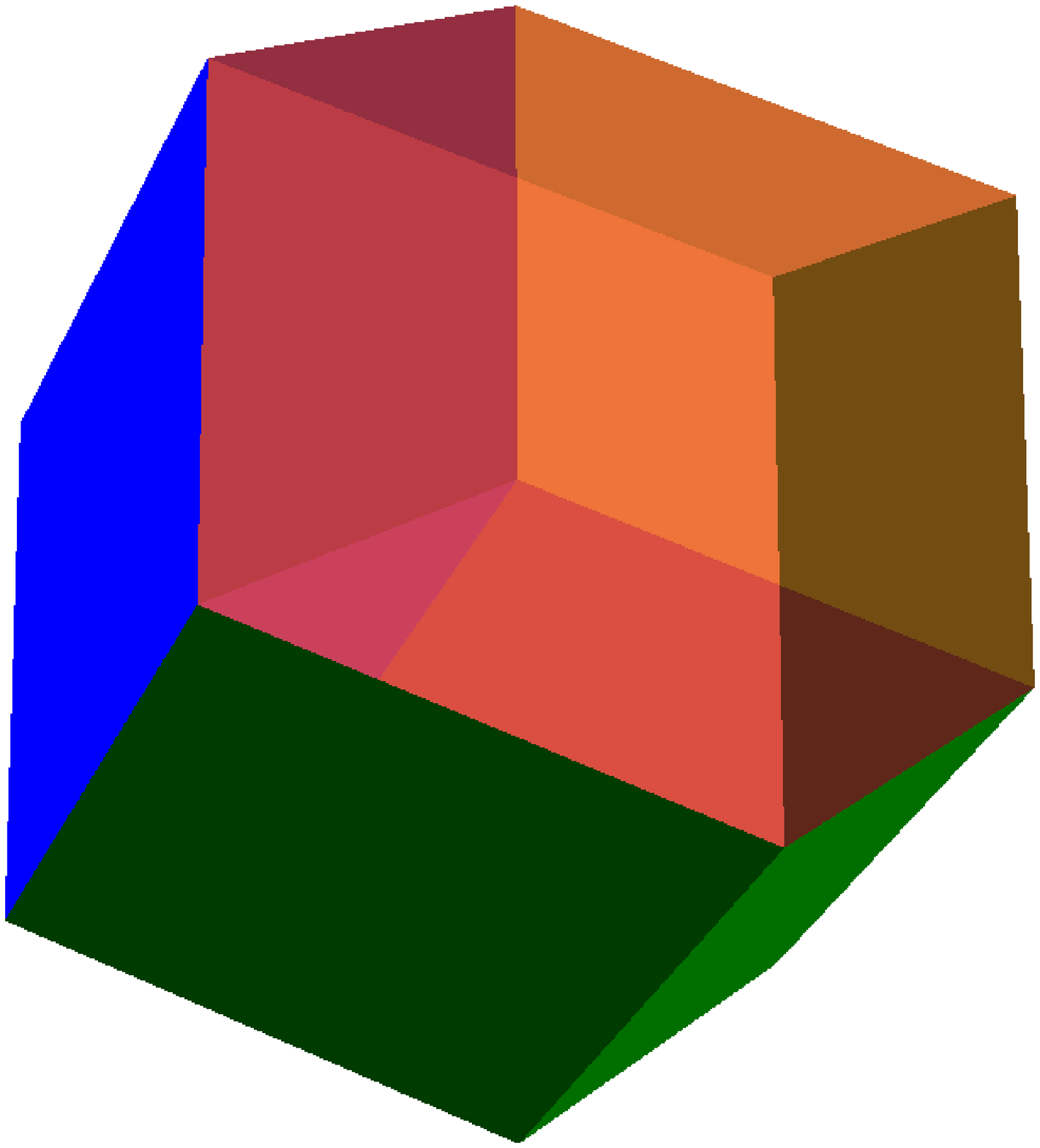}
}
\subfigure{
  \includegraphics[height=0.23\textwidth]{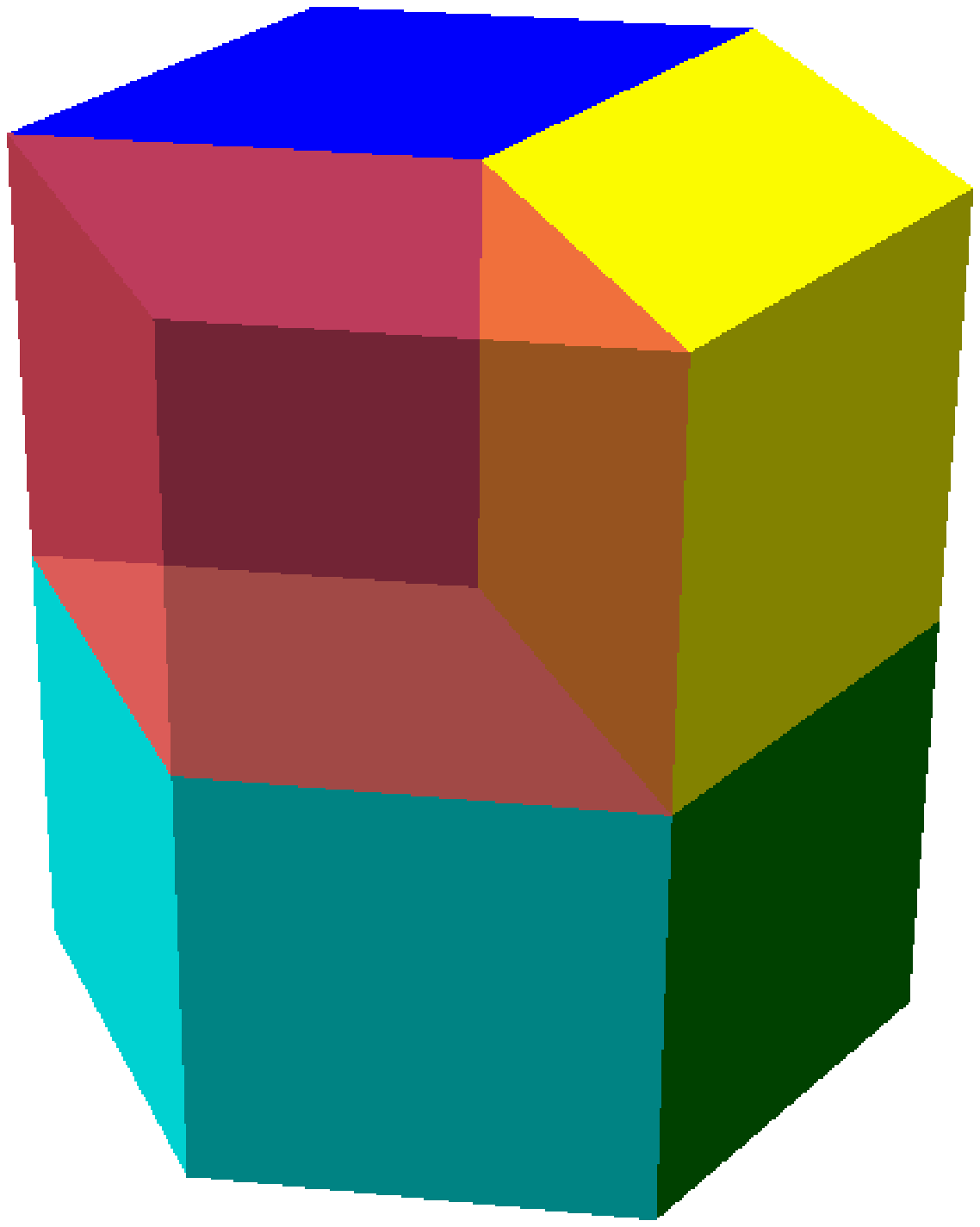}
}
\subfigure{
  \includegraphics[height=0.23\textwidth]{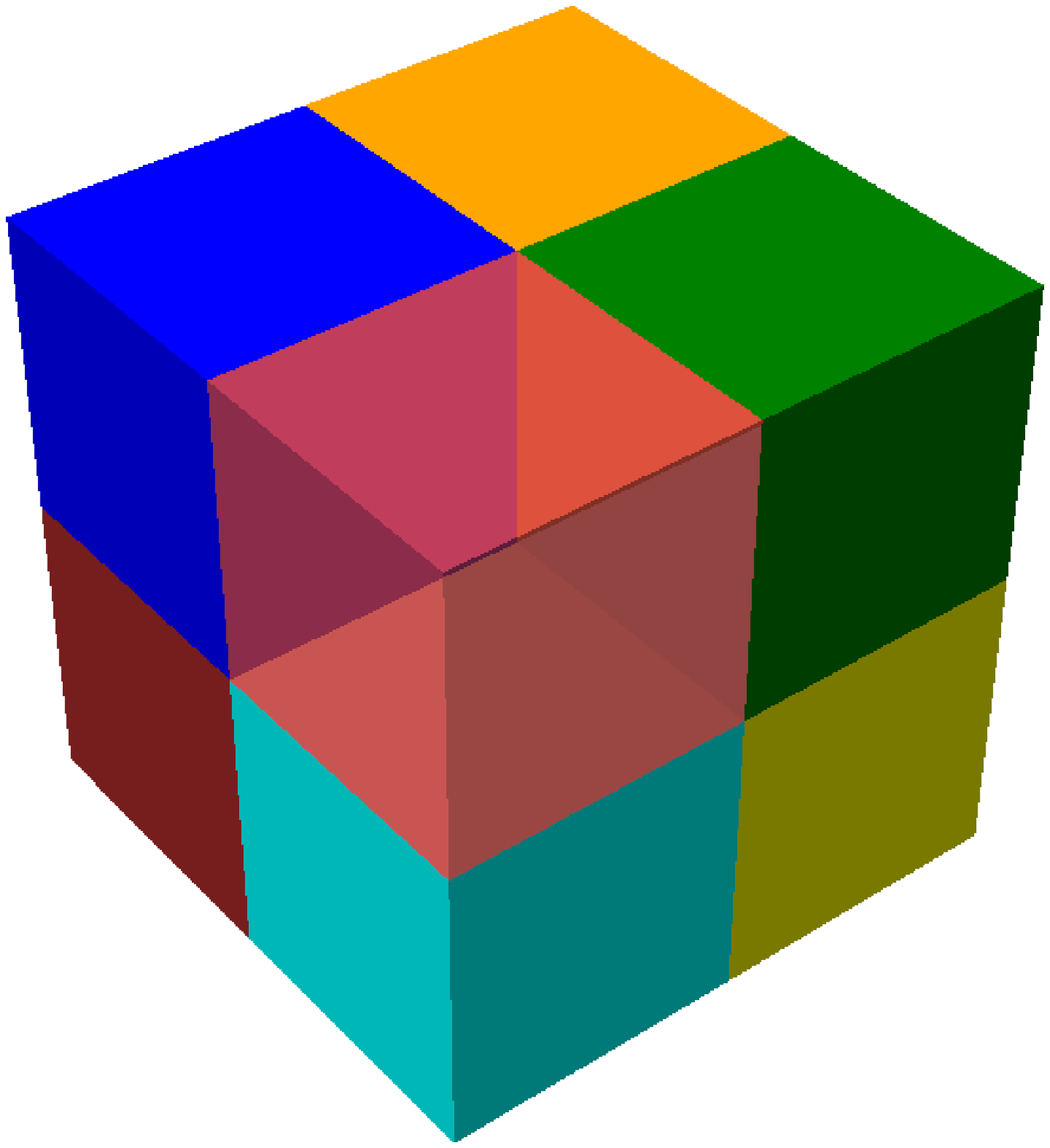}
} 
\subfigure{
  \includegraphics[height=0.23\textwidth]{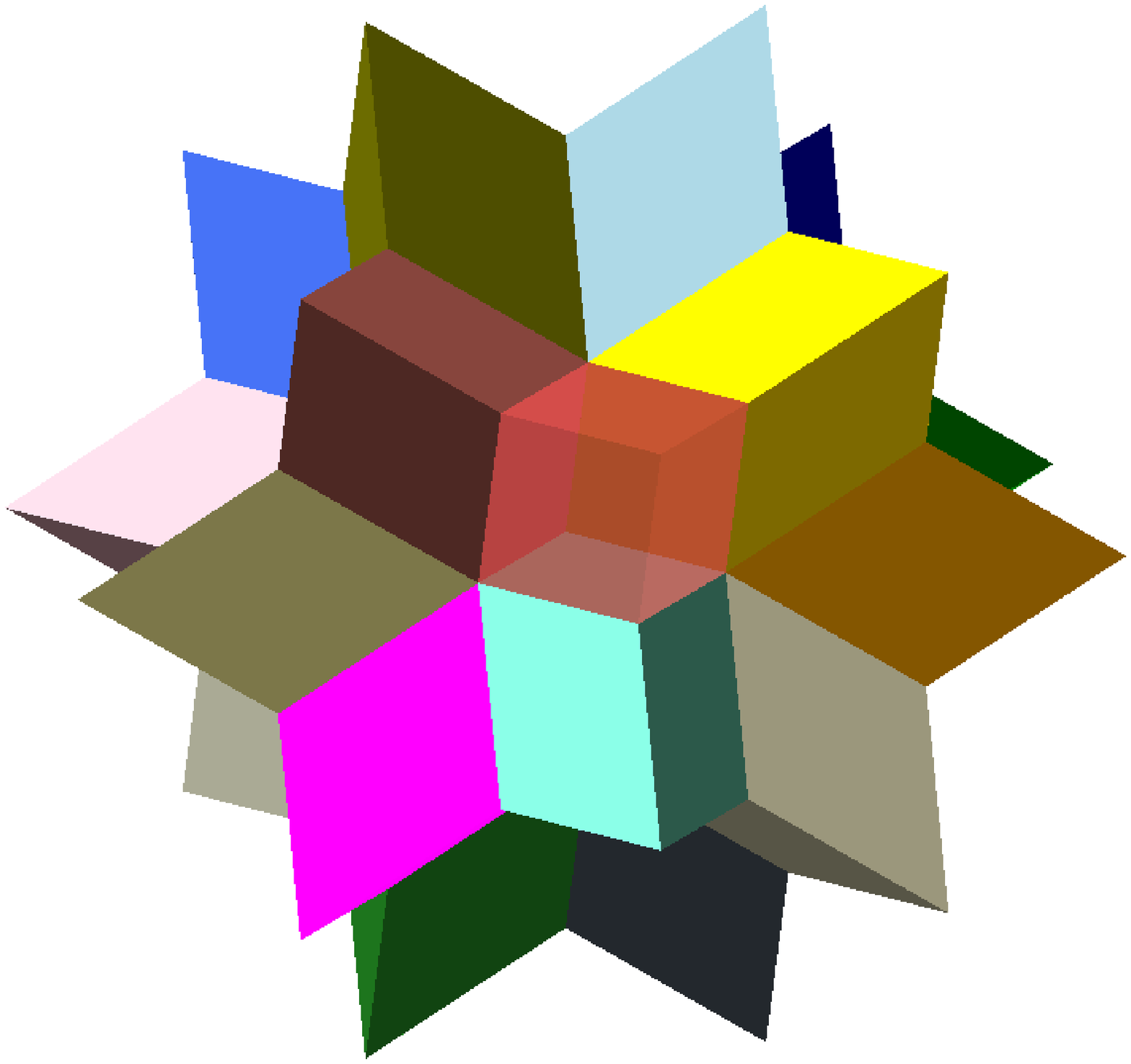}
}
\caption{\label{f:Star_Shaped_Domains} Examples of star-shaped
  domains, $\mathcal{S}_I$, in three dimensions consisting of four,
  six, eight and twenty cubic regions, respectively, that intersect at
  their central vertex points.  The cubic regions in each example have
  been distorted so they fit together smoothly with the flat metric of
  the $\Rn$ in which they are shown.  One (red colored) cubic region
  in each example has been made semi-transparent to allow the internal
  structures of these domains to be seen more clearly.}
\end{figure*}

In this first step, the multicube structure consisting of a collection
of cubic regions, $\mathcal{B}_A$, is enhanced by defining a set of
domains, called the star-shaped domains, $\mathcal{S}_I$, that overlap
the boundaries between the primary cubic regions.  One star-shaped
domain surrounds each distinct vertex of the multicube structure. It
is constructed from (copies of) all the cubic regions that intersect
at that vertex point.  (A particular cubic region $\mathcal{B}_A$ may
be included more than once in a star-shaped domain if two or more of
its vertices are identified with each other.)  Each of the star-shaped
domains, $\mathcal{S}_I$, has the topology of an open ball in $\Rn$.
The index {\scriptsize $A$} is used to label the cubes $\mathcal{B}_A$
in the multicube structure, while the index {\scriptsize $I$} labels
the star-shaped domains $\mathcal{S}_I$, or equivalently the distinct
vertices in the multicube structure.  The structures of the individual
star-shaped domains depend on the global properties of the multicube
structure, in particular on how many cube vertices intersect in the
manifold at the center of each $\mathcal{S}_I$.
Figure~\ref{f:Star_Shaped_Domains} illustrates several examples of
star-shaped domains having different numbers of cubic regions
intersecting at their central vertex points.

A code designed to use multicube structures can be enhanced to
assemble the $\mathcal{S}_I$ in a fairly straightforward way: Any
multicube structure code must include the cube face identification
maps.  Starting at one vertex of one cubic region, the identities of
the three cubes whose faces are identified with the faces of
$\mathcal{B}_A$ adjacent to this vertex are determined from the
multicube maps.  This can be done, for example, by following the
interface identification maps for points near this vertex on each of
the three faces that meet at that point. Copies of the three cubes
identified as neighbors in this way are added to $\mathcal{S}_I$.
This identification step is repeated for the adjacent faces of each of
the additional cubes, and then iterated until (copies of) all the cube
vertices that intersect the original vertex point are included in
$\mathcal{S}_I$.  Once a star-shaped domain $\mathcal{S}_I$ is
complete, if some cube vertices in the full multicube structure
remain un-assigned to the center of some star-shaped domain, then a
new star-shaped domain $\mathcal{S}_{I+1}$ is constructed around this
vertex using the same procedure.  The process terminates when all the
cube vertices have been included at the center of some star-shaped
domain.  There are a finite number of cube vertices in any multicube
structure (that can be used for practical numerical work), so in
practice this process always terminates after a finite number of
steps.

\subsection{Step 2: Constructing Semi-Local Flat Metrics.}
\label{s:Step2}

The second step in the procedure to construct global $C^{\,0}$
reference metrics builds a flat metric in each of the star-shaped
domains, $\mathcal{S}_I$, introduced in Sec.~\ref{s:Step1}.  Each
$\mathcal{S}_I$ consists of a cluster of cubes that intersect at its
central point.  If these cubes are appropriately distorted into
parallelograms (by adjusting the dihedral angles between the cube
faces), they can be fitted together (without overlapping and without
leaving gaps between them) to form an isometric subset of $\Rn$, and
thus inherit a natural flat metric.
Figure~\ref{f:Star_Shaped_Domains} illustrates several simple examples
of star-shaped domains isometrically embedded in $\Rn$.

\begin{figure}[!htb] 
\centering 
\subfigure{
  \includegraphics[width=0.3\textwidth]{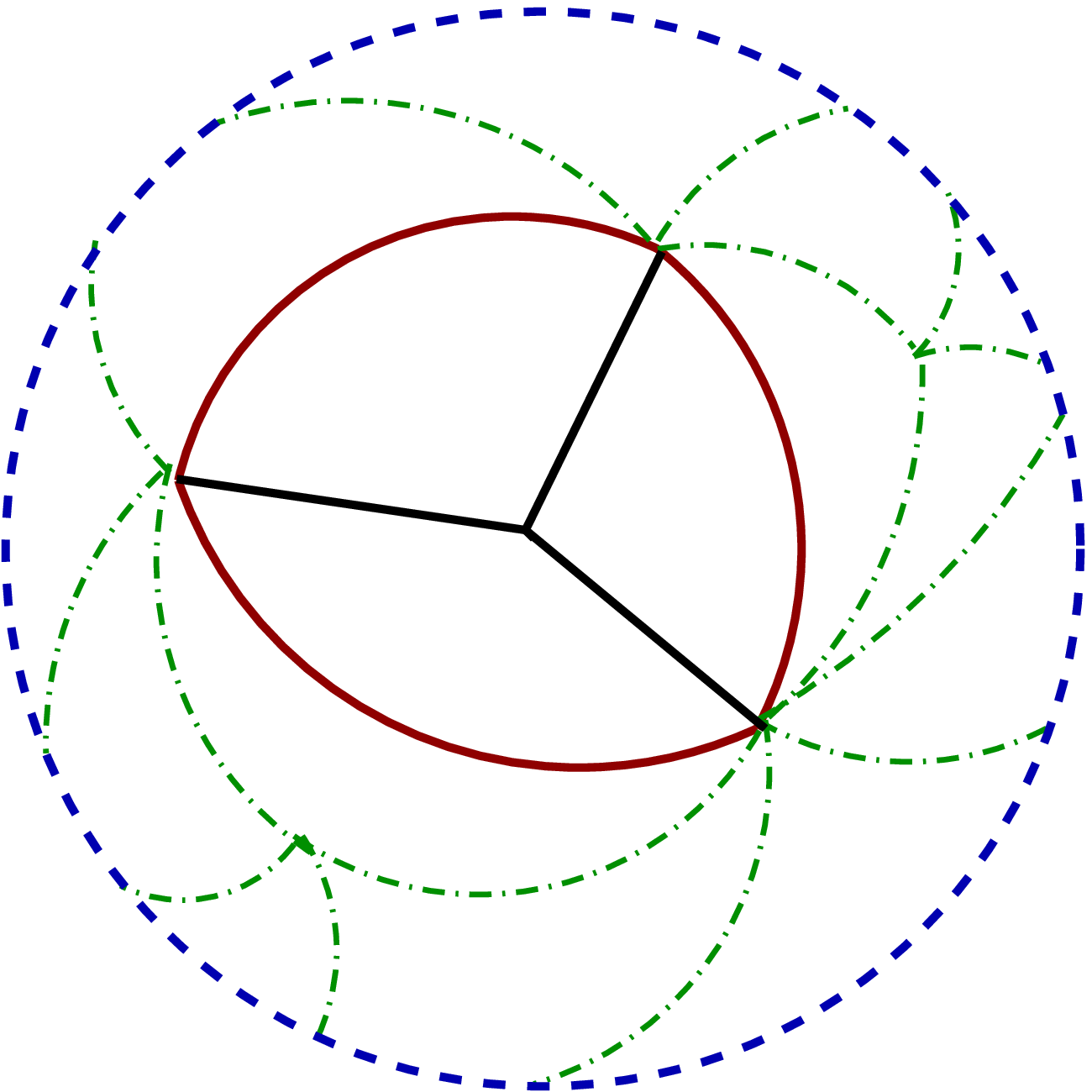}
  \label{f:ThreeDWedgeA}}
\hspace{0.04\textwidth}
\subfigure{\vspace{0.1cm}
  \includegraphics[width=0.3\textwidth]{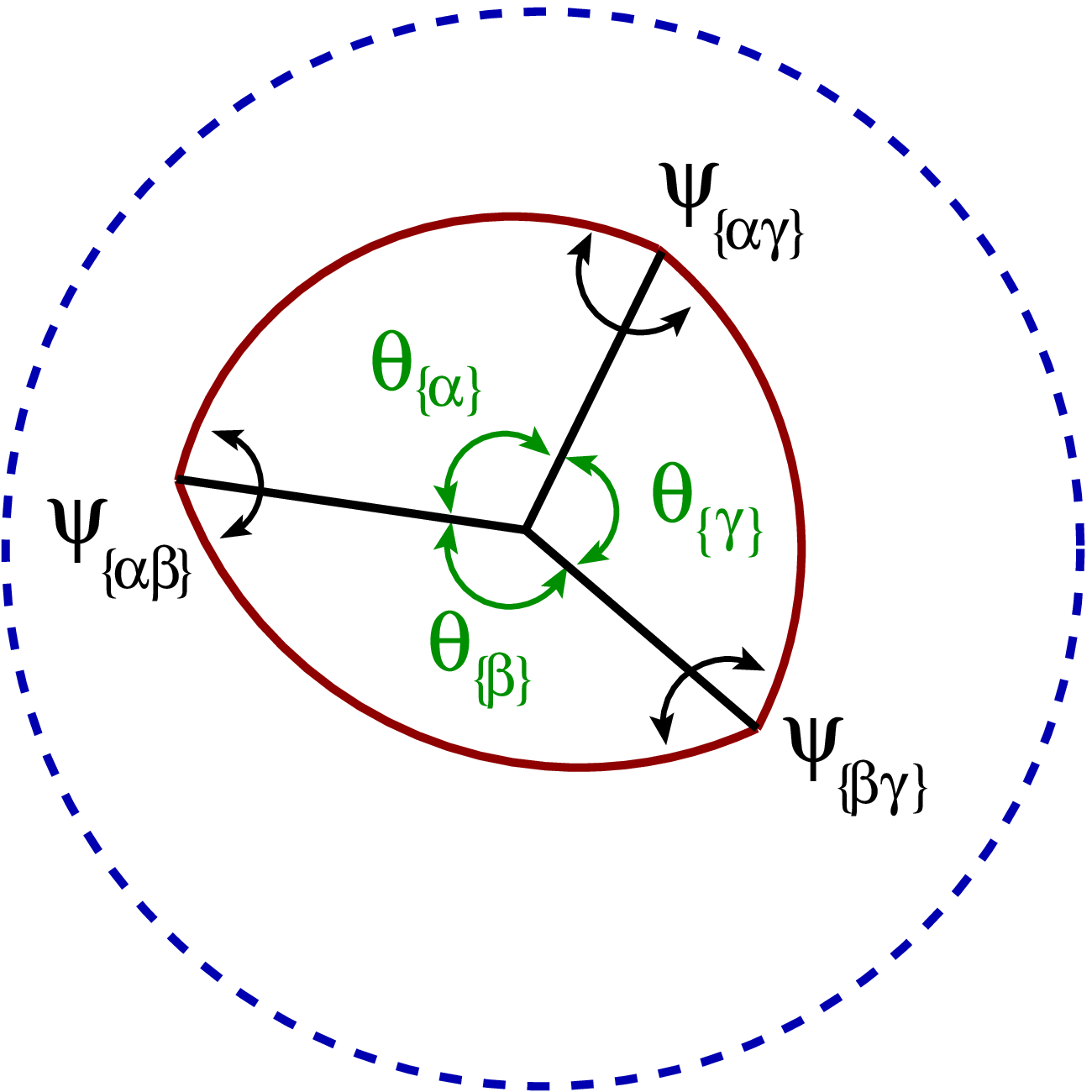}
           \label{f:ThreeDWedgeB}}
\hspace{0.04\textwidth}

\caption{Left illustration, \ref{f:ThreeDWedgeA}, shows the
  intersection between the corner of a cubic region and a small sphere
  centered on the vertex of one of the star-shaped domains.  This
  sphere is depicted as the dashed (blue) curve, the intersections
  between this cubic region and the sphere are shown as solid (red)
  curves.  The dash-dot (green) curves represent the intersections of
  nearby cubic regions in the star-shaped domain.  Right illustration,
  \ref{f:ThreeDWedgeB}, labels the angles that characterize the
  spherical triangle formed by the intersection of a cubic region and
  a small sphere centered at its vertex. The $\psi_{\{\alpha\beta\}}$,
  etc.  are the dihedral angles (in the local flat metric) between the
  faces of this cubic region.  These $\psi_{\{\alpha\beta\}}$ are also
  the angles of the spherical triangle.  The $\theta_{\{\alpha\}}$,
  etc. are the angles between the edges of the cubic region.  These
  $\theta_{\{\alpha\}}$ are also the arc lengths of the sides of the
  spherical triangle.
\label{f:ThreeDWedge} }
\end{figure}

To understand whether the cubic regions of a multicube structure can
always be fitted together in this way, consider a small two-sphere
surrounding the central vertex of $\mathcal{S}_I$.  This sphere
intersects all the cubic regions that meet at this central point. The
intersections of this sphere with the faces and edges of each cube
form triangles on this sphere. Figure~\ref{f:ThreeDWedgeA} illustrates
the spherical triangles that result from these intersections.  The
intersection of one of these cubes, $\mathcal{B}_A$, is displayed as
the spherical triangle with solid (red) line edges.  The intersections
of other nearby cubes in $\mathcal{S}_I$ are displayed with dash-dot
(green) line edges.  Together the intersections from all the cubic
regions in $\mathcal{S}_I$ form a triangulation of this two-sphere.

Any triangulation on a two-sphere can be realized geometrically in an
infinite number of ways.  Given any one realization, an infinite
number of others can be created simply by moving the vertices of the
triangles around on the sphere by small amounts (i.e., much smaller
than the sizes of the triangles), and then replacing the edges with
geodesics (great circles) between vertices.  Each spherical triangle
with geodesic edges represents the intersection of a parallelogram
(whose dihedral angles match the angles of the triangle) with the
two-sphere.  Thus there are an infinite number of ways to construct
distorted parallelograms that fit together in the correct way to
represent $\mathcal{S}_I$ as an isometric subset of $\Rn$.

An algorithm designed to compute a flat metric on $\mathcal{S}_I$ must
choose from among the infinite possibilities in some way.  Making and
implementing that choice is expected to be a complicated optimization
problem that we plan to analyze fully in a future study.  For the
purposes of the present study, however, we have chosen to adopt a
simple pragmatic approach: choosing the dihedral angles to have
uniform sizes around each edge.  This simple approach limits the class
of multicube structures to which it can be applied.  However, it is
general enough that we have been able to construct compliant examples
(see Sec.~\ref{s:NumericalExamples}) from most of the Thurston
geometrization classes, plus examples of several well known manifolds
like the Poincar\'e dodecahedral space~\cite{ThrelfallSeifert1931},
Seifert-Weber space~\cite{SeifertWeber1933}, and all six compact
orientable three-manifolds that admit flat metrics~\cite{riazuelo2004,
  Hitchman2018}, including the Hantzsche-Wendt
space~\cite{Hantzsche1934} (E6).

Before proceeding with the details of constructing flat metrics on the
$\mathcal{S}_I$ in these simple multicube structures, it will be
helpful to establish some basic notation.  The notation
$\partial_\alpha\mathcal{B}_A$ (or more compactly $A\{\alpha\}$) is
used to refer to the $\alpha$ face of cubic region
$\mathcal{B}_A$.  The index $\alpha$ can have the values $\{-x$, $+x$,
$-y$, $+y$, $-z$, $ +z\}$.  The edge of region $\mathcal{B}_A$ formed
by the intersection of the $A\{\alpha\}$ and $A\{\beta\}$ faces is
referred to as $A\{\alpha\beta\}$, and the vertex formed by the
intersections of the $A\{\alpha\}$, $A\{\beta\}$, and $A\{\gamma\}$
faces is referred to as $A\{\alpha\beta\gamma\}$.  The dihedral angle
between the $A\{\alpha\}$ and $A\{\beta\}$ faces is denoted
$\psi_{A\{\alpha\beta\}}$, while the angle between the axes at the
edges of the $A\{\alpha\}$ face is denoted $\theta_{A\{\alpha\}}$.
The $\psi_{A\{\alpha\beta\}}$ are also equal to the angles of the
spherical triangle created by the intersection of cube $\mathcal{B}_A$
with a small sphere (see Fig.~\ref{f:ThreeDWedgeB}), and the
$\theta_{A\{\alpha\}}$ are also equal to the arc lengths of the edges
of this triangle.

The uniform dihedral angle spacing assumption adopted here requires
the dihedral angles of all the cubic regions that intersect along an
edge to be the same.  In addition to being reasonably simple to
impose, it has the advantage of imposing a rigid uniformity that
prevents any cubic region from being more distorted than its
neighbors.  To prevent conical singularities along the cube edges, the
sum of the dihedral angles around each edge must be exactly $2\pi$.
The uniform dihedral angle assumption therefore implies that the
dihedral angle at the $A\{\alpha\beta\}$ edge must be given by
\begin{equation}
  \psi_{A\{\alpha\beta\}} = \frac{2\pi}{K_{A\{\alpha\beta\}}},
  \label{e:DihedralAngleDef}
\end{equation}
where $K_{A\{\alpha\beta\}}$ is the number of cubic regions that
intersect along this edge.

\begin{figure}[!htb] 
\centering 
\subfigure{\vspace{0.1cm} 
           \includegraphics[height=0.3\textwidth]{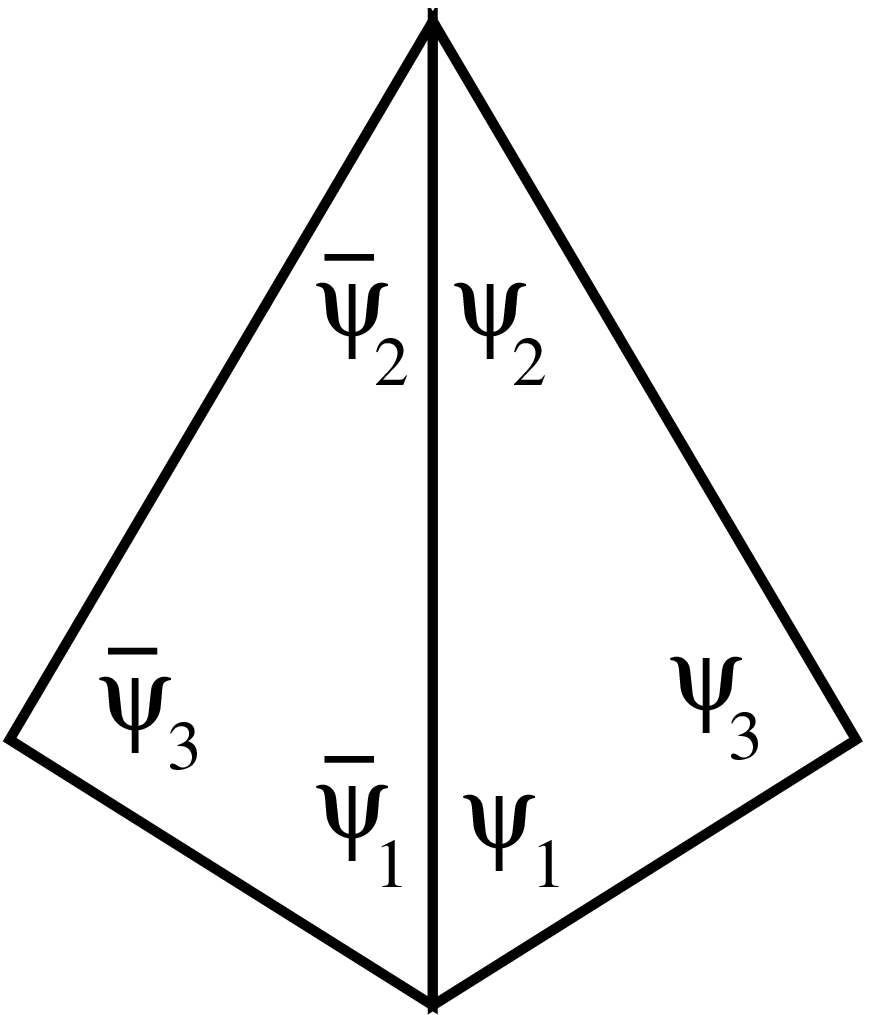}
           \label{f:AnglesA}}
\hspace{0.1cm}
\subfigure{\vspace{0.1cm} 
           \includegraphics[height=0.3\textwidth]{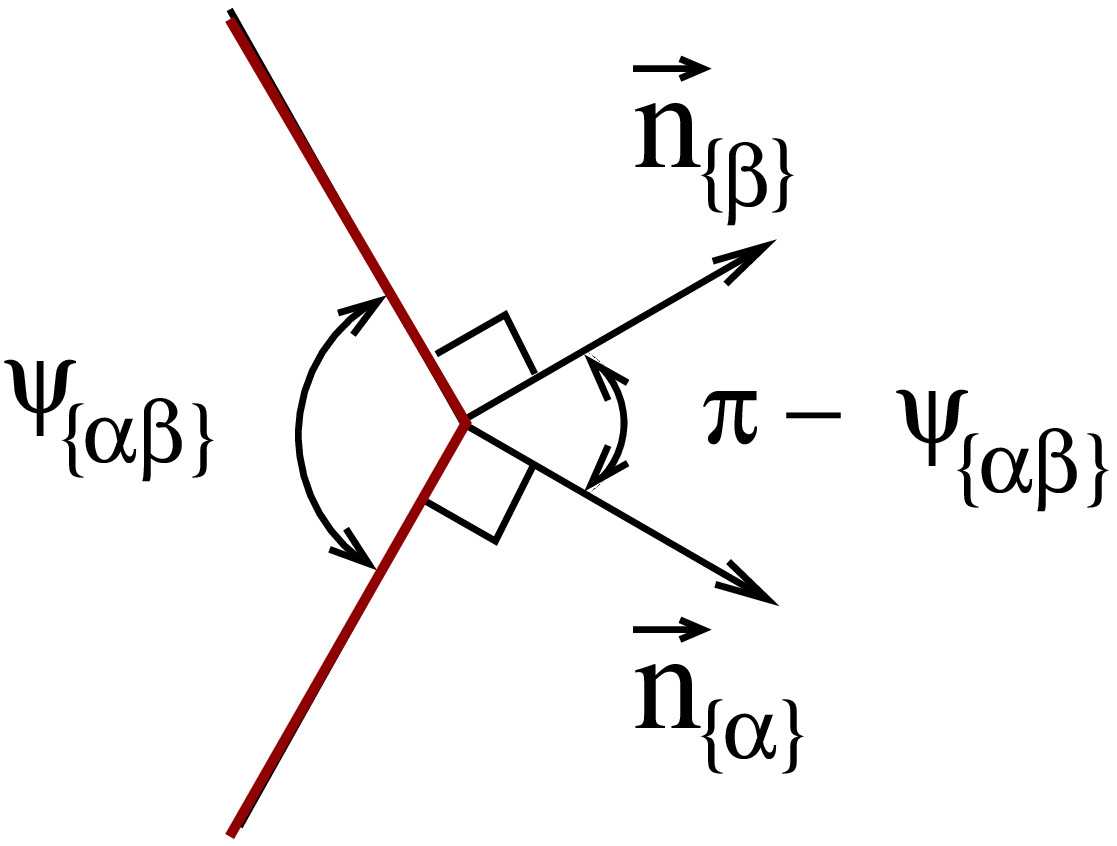}
           \label{f:AnglesB}}

\caption{Left illustration, \ref{f:AnglesA}, shows two neighboring
  triangles in a two-sphere triangulation.  The uniform dihedral angle
  assumption requires $\psi_1=\bar \psi_1$ and $\psi_2=\bar \psi_2$,
  and this in turn implies $\psi_3=\bar \psi_3$ for every pair of
  neighboring triangles.  Right illustration, \ref{f:AnglesB}, shows
  the relationship between the dihedral angle $\psi_{\{\alpha\beta\}}$
  between two faces of a cubic region, and the angle between the
  outward directed normals to the cube faces: $\vec
  n_{\{\alpha\}}\cdot\vec n_{\{\beta\}}=
  \cos(\pi-\psi_{\{\alpha\beta\}}) = - \cos(\psi_{\{\alpha\beta\}})$.
  \label{f:Angles}
}
\end{figure}

The uniform dihedral angle assumption also implies that the
triangulations of the two-sphere at the center of a star-shaped
domain, $\mathcal{S}_I$, must have a special local reflection
symmetry.  Figure~\ref{f:AnglesA} illustrates two neighboring
triangles in one of these triangulations.  If the uniform dihedral
angle assumption has been imposed then $\psi_1=\bar\psi_1=2\pi/K_1$
and $\psi_2=\bar\psi_2=2\pi/K_2$.  The spherical geometry analog of
the angle-side-angle congruence theorem from Euclidean geometry then
implies that $\psi_3=\bar\psi_3$.  With the uniform dihedral angle
assumption, this means that $K_3 = \bar K_{\bar 3}$, i.e. the number
of edges that meet at vertex $3$ in this triangulation must be the
same as the number that meet at vertex, $\bar 3$, of the neighboring
triangle.  This symmetry must apply to every edge of every triangle in
the triangulations of the two-spheres at the centers of each
star-shaped domain $\mathcal{S}_I$.  Therefore, this simple assumption
is quite limiting, and is not satisfied by most two-sphere
triangulations and consequently most multicube structures.

Any open subset of $\Rn$ inherits the flat Euclidean metric of $\Rn$.
Thus any $\mathcal{S}_I$ constructed from parallelograms whose
dihedral angles satisfy the uniform dihedral angle assumption will
naturally inherit a flat metric.  The illustration in
Fig.~\ref{f:AnglesB} shows the relationship between the dihedral angle
$\psi_{\{\alpha\beta\}}$ and the angle between unit normals to the
cube faces, $\vec n_{A\{\alpha\}}= c_{\{\alpha\}}\vec \nabla x$ for
$\alpha=\pm x$, $\vec n_{A\{\alpha\}}= c_{\{\alpha\}}\vec \nabla y$
for $\alpha=\pm y$, and $\vec n_{A\{\alpha\}}= c_{\{\alpha\}}\vec
\nabla z$ for $\alpha=\pm z$.  The constants $c_{\{\alpha\}}$ are
chosen to ensure that the $\vec n_{A\{\alpha\}}$ are the outgoing unit
normals.  The inner products of the outgoing unit normals are related
to the dihedral angles by $\vec n_{A\{\alpha\}}\cdot \vec
n_{A\{\beta\}} = \cos\bigl(\pi-\psi_{A\{\alpha\beta\}}\bigr) =
-\cos{\psi_{A\{\alpha\beta\}}}$.  The inner products of the coordinate
gradients also determine the coordinate components of the inverse
metric: $e^{\,ab}= \vec \nabla x^a \cdot \vec\nabla x^b$.  Therefore
the flat inverse metric, $e^{\,ab}_{A\{\alpha\beta\gamma\}}$,
associated with the vertex $A\{\alpha\beta\gamma\}$, expressed in
terms of the local Cartesian coordinates of region $\mathcal{B}_A$, is
given by
\begin{eqnarray}
  \partial s^{2}_{A\{\alpha\beta\gamma\}}
  &=&e_{A\{\alpha\beta\gamma\}}^{\,ab}\,{\partial_a}\,
  \partial_b,\nonumber\label{e:GlobalFlatMetricS1}
\\
&=&
  \partial_x^{\,2}  +\partial_y^{\,2} +\partial_z^{\,2}
  -2c_{\{\alpha\}}c_{\{\beta\}}\cos{\psi_{A\{\alpha\beta\}}}\,\partial_x\,\partial_y
  \nonumber\\
  &&
  -2c_{\{\alpha\}}c_{\{\gamma\}}\cos{\psi_{A\{\alpha\gamma\}}}\,\partial_x\,\partial_z
  -2c_{\{\beta\}}c_{\{\gamma\}}\cos{\psi_{A\{\beta\gamma\}}}\,\partial_y\,\partial_z,
  \label{e:GlobalFlatMetricS2}
\end{eqnarray} 
where the constants $c_{\{+x\}}=-c_{\{-x\}}=
c_{\{+y\}}=-c_{\{-y\}}=c_{\{+z\}}=-c_{\{-z\}}=1$ ensure the unit
normals are outgoing.  In Eq.~(\ref{e:GlobalFlatMetricS2}) $\alpha=\pm
x$, $\beta=\pm y$ and $\gamma=\pm z$.  These metrics have the correct
dihedral angles between coordinate faces to allow them to fit together
smoothly with the metrics in neighboring regions.  Since the
derivatives of these metrics vanish throughout each region, they are
all flat.

The final step in constructing a flat metric on the $\mathcal{S}_I$
domain is to show that the intrinsic parts of the metrics constructed
in Eq.~(\ref{e:GlobalFlatMetricS2}) are continuous across the
interface boundaries between the cubic regions in $\mathcal{S}_I$.
Equation~(\ref{e:GlobalFlatMetricS2}) shows that these metrics depend
only on the dihedral angles of the edges of the cubic region.  The
simple uniform dihedral angle assumption adopted for this study
implies the local reflection symmetry of the triangulations described
above.  This symmetry guarantees that the dihedral angles of each
cubic region $\mathcal{B}_A$ are the same as those of the neighboring
cubic regions.  The metrics in two neighboring regions will therefore
be related to one another by the local reflection symmetry across the
interface boundary between them.  It follows that the intrinsic parts
of the metrics must be continuous across the interface boundary.  In
general the gauge components of the metric will not be continuous when
expressed in the Cartesian coordinates of the multicube structure.

\subsection{Step 3: Constructing $\hat g_{ab}$.}
\label{s:Step3}

The next step in our procedure for constructing a reference metric is
to build a partition of unity that can be used to combine the flat
metrics from the various overlapping domains into a global
non-singular metric that is smooth within each cubic region, and whose
intrinsic parts are continuous across the interfaces with each
neighboring region.  The needed partition of unity function
$u_{A\{\alpha\beta\gamma\}}(\vec x)\geq 0$ has the value $1$ at the
$A\{\alpha\beta\gamma\}$ vertex of domain $\mathcal{B}_A$, and falls
smoothly to zero on the faces of $\mathcal{B}_A$ that do not intersect
this vertex. The $u_{A\{\alpha\beta\gamma\}}(\vec x)$ are positive
within the star-shaped domain $\mathcal{S}_I$ centered on the
$A\{\alpha\beta\gamma\}$ vertex, and vanish on its outer boundary.
They are used as weight functions to compute averages of the flat
inverse metrics $e_{A\{\alpha\beta\gamma\}}^{\,ab}$ defined on the
$\mathcal{S}_I$ domains in Eq.~(\ref{e:GlobalFlatMetricS2}).  The
inverse of the resulting average, $\hat g_{ab}$, is the global
$C^{\,0}$ reference metric.

First introduce a set of non-negative weight functions,
$w_{A\{\alpha\beta\gamma\}}(\vec x)\geq 0$, whose support is centered
on the vertex $A\{\alpha\beta\gamma\}$.  In the two-dimensional
case~\cite{Lindblom2015} simple separable functions of the global
Cartesian coordinates were used successfully for these weight
functions.  The three-dimensional analogs of those two-dimensional
functions are
\begin{eqnarray}
w_{A{\{\alpha\beta\gamma\}}}(\vec x)&=& 
h\left(\frac{{\Delta x}_{A{\{\alpha\beta\gamma\}}}}{L}\right)
h\left(\frac{{\Delta y}_{A{\{\alpha\beta\gamma\}}}}{L}\right)
h\left(\frac{{\Delta z}_{A{\{\alpha\beta\gamma\}}}}{L}\right),
\label{e:PartitionOfUnity}
\end{eqnarray} 
where the index $A{\{\alpha\beta\gamma\}}$ refers to the vertex of the
cubic region $\mathcal{B}_A$, and $L$ is the coordinate size of the
regions.  The vectors $\vec {\Delta
  x}_{A{\{\alpha\beta\gamma\}}}=\left(\Delta
x_{A\{\alpha\beta\gamma\}}, \Delta y_{A\{\alpha\beta\gamma\}},\Delta
z_{A\{\alpha\beta\gamma\}}\right)$ are defined by
\begin{eqnarray}
  \vec{\Delta x}_{A{\{\alpha\beta\gamma\}}}
  =\vec x - \vec c_A - \vec \nu_{\{\alpha\beta\gamma\}},
    \label{e:DeltaXVertexDef}
\end{eqnarray}
where $\vec x$\, are the global Cartesian coordinates of the multicube
structure that are aligned with the cube faces, and where $\vec
c_A+\vec\nu_{\{\alpha\beta\gamma\}}$ are the coordinates of the vertex
$A{\{\alpha\beta\gamma\}}$. The values of $\vec c_A$, the locations of
the centers of regions $\mathcal{B}_A$, are specified as part of the
definition of the multicube structure; and the values of $\vec
\nu_{\{\alpha\beta\gamma\}}$, the positions of the vertices relative
to the centers are given in Table~\ref{t:TableII}.  They are the same
for all the regions.
\begin{table}[!htb]
  \centering
  \caption{ This table gives the coordinates of each of the eight cube
    vertices $\vec \nu_{\{\alpha\beta\gamma\}}$ with respect to the
    center of ${\cal B}_A$.
    \label{t:TableII} }
  \setlength\tabcolsep{12pt} 
  \begin{tabular}{@{}cc|cc@{}} \toprule
    $\{\alpha\beta\gamma\}$ 
    & $\vec v_{\{\alpha\beta\gamma\}}$
    & $\{\alpha\beta\gamma\}$ 
    & $\vec v_{\{\alpha\beta\gamma\}}$ \\ \midrule
    $\{-x-y-z\}$  & $\tfrac{1}{2}L(-1,-1,-1)$
    & $\{-x-y+z\}$  & $\tfrac{1}{2}L(-1,-1,+1)$ \\ [2pt]
    $\{-x+y-z\}$  & $\tfrac{1}{2}L(-1,+1,-1)$
    & $\{-x+y+z\}$  & $\tfrac{1}{2}L(-1,+1,+1)$ \\ [2pt]
    $\{+x-y-z\}$  & $\tfrac{1}{2}L(+1,-1,-1)$
    & $\{+x-y+z\}$  & $\tfrac{1}{2}L(+1,-1,+1)$ \\ [2pt]
    $\{+x+y-z\}$  & $\tfrac{1}{2}L(+1,+1,-1)$
    & $\{+x+y+z\}$  & $\tfrac{1}{2}L(+1,+1,+1)$ \\ [2pt]
    \bottomrule
  \end{tabular}
\end{table}

\begin{figure}[!ht]
\centering
\includegraphics[width=0.32\textwidth]{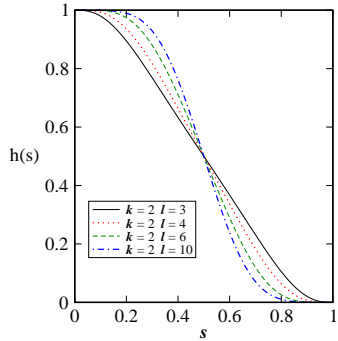}
\caption{Examples of the
  $h(s)=\tfrac{1}{2}\{1+(1-s^{2k})^\ell-[1-(1-|s|)^{2k}]^\ell\}$
  functions defined in Eq.~(\ref{e:PartitionOfUnityhAlt}) used to
  construct the partitions of unity for the $C^0$ reference metrics,
  and in the construction of the $C^0$ and $C^1$ reference metrics.
\label{f:PartitionOfUnityAlt} }
\end{figure}

The functions $h(s)$ used in Eq.~(\ref{e:PartitionOfUnity}) are chosen
to have the values $h(0)=1$ and $h(\pm 1)=0$.  With the arguments
specified in Eq.~(\ref{e:PartitionOfUnity}) this corresponds to
setting $w_{A{\{\alpha\beta\gamma\}}}=1$ at the vertex point
$A\{\alpha\beta\gamma\}$, and $w_{A{\{\alpha\beta\gamma\}}}=0$ on the
boundary faces of $\mathcal{B}_A$ that do not intersect this
vertex. Each of the functions $w_{A{\{\alpha\beta\gamma\}}}(\vec x)$
is also continuous across the $A\{\alpha\}$, $A\{\beta\}$, and
$A\{\gamma\}$ interfaces with the corresponding functions in the
neighboring domains centered on this same vertex. We find that the
functions
\begin{eqnarray}
  h(s)=h(-s)=\tfrac{1}{2}\left\{1+\left(1-s^{2k}\right)^\ell
  -\left[1-(1-|s|)^{2k}\right]^\ell\right\},
\label{e:PartitionOfUnityhAlt}
\end{eqnarray}
with integers $k>0$ and $\ell>0$, work quite well in practice.  Some
of these functions are illustrated in
Fig.~\ref{f:PartitionOfUnityAlt}, with integer values in the range
that worked best in our numerical tests.

The final task in constructing special partitions of unity for the
region $\mathcal{B}_A$ is to construct the normalizing functions
$W_A(\vec x)$:
\begin{eqnarray}
  W_A(\vec x)=\sum_{\{\alpha\beta\gamma\}} w_{A{\{\alpha\beta\gamma\}}}(\vec x),
\end{eqnarray}
where $w_{A{\{\alpha\beta\gamma\}}}(\vec x)$ is defined in
Eq.~(\ref{e:PartitionOfUnity}).  These $W_A(\vec x)$ are strictly
positive, so they can be used to define the partition of unity
functions:
\begin{eqnarray}
  u_{A{\{\alpha\beta\gamma\}}}(\vec x)&=&\frac{w_{A{\{\alpha\beta\gamma\}}}(\vec x)}
  {W_A(\vec x)}.\label{e:PartitionS}
\end{eqnarray}
This normalization ensures that these functions satisfy the
inequalities $0\leq u_{A{\{\alpha\beta\gamma\}}}(\vec x)\leq 1$, and
also the usual partition of unity normalization condition
\begin{eqnarray}
  1=\sum_{\{\alpha\beta\gamma\}} u_{A{\{\alpha\beta\gamma\}}}(\vec x),
\end{eqnarray}
for each $\vec x$\, in each $\mathcal{B}_A$.

A global reference inverse metric $\hat g^{\,ab}(\vec x)$ for $\vec x$\,
in region $\mathcal{B}_A$ can now be constructed by combining the flat
metrics $e_{A{\{\alpha\beta\gamma\}}}^{\,ab}$ defined in
Eq.~(\ref{e:GlobalFlatMetricS2}) with the partition of unity functions
defined in that region by Eq.~(\ref{e:PartitionS}):
\begin{eqnarray}
\hat g^{\,ab}({\vec x})= \sum_{\{\alpha\beta\gamma\}}
u_{A{\{\alpha\beta\gamma\}}}({\vec x})\,
e_{A{\{\alpha\beta\gamma\}}}^{\,ab}({\vec x}).
\label{e:CompleteMetric}
\end{eqnarray}
The sum is over the eight vertices of region $\mathcal{B}_A$.  This
inverse metric is positive definite since it is a linear combination
of positive definite inverse metrics,
$e_{A{\{\alpha\beta\gamma\}}}^{\,ab}$, using the non-negative weight
functions $u_{A{\{\alpha\beta\gamma\}}}(\vec x)$.  A global continuous
reference metric, $\hat g_{ab}(\vec x)$, is then obtained by inverting
$\hat g^{\,ab}(\vec x)$ at each point $\vec x$. The metric $\hat
g_{ab}$ has continuous intrinsic parts across all of the multicube
interface boundaries because it is constructed from flat metrics and
partition of unity functions that are each appropriately continuous
across those interfaces.
  

\section{Constructing a $C^{\,1}$ Three-Dimensional Reference Metric}
\label{s:C1ReferenceMetrics}

In this section the $C^0$ metric, $\hat g_{ab}$, constructed in
Sec.~\ref{s:C0ReferenceMetrics} is transformed into a $C^1$ metric in
three steps.  In the first of these steps, in Sec.~\ref{s:Step3.1}, a
conformal transformation is applied to $\hat g_{ab}$ producing a new
metric, $\bar g_{ab}$, that makes all the edges of each cubic region
into geodesics while keeping the intrinsic parts of $\bar g_{ab}$
continuous across the interface boundaries. This transformation also
fixes one component of the associated extrinsic curvatures along the
edges: $\bar K^{\{\alpha\}}_{ab}\gamma^a\gamma^b=0$ where $\gamma^a$
is tangent to the edge. (See Sec.~\ref{s:Step3.2} for details.) In the
second step, in Sec.~\ref{s:Step3.2}, the metric $\bar g_{ab}$ is
transformed to produce a new metric, $\bar{\bar g}_{ab}$, whose
intrinsic components on each cube face are identical to those of $\bar
g_{ab}$, but whose extrinsic curvatures, $\bar{\bar
  K}^{\{\alpha\}}_{ab}$, vanish identically on each edge of each cubic
region.  In the third step, in Sec.~\ref{s:Step3.3}, the metric
$\bar{\bar g}_{ab}$ is adjusted in the interiors of each cubic region
(keeping its boundary values fixed) in such a way that the resulting
metric $\tilde g_{ab}$ has extrinsic curvatures $\tilde
K^{\{\alpha\}}_{ab}$ that vanish identically on each cube face.  The
$\tilde g_{ab}$ metric constructed by these three steps preserves the
continuity of the intrinsic components of the metric $\hat g_{ab}$
across each interface boundary.  This intrinsic metric continuity
together with the continuity of the extrinsic curvatures across the
interface boundaries (which vanish identically on those boundaries in
this case) are the geometric Israel junction
conditions~\cite{Israel1966} needed to ensure that the metric $\tilde
g_{ab}$ is $C^{\,1}$ across those interfaces.

The methods used in this section are quite general.  In particular,
they do not depend on the uniform dihedral angle assumption used in
Sec.~\ref{s:C0ReferenceMetrics} to construct $\hat g_{ab}$.  The only
specialized assumption used here, specifically in
Sec.~\ref{s:Step3.2}, requires that the dihedral angles are constants
along each cube edge.  This additional assumption is satisfied
automatically by the $\hat g_{ab}$ metrics constructed in
Sec.~\ref{s:C0ReferenceMetrics} using the uniform dihedral angle
assumption, but it will not be satisfied in the most general case.
All that is required to avoid conical singularities is the sum of the
dihedral angles (from each of the cubes that intersect at an edge)
equal $2\pi$ at each point along the edge.  So the dihedral angle in
any one cube may (and perhaps will) in the general case vary along an
edge, so long as the global sum constraint is satisfied.  In this
general case, the construction of the $C^1$ metrics described here,
particularly the parts described in Sec.~\ref{s:Step3.2}, will have to
be generalized as well.

\subsection{Step 1: Converting $\hat g_{ab}$ into $\bar g_{ab}$.}
\label{s:Step3.1}

This section constructs a conformal factor, $e^\phi$, that is used to
transform the reference metric $\hat g_{ab}$ constructed in
Sec.~\ref{s:Step3}:
\begin{equation}
  \bar g_{ab}=e^\phi \hat g_{ab}.\label{e:bargDef}
\end{equation}
The geodesic equation for the curve $x^a(s)$ in the $\bar g_{ab}$
metric is given by
\begin{equation}
  \frac{\mathrm{d}^2 x^a}{\mathrm{d}s^2}+\bar \Gamma^a_{bc}
  \frac{\mathrm{d}x^b}{\mathrm{d}s}\frac{\mathrm{d}x^c}{\mathrm{d}s}
  = A(s) \frac{\mathrm{d}x^a}{\mathrm{d}s},\label{e:geodesicEq}
\end{equation}
where $s$ is an arbitrary parameterization of this curve, where $\bar
\Gamma^a_{bc}$ are the Christoffel symbols of the second kind for this
metric, and where $A(s)$ is a parameter dependent function.  The idea
is to choose a conformal factor $e^\phi$ in Eq.~(\ref{e:bargDef})
having two properties: {\sl a)} it makes each edge of each cubic
region into a geodesic of the metric $\bar g_{ab}$, and {\sl b)} it is
continuous across each cube interface.

Consider the cubic region, $\mathcal{B}_A$, whose Cartesian
coordinates are labeled $x^a=\{x^\alpha, x^\beta,x^\gamma\}$, and
consider the $A\{\alpha\beta\}$ edge of this region where the
$A\{\alpha\}$ and $A\{\beta\}$ faces intersect.  This edge is a curve
with tangent vector $\mathrm{d}x^a/\mathrm{d}s = \{0,0,1\}$, where the
parameter $s$ has been chosen to be $s=x^\gamma$.  An equivalent form
of this equation, more convenient for these purposes, is given by
\begin{equation}
  \bar g_{ab}\frac{\mathrm{d}^2 x^b}{\mathrm{d}s^2} +\bar
  \Gamma_{abc}\frac{\mathrm{d}x^b}{\mathrm{d}s}\frac{\mathrm{d}x^c}{\mathrm{d}s}
  = A(s)\, \bar g_{ab}\frac{\mathrm{d}x^b}{\mathrm{d}s},\label{e:geodesicEq1}
\end{equation}
where the $\bar \Gamma_{abc}$ are the Christoffel symbols of the first
kind.  The three components of this equation can be reduced to
\begin{eqnarray}
  \partial_\gamma \bar g_{\alpha\gamma}
  - \tfrac{1}{2}\partial_\alpha \bar g_{\gamma\gamma}
  &=& A(s)\, \bar g_{\alpha\gamma},\label{e:geodesicEqA}\\
  \partial_\gamma \bar g_{\beta\gamma}
  - \tfrac{1}{2}\partial_\beta \bar g_{\gamma\gamma}
  &=& A(s)\, \bar g_{\beta\gamma},\label{e:geodesicEqB}\\
  \tfrac{1}{2}\partial_\gamma \bar g_{\gamma\gamma}
  &=& A(s)\, \bar g_{\gamma\gamma}.\label{e:geodesicEqC}
\end{eqnarray}
As an interesting aside, note that Eq.~(\ref{e:geodesicEqA}) depends
only on the intrinsic metric on face $A\{\beta\}$, and together with
Eq.~(\ref{e:geodesicEqC}) forms the intrinsic geodesic equation on
this face. Similarly Eq.~(\ref{e:geodesicEqB}) depends only on the
intrinsic metric on face $A\{\alpha\}$, and together with
Eq.~(\ref{e:geodesicEqC}) forms the intrinsic geodesic equation on
this face. Thus the curve formed by the intersection of two surfaces
is a geodesic of the full three-dimensional space if and only if it is
a geodesic of the intrinsic geometry of each surface separately.

The idea now is to choose a conformal factor $\phi$ that transforms
$\hat g_{ab}$, using Eq.~(\ref{e:bargDef}), so the resulting $\bar
g_{ab}$ satisfies Eqs.~(\ref{e:geodesicEqA})--(\ref{e:geodesicEqC}) on
the edges of each cubic region.  The intrinsic parts of the resulting
$\bar g_{ab}$ will be continuous across the interfaces between regions
if and only if the conformal factor $\phi$ is continuous across those
interfaces.  First set $\phi=0$ along the edges of the cubic region to
ensure that $\bar g_{ab}$ is continuous there.  In this case
Eqs.~(\ref{e:geodesicEqA})--(\ref{e:geodesicEqC}) can be re-written in
terms of $\hat g_{ab}$ and $\phi$ for points along the
$A\{\alpha\beta\}$ edge:
\begin{eqnarray}
  \partial_\gamma \hat g_{\alpha\gamma}
  - \tfrac{1}{2}\partial_\alpha \hat g_{\gamma\gamma}
  -\tfrac{1}{2} \hat g_{\gamma\gamma}\partial_\alpha \phi
  &=& A(s)\, \hat g_{\alpha\gamma},\label{e:geodesicEqA.1}\\
  \partial_\gamma \hat g_{\beta\gamma}
  - \tfrac{1}{2}\partial_\beta \hat g_{\gamma\gamma}
  -\tfrac{1}{2} \hat g_{\gamma\gamma}\partial_\beta \phi
  &=& A(s)\, \hat g_{\beta\gamma},\label{e:geodesicEqB.1}\\
  \tfrac{1}{2}\partial_\gamma
  \hat g_{\gamma\gamma} &=& A(s)\, \hat g_{\gamma\gamma}.
  \label{e:geodesicEqC.1}
\end{eqnarray}
The terms involving $\partial_\gamma\phi$ all vanish in these
equations because $\phi=0$ along this edge.  These equations place
constraints on $\partial_\alpha\phi$ and $\partial_\beta\phi$.  In
particular, Eq.~(\ref{e:geodesicEqC.1}), determines $A(s)$ in terms of
$\hat g_{\gamma\gamma}$, while Eqs.~(\ref{e:geodesicEqA.1}) and
(\ref{e:geodesicEqB.1}), can be re-written as boundary conditions for
$\partial_a\phi$ along the $A\{\alpha\beta\}$ edge:
\begin{eqnarray}
  \partial_\alpha \, \phi_{\{A\alpha\beta\}}
   &=& \frac{2 \, \partial_\gamma \hat g_{\alpha\gamma}
     - \partial_\alpha \hat g_{\gamma\gamma}}
   {\hat g_{\gamma\gamma}}
         - \frac{\hat g_{\alpha\gamma}\,\partial_\gamma\hat g_{\gamma\gamma}}
        {\bigl(\hat g_{\gamma\gamma}\bigr)^2},\label{e:psi_bc1}\\
  \partial_\beta\, \phi_{\{A\alpha\beta\}} 
   &=& \frac{2 \, \partial_\gamma \hat g_{\beta\gamma}
     - \partial_\beta \hat g_{\gamma\gamma}}
   {\hat g_{\gamma\gamma}}
         -\frac{ \hat g_{\beta\gamma}\,\partial_\gamma\hat g_{\gamma\gamma}}
        {\bigl(\hat g_{\gamma\gamma}\bigr)^2}.\label{e:psi_bc2}
\end{eqnarray}
Note that these expressions imply that the conformal factor $\phi$
will not vanish everywhere on the cube faces unless $0=\hat
\Gamma_{a\gamma\gamma}$ on the edges, in which case those edges would
already be geodesics of the $\hat g_{ab}$ metric.  Also note that the
metric $\hat g_{ab}$ constructed in Sec.~\ref{s:C0ReferenceMetrics}
rapidly approaches a constant flat metric at each vertex of the cubic
region.  It follows that the connection $\hat \Gamma_{abc}$ and the
gradient $\partial_a\phi$ all vanish at these vertex points.

The next step is to extend the conformal factor $\phi$ across the
faces of the cubic region $\mathcal{B}_A$ in a way that {\sl a)}
satisfies the boundary conditions given in Eqs.~(\ref{e:psi_bc1}) and
(\ref{e:psi_bc2}) along each edge, and {\sl b)} ensures that it is
continuous across the interface with the neighboring cubic region.
The conformal factor $\phi_{A\{\alpha\}}$ on the $A\{\alpha\}$ face
satisfies the boundary conditions $\phi_{A\{\alpha\beta\}}=0$ and
Eq.~(\ref{e:psi_bc2}) along the $A\{\alpha\beta\}$ edge. Analogous
conditions must also be imposed on each edge of this face.  Together,
these conditions constitute Dirichlet and Neumann boundary conditions
for $\phi_{A\{\alpha\}}$ on the $A\{\alpha\}$ face.  One convenient
way to find $\phi_{A\{\alpha\}}$ that satisfy these boundary
conditions is to solve the bi-harmonic equation for
$\phi_{A\{\alpha\}}$ on this face:
\begin{eqnarray}
\left(\partial_\beta^{\,4} + 2\, \partial_\beta^{\,2}\,\partial_\gamma^{\,2}
+\partial_\gamma^{\,4}\right)\phi_{A\{\alpha\}}&=&0.
\label{e:Biharmonic2Da0}
\end{eqnarray}
Solutions to the bi-harmonic equation are uniquely determined by
specifying both Dirichlet and Neumann conditions on the boundary of a
compact domain~\cite{Barton2014}.  This approach can
then be used to determine the surface values of $\phi$ on each face of
each region in the multicube structure.  The pseudo-spectral
numerical methods used to solve this equation for this study are
described briefly in \ref{s:BiharmonicMethods}.

The boundary conditions that determine the solution to
Eq.~(\ref{e:Biharmonic2Da0}) only depend on the intrinsic components
of the metric, $\hat g_{\beta\beta}$, $\hat g_{\beta\gamma}$ and $\hat
g_{\gamma\gamma}$, on the $A\{\alpha\}$ face.  These
intrinsic metric components were constructed to be continuous across
this interface boundary in Sec.~\ref{s:C0ReferenceMetrics}.  It
follows that boundary conditions used to determine
$\phi_{A\{\alpha\}}$ will be the same on both sides of the interface
boundary.  Since the solution to Eq.~(\ref{e:Biharmonic2Da0}) with
Dirichlet and Neumann boundary conditions is unique~\cite{Barton2014},
it follows that the $\phi_{A\{\alpha\}}$ determined in this way will be
the same on both sides of the interface.

The method describe above can be used to determine the surface values
$\phi_{A\{\alpha\}}$ on each face of each cubic region. These
solutions provide Dirichlet boundary conditions for the full conformal
factor $\phi_A$ within each region.  The normal derivatives of
$\phi_{A\{\alpha\}}$ are unconstrained, however, beyond the
requirement that those derivatives agree along the edges with the
tangential derivatives from the neighboring faces.  The conformal
factor $\phi$ within the cubic region can therefore be determined in
any number of ways.  For example, it could be determined by solving
the three-dimensional Laplace equation with the Dirichlet boundary
conditions $\phi_{A\{\alpha\}}$.

A computationally more efficient approach has been adopted for this
study.  Begin by defining a set of coordinates $s^\alpha_A$ that
measure the relative distance between a point inside region
$\mathcal{B}_A$ and the $\partial_{-\alpha}\mathcal{B}_A$ face of that
region.  The $s^\alpha_A$ are normalized so that $s^\alpha_A=0$ on the
$\partial_{\alpha}\mathcal{B}_A$ face, while $s^\alpha_A=1$ on the
opposite $\partial_{-\alpha}\mathcal{B}_A$ face.  In particular
\begin{eqnarray}
  s^{\pm x}_A = \left|\frac{x - c_A^x}{L} \mp\frac{1}{2}\right|, \qquad
  s^{\pm y}_A = \left|\frac{y - c_A^y}{L} \mp\frac{1}{2}\right|, \qquad
  s^{\pm z}_A = \left|\frac{z - c_A^z}{L} \mp\frac{1}{2}\right|,
  \label{e:s_alphaDef}
\end{eqnarray}
where $\vec x=(x,y,z)$ are the global Cartesian coordinates of the
multicube structure, $\vec c_A=(c^x_A,c^y_A,c^z_A)$ are the
coordinates of the center, and $L$ is the coordinate size of region
$\mathcal{B}_A$.
 
The conformal factor $\phi_{A\{\alpha\}}$ on the $A\{\alpha\}$ face,
constructed by solving Eq.~(\ref{e:Biharmonic2Da0}), can now be
extrapolated into the interior using the $h(s)$ functions defined in
Eq.~(\ref{e:PartitionOfUnityhAlt}).  Consider the extrapolation
$\phi^\alpha_A=h(s^\alpha_A)\,\phi_{A\{\alpha\}}$.  The
$\phi_{A\{\alpha\}}$ vanish identically along the edges because of the
boundary conditions used to solve Eq.~(\ref{e:Biharmonic2Da0}).
Therefore the $\phi_{A\{\alpha\}}$ extrapolated in this way do not
modify the $\phi_{A\{\beta\}}$ on the adjacent faces. It does not
modify $\phi_{A\{-\alpha\}}$ on the opposite face, either, because
$h(s^\alpha_A)$ vanishes there.  The complete conformal factor
$\phi_A$ in the interior of region $\mathcal{B}_A$ can therefore be
determined by combining the extrapolations from all the cube faces:
\begin{equation}
  \phi_A=\sum_\alpha\phi^\alpha_A=\sum_\alpha
  h(s^\alpha_A)\,\phi_{A\{\alpha\}}.
\end{equation}
The resulting $\phi_A$ automatically satisfies the Dirichlet conditions
$\phi_{A\{\alpha\}}$ on each of the faces.

The conditions in Eqs.~(\ref{e:psi_bc1}) and (\ref{e:psi_bc2}) ensure
that the edges of each cubic region are geodesics of the metric $\bar
g_{ab}=e^\phi\hat g_{ab}$.  The continuity across the interface
boundaries of $\phi_A$ ensures that the global solution for $\phi$ is
continuous across those boundaries.  And this in turn ensures that the
intrinsic components of the $\bar g_{ab}$ metric are continuous across
all the interface boundaries as well.

\subsection{Step 2: Converting $\bar g_{ab}$ into $\bar{\bar g}_{ab}$.}
\label{s:Step3.2}

Let $\bar g_{ab}$ denote the global $C^{\,0}$ metric constructed in
Sec.~\ref{s:Step3.1}.  The goal of this second step is to convert
$\bar g_{ab}$ into a metric $\bar{\bar g}_{ab}$ having two important
properties: first, the extrinsic curvatures $\bar{\bar
  K}^{\{\alpha\}}_{ab}$ associated with $\bar{\bar g}_{ab}$ must
vanish identically on each edge of each cubic region; and second, the
intrinsic parts of $\bar{\bar g}_{ab}$ must be identical to those of
$\bar g_{ab}$.  We note that while $\bar g_{ab}$ was constructed to
have intrinsic parts that are only $C^{\,0}$ across the interface
boundaries, within each region $\bar g_{ab}$ is actually smooth.

Consider the interface boundary $A\{\alpha\}$ of cubic region
$\mathcal{B}_A$.  In the global Cartesian coordinates of our multicube
structure, the boundary $A\{\alpha\}$ is a level surface of the
coordinate $x^\alpha$.  The unit normal co-vector field to the
foliation of constant $x^\alpha$ surfaces, $\bar{\bar
  n}^{\{\alpha\}}_a$, is given by
\begin{eqnarray}
 \bar{\bar n}^{\{\alpha\}}_a=\bar{\bar N}^{\{\alpha\}}\partial_a
 x^\alpha,\label{e:nbarbarDef}
\end{eqnarray}
where $\bar{\bar N}^{\{\alpha\}}=\epsilon^{\{\alpha\}}\left( \bar{\bar
  g}^{\alpha\alpha}\right)^{-1/2}$.  The constant
$\epsilon^{\{\alpha\}}=\pm 1$ determines the sign of $\bar{\bar
  N}^{\{\alpha\}}$, and is chosen to ensure that $\bar{\bar
  n}^{\{\alpha\}}_a$ is the outward directed unit normal on the
$A\{\alpha\}$ face.  The extrinsic curvature $\bar{\bar
  K}^{\{\alpha\}}_{ab}$ of this surface is given by
\begin{eqnarray}
  \bar{\bar K}^{\{\alpha\}}_{ab}
  &=&\bar{\bar P}^{\{\alpha\}\, c}_a\bar{\bar
  P}^{\{\alpha\}\, d}_b\, \bar{\bar \nabla}_c\,\bar{\bar
  n}^{\{\alpha\}}_d\label{e:KbarbarDefA} \\
  &=&\bar{\bar P}^{\{\alpha\}\, c}_a\bar{\bar P}^{\{\alpha\}\,
    d}_b \left(\partial_c\,\bar{\bar
  n}^{\{\alpha\}}_d-\bar{\bar n}^{\{\alpha\}}_e\,\bar{\bar \Gamma}^e_{cd}\right)
  \label{e:KbarbarDefB} \\
  &=&\tfrac{1}{2} \bar{\bar P}^{\{\alpha\}\, c}_a\bar{\bar
  P}^{\{\alpha\}\, d}_b \,\bar{\bar n}^{\{\alpha\}\,e}\left(\partial_e\,\bar{\bar
  g}_{cd}-\partial_c\,\bar{\bar g}_{ed} -\partial_d\,\bar{\bar
  g}_{ec}\right),\label{e:KbarbarDefC}
\end{eqnarray}
where $\bar{\bar \nabla}_c$ is the $\bar{\bar g}_{ab}$
metric-compatible covariant derivative, and $\bar{\bar
  P}^{\{\alpha\}\,c}_{a}$ is the projection tensor $\bar{\bar
  P}^{\{\alpha\}\,c}_{a}=\delta_{a}^c-\bar{\bar
  n}^{\{\alpha\}}_a\bar{\bar n}^{\{\alpha\}\,c}$.  Note that the term
proportional to $\partial_c \bar{\bar
  n}^{\{\alpha\}}_d=\partial_c\bar{\bar
  N}^{\{\alpha\}}\partial_dx^\alpha$ in Eq.~(\ref{e:KbarbarDefB})
vanishes identically because $\bar{\bar
  P}^{\{\alpha\}\,d}_b\partial_dx^\alpha=0$.

We define the difference between the $\bar{\bar g}_{ab}$ and the $\bar
g_{ab}$ metrics, $\delta \bar g_{ab}$, and the associated differences
between extrinsic curvatures, $\delta \bar K^{\{\alpha\}}_{ab}$:
\begin{eqnarray}
  \delta \bar g_{ab}&=& \bar{\bar g}_{ab}-\bar g_{ab},
    \label{e:DeltagbarDef}\\
    \delta \bar
  K^{\{\alpha\}}_{ab}&=&\bar{\bar K}^{\{\alpha\}}_{ab} -\bar
  K^{\{\alpha\}}_{ab}.\label{e:DeltaKbarDef}
\end{eqnarray}
Note that these differences are not necessarily infinitesimal.  To
ensure that the intrinsic parts of $\bar{ \bar g}_{ab}$ are identical
to those of $\bar g_{ab}$, we choose $\delta \bar g_{ab}$ to be a
smooth tensor in the interior of each cubic region that satisfies,
\begin{eqnarray}
  \bar{ \bar P}_a^{\{\alpha\}\,c} \bar{\bar P}_b^{\{\alpha\}\,d}\delta
  \bar g_{cd}=0
\end{eqnarray}
on each cube face $A\{\alpha\}$.  Note that the projection tensor
$\bar{\bar P}_a^{\{\alpha\}\,c}$ is identical to $\bar
P_a^{\{\alpha\}\,c}$ on $A\{\alpha\}$ because
$\bar{\bar P}_a^{\{\alpha\}\,c}\partial_cx^\alpha =\bar
P_a^{\{\alpha\}\,c}\partial_cx^\alpha=0$.  Therefore the metric
continuity condition on $\delta\bar g_{ab}$ is equivalent to
\begin{eqnarray}
  { \bar P}_a^{\{\alpha\}\,c} {\bar P}_b^{\{\alpha\}\,d}\delta \bar
  g_{cd}=0,
  \label{e:IntrinsicMetricContinuity}
\end{eqnarray}
which is easier to enforce since the metric $\bar g_{ab}$, and
consequently the normal vector $\bar n^{\{\alpha\}}_a$, is already
known.  The condition that the extrinsic curvature $\bar{\bar
  K}^{\{\alpha\}}_{ab}$ vanish along each edge of the
$A\{\alpha\}$ face can be expressed as the following
condition on $\delta\bar K^{\{\alpha\}}_{ab}$,
\begin{eqnarray}
    \delta \bar
  K^{\{\alpha\}}_{ab}&=&-\bar K^{\{\alpha\}}_{ab},\label{e:DeltaKbarCondition}
\end{eqnarray}
on each edge of each cube face $A\{\alpha\}$.

To determine exactly what restrictions are placed on $\delta \bar
g_{ab}$ by the intrinsic metric and extrinsic curvature continuity
conditions,
Eqs.~(\ref{e:IntrinsicMetricContinuity})--(\ref{e:DeltaKbarCondition}),
we examine those conditions expressed in the Cartesian coordinates
$\{x^\alpha,x^\beta,x^\gamma\}$ of region $\mathcal{B}_A$.  The
$A\{\alpha\}$ face of this region is defined by an $x^\alpha$=constant
surface, while the $x^\beta$ and $x^\gamma$ coordinates label points
on that face.  In these coordinates the intrinsic metric continuity
condition, Eq.~(\ref{e:IntrinsicMetricContinuity}), implies that all
the $x^\beta$ and $x^\gamma$ components of the metric perturbation
vanish everywhere on that face: $\delta \bar g_{\beta\beta}=\delta
\bar g_{\beta\gamma}=\delta \bar g_{\gamma\gamma}=0$.  Similarly on
the adjacent $A\{\beta\}$ face, all the $x^\alpha$ and $x^\gamma$
components of the metric perturbation vanish 
$\delta \bar g_{\alpha\alpha}=\delta \bar g_{\alpha\gamma}=\delta \bar
g_{\gamma\gamma}=0$.  It follows that all the components of $\delta
\bar g_{ab}$ except $\delta \bar g_{\alpha\beta}$ must vanish along
the $A\{\alpha\beta\}$ edge: $\delta \bar g_{\alpha\alpha}=\delta
\bar g_{\alpha\gamma}= \delta \bar g_{\beta\beta}=\delta \bar
g_{\beta\gamma}=\delta \bar g_{\gamma\gamma}=0$.

The $C^{\,0}$ metric $\bar g_{ab}$ was constructed so that the
dihedral angle between the $A\{\alpha\}$ and $A\{\beta\}$ faces, $\cos
\psi_{\{A\alpha\beta\}} = -\bar g^{\alpha\beta}/\sqrt{\bar
  g^{\alpha\alpha}\bar g^{\beta\beta}}$, is constant along the edge
between these faces.  This was done to ensure there is no conical
singularity along this edge.  To ensure that the $\bar{\bar g}_{ab}$
metric has no conical singularity there, we keep this dihedral angle
fixed in this metric along this edge as well.  This can be done by
imposing the additional constraint $\delta \bar g_{\alpha\beta}=0$
along this edge.  This makes all the components $\delta \bar
g_{ab}=0$, and consequently $\delta\bar g^{ab}=0$ which keeps the
dihedral angles fixed.  Thus the intrinsic metric continuity
conditions, along with the conditions to ensure there are no conical
singularities along the edges, require that all the components of the
metric perturbation vanish along each edge: $\delta \bar g_{ab}=0$.

Exact expressions for the $x^\beta$ and $x^\gamma$ components of the
extrinsic curvature $\bar{\bar K}^{\{\alpha\}}_{ab}$ on the
$A\{\alpha\}$ face are obtained from
Eqs.~(\ref{e:KbarbarDefC}) and (\ref{e:DeltagbarDef}):
\begin{eqnarray}
  \bar{\bar K}^{\{\alpha\}}_{\beta\beta} &=&
  \tfrac{1}{2} \bar{\bar n}^{\{\alpha\}\,a}
  \left(\partial_a\,\bar g_{\beta\beta}-2\partial_\beta\,\bar g_{a\beta}\right)+
  \tfrac{1}{2} \bar{\bar n}^{\{\alpha\}\,a}
  \left(\partial_a\,\delta \bar g_{\beta\beta}
  -2\partial_\beta\,\delta\bar g_{a\beta} \right).\label{e:DeltaKMuNuDef1}\\
  \bar{\bar K}^{\{\alpha\}}_{\beta\gamma} &=&
  \tfrac{1}{2} \bar{\bar n}^{\{\alpha\}\,a}
  \left(\partial_a\,\bar g_{\beta\gamma}
  -\partial_\beta\,\bar g_{a\gamma} -\partial_\gamma\,\bar g_{a\beta}\right)+
  \tfrac{1}{2} \bar{\bar n}^{\{\alpha\}\,a}
  \left(\partial_a\,\delta \bar g_{\beta\gamma}
  -\partial_\beta\,\delta \bar g_{a\gamma}
  -\partial_\gamma\,\delta \bar g_{a\beta} \right).\label{e:DeltaKMuNuDef2}\\
  \bar{\bar K}^{\{\alpha\}}_{\gamma\gamma} &=&
  \tfrac{1}{2} \bar{\bar n}^{\{\alpha\}\,a}
  \left(\partial_a\,\bar g_{\gamma\gamma}-2\partial_\gamma\,\bar g_{a\gamma} \right)+
  \tfrac{1}{2} \bar{\bar n}^{\{\alpha\}\,a}
  \left(\partial_a\,\delta \bar g_{\gamma\gamma}
  -2\partial_\gamma\,\delta \bar g_{a\gamma} \right).\label{e:DeltaKMuNuDef3}
\end{eqnarray}
On the $A\{\alpha\beta\}$ edge, where the $A\{\alpha\}$
and $A\{\beta\}$ faces intersect, $\delta \bar
g_{ab}=0$, so $\bar{\bar n}^{\{\alpha\}\,a}=\bar n^{\{\alpha\}\,a}$.
Consequently Eqs.~(\ref{e:DeltaKMuNuDef1})--(\ref{e:DeltaKMuNuDef3})
can be re-written in the simpler form:
\begin{eqnarray}
  \delta \bar K^{\{\alpha\}}_{\beta\beta}  &=& 
  \bar{\bar K}^{\{\alpha\}}_{\beta\beta}-
  \bar K^{\{\alpha\}}_{\beta\beta} =  \tfrac{1}{2} \bar
  n^{\{\alpha\}\,a}\left(\partial_a\,\delta \bar g_{\beta\beta}-
  2\partial_\beta\,\delta \bar g_{a\beta} \right),\label{e:DeltaKExp1}\\
  \delta \bar K^{\{\alpha\}}_{\beta\gamma}  &=&
  \bar{\bar K}^{\{\alpha\}}_{\beta\gamma}
  - \bar K^{\{\alpha\}}_{\beta\gamma} =  \tfrac{1}{2} \bar
  n^{\{\alpha\}\,a}\left(\partial_a\,\delta \bar g_{\beta\gamma}
  -\partial_\beta\,\delta \bar g_{a\gamma}
  -\partial_\gamma\,\delta \bar g_{a\beta}\right),\label{e:DeltaKExp2}\\
  \delta \bar K^{\{\alpha\}}_{\gamma\gamma}  &=&
  \bar{\bar K}^{\{\alpha\}}_{\gamma\gamma}
  -\bar K^{\{\alpha\}}_{\gamma\gamma} =  \tfrac{1}{2} \bar
  n^{\{\alpha\}\,a}\left(\partial_a\,\delta \bar g_{\gamma\gamma}
  -2\partial_\gamma\,\delta \bar g_{a\gamma} \right).\label{e:DeltaKExp3}
\end{eqnarray}

Since $\delta \bar g_{ab}=0$ along the $A\{\alpha\beta\}$ edge, it
follows that $\partial_\gamma\delta \bar g_{ab}=0$ there.  Since
$\delta \bar g_{\beta\beta}= \delta \bar g_{\beta\gamma}=\delta \bar
g_{\gamma\gamma}=0$ everywhere on the $A\{\alpha\}$ face, it follows
that $\partial_\beta\delta \bar g_{\beta\beta}=\partial_\beta\delta
\bar g_{\beta\gamma} =\partial_\beta\delta \bar g_{\gamma\gamma}=0$
along the $A\{\alpha\beta\}$ edge as well.  Finally, $\delta \bar
g_{\gamma\gamma}$ on the adjacent $\partial_\beta \mathcal{B}_A$ face
represents a perturbation of the intrinsic metric, so
$\partial_\alpha\delta \bar g_{\gamma\gamma}=0$ along the
$A\{\alpha\beta\}$ edge as well.  The components of $\bar
n^{\{\alpha\}\,a}$ in these coordinates are given by $\bar
n^{\{\alpha\}\,a}=\bar N^{\{\alpha\}} \{\bar g^{\alpha\alpha},\bar
g^{\alpha\beta},\bar g^{\alpha\gamma}\}$, so the expressions for
$\delta \bar K^{\{\alpha\}}_{ab}$ from
Eqs.~(\ref{e:DeltaKExp1})--(\ref{e:DeltaKExp3}) can be simplified
further:
\begin{eqnarray}
  \delta \bar K^{\{\alpha\}}_{\beta\beta}&=&  \tfrac{1}{2}\, \bar N^{\{\alpha\}}\,
  \bar g^{\alpha\alpha}\left(\partial_\alpha\,\delta \bar g_{\beta\beta}
  -2\,\partial_\beta\,\delta \bar g_{\alpha\beta}
  \right),\label{e:dKalpha1}\\ 
  \delta \bar K^{\{\alpha\}}_{\beta\gamma}&=&  \tfrac{1}{2}\, \bar N^{\{\alpha\}}\,
  \bar g^{\alpha\alpha} \left(\partial_\alpha\,\delta \bar g_{\beta\gamma}
  -\partial_\beta\,\delta \bar g_{\alpha\gamma}
  \right),\label{e:dKalpha2}\\
  \delta \bar K^{\{\alpha\}}_{\gamma\gamma}&=&  0.\label{e:dKalpha3}
\end{eqnarray}
The analogous expressions for $\delta\bar K^{\{\beta\}}_{ab}$, the
extrinsic curvature of the adjacent $A\{\beta\}$
face, along this edge can be obtained by interchanging the roles of
$x^\alpha$ and $x^\beta$ in
Eqs.~(\ref{e:dKalpha1})--(\ref{e:dKalpha3}):
\begin{eqnarray}
  \delta \bar K^{\{\beta\}}_{\alpha\alpha}&=&  \tfrac{1}{2}\, \bar N^{\{\beta\}}\,
  \bar g^{\beta\beta}\left(\partial_\beta\,\delta \bar g_{\alpha\alpha}
  -2\,\partial_\alpha\,\delta \bar g_{\alpha\beta}
  \right),\label{e:dKmu1}\\ 
  \delta \bar K^{\{\beta\}}_{\alpha\gamma}&=&  \tfrac{1}{2}\, \bar N^{\{\beta\}}\,
  \bar g^{\beta\beta} \left(\partial_\beta\,\delta \bar g_{\alpha\gamma}
  -\partial_\alpha\,\delta \bar g_{\beta\gamma}
  \right),\label{e:dKmu2}\\
  \delta \bar K^{\{\beta\}}_{\gamma\gamma}&=&  0.\label{e:dKmu3}
\end{eqnarray}
These expressions for $\delta \bar K^{\{\alpha\}}_{\beta\gamma}$ and
$\delta \bar K^{\{\beta\}}_{\alpha\gamma}$ define the required
boundary conditions on the derivatives $\partial_a\delta \bar g_{ab}$
along the $A\{\alpha\beta\}$ edge, where the $A\{\alpha\}$ and
$A\{\beta\}$ faces intersect.  Since $\delta \bar K^{\{\alpha\}}_{ab}$
must satisfy Eq.~(\ref{e:DeltaKbarCondition}), we see that these
boundary conditions imply that the gauge components of the metric
(i.e. the non-intrinsic components) $\delta \bar g_{\alpha\alpha}$,
$\delta \bar g_{\alpha\beta}$, and $\delta \bar g_{\alpha\gamma}$
cannot simply be set to zero on the $A\{\alpha\}$ face.

The expressions for $\delta\bar K^{\{\alpha\}}_{\gamma\gamma}$ and
$\delta \bar K^{\{\beta\}}_{\gamma\gamma}$, Eqs.~(\ref{e:dKalpha3})
and (\ref{e:dKmu3}), along the $A\{\alpha\beta\}$ edge, imply that no
discontinuity in $\bar K^{\{\alpha\}}_{\gamma\gamma}$ or $\bar
K^{\{\beta\}}_{\gamma\gamma}$ along this edge can be removed by any
$\delta \bar g_{ab}$ allowed by our constraints.  To understand what
that means, let $\gamma^a$ be the components of the vector
$\partial_\gamma = \gamma^a\partial_a$.  This vector is orthogonal to
the surface normals: $0=\bar n^{\{\alpha\}}_a\gamma^a=\bar
n^{\{\beta\}}_a\gamma^a$.  It follows that $\bar
K^{\{\alpha\}}_{\gamma\gamma}=\gamma^a\gamma^b\bar \nabla_a \bar
n^{\{\alpha\}}_b =\gamma^a\bar \nabla_a\left(\bar 
n^{\{\alpha\}}_b\gamma^b\right) - \bar n^{\{\alpha\}}_b\gamma^a \bar
\nabla_a \gamma^b=- \bar n^{\{\alpha\}}_b\gamma^a \bar \nabla_a
\gamma^b$.  Since the metric $\bar g_{ab}$ constructed in
Sec.~\ref{s:Step3.1} has the property that each edge is a geodesic,
$\gamma^a\bar\nabla_a\gamma^b=A(x^\gamma)\gamma^b$ (for some function
$A(x^\gamma)$ along this edge), it follows that $\bar
K^{\{\alpha\}}_{\gamma\gamma}=0$ and similarly $\bar
K^{\{\beta\}}_{\gamma\gamma}=0$ along the $A\{\alpha\beta\}$ edge.
This component of the extrinsic curvature continuity condition
Eq.~(\ref{e:DeltaKbarCondition}) is therefore satisfied automatically
along this edge, so $\delta \bar K^{\{\alpha\}}_{\gamma\gamma} =
\delta \bar K^{\{\beta\}}_{\gamma\gamma}=0$ are the appropriate
corrections there.

The right sides of Eqs.~(\ref{e:dKalpha2}) and (\ref{e:dKmu2}) are
both proportional to $\partial_\beta\,\delta \bar g_{\alpha\gamma}
-\partial_\alpha\,\delta \bar g_{\beta\gamma}$, so the extrinsic
curvature perturbations on the left sides must also be related: $\bar
N^{\{\beta\}}\,\delta \bar K^{\{\beta\}}_{\alpha\gamma} =-\bar
N^{\{\alpha\}}\,\delta \bar K^{\{\alpha\}}_{\beta\gamma}$, obtained by
simplifying using $\bar N^{\{ \alpha \}} = \epsilon^{\{\alpha\}}
(g^{\alpha\alpha})^{-1/2}$ and $\bar N^{\{ \beta \}} =
\epsilon^{\{\beta\}} (g^{\beta\beta})^{-1/2}$.  This condition is
inconsistent with Eq.~(\ref{e:DeltaKbarCondition}) unless
\begin{equation}
\bar N^{\{\beta\}}\, \bar K^{\{\beta\}}_{\alpha\gamma} =-\bar
N^{\{\alpha\}}\, \bar K^{\{\alpha\}}_{\beta\gamma}.
\label{e:KbarCoordinateIdentity}
\end{equation}
The simple proof of this identity is given in \ref{s:IdentityProof}.
This identity shows that the edge constraints given in
Eqs.~(\ref{e:dKalpha2}) and (\ref{e:dKmu2}) for $\delta
K^{\{\alpha\}}_{\beta\gamma}$ and $\delta
K^{\{\beta\}}_{\alpha\gamma}$ are self-consistent.
 
Equation~(\ref{e:DeltaKbarCondition}) together with
Eqs.~(\ref{e:dKalpha1})--(\ref{e:dKmu3}) place the following
constraints on the derivatives of certain components of $\delta \bar
g_{ab}$ along the $A\{\alpha\beta\}$ edge,
\begin{eqnarray}
  \partial_\alpha\delta \bar g_{\beta\beta}
  -2\partial_\beta\delta \bar g_{\alpha\beta}&=&
  -2 \bar N^{\{\alpha\}} \bar K^{\{\alpha\}}_{\beta\beta},
  \label{e:GaugeConstraint1}\\
  \partial_\alpha\delta \bar g_{\beta\gamma}
  -\partial_\beta\delta \bar g_{\alpha\gamma}&=&
  - 2 \bar N^{\{\alpha\}} \bar K^{\{\alpha\}}_{\beta\gamma}
  =2 \bar N^{\{\beta\}} \bar K^{\{\beta\}}_{\alpha\gamma},
  \label{e:GaugeConstraint2}\\
  \partial_\beta\delta \bar g_{\alpha\alpha}
  -2\partial_\alpha\delta \bar g_{\alpha\beta}&=&
  -2 \bar N^{\{\beta\}} \bar K^{\{\beta\}}_{\alpha\alpha}.
  \label{e:GaugeConstraint3}
\end{eqnarray}
The metric perturbation components $\delta \bar g_{\alpha\alpha}$,
$\delta \bar g_{\alpha\beta}$, and $\delta \bar g_{\alpha\gamma}$ do
not affect the intrinsic metric on the $A\{\alpha\}$ face, while the
$\delta \bar g_{\beta\beta}$, $\delta \bar g_{\beta\alpha}$ and
$\delta \bar g_{\beta\gamma}$ components do not affect the intrinsic
metric on the $A\{\beta\}$ face.  These components therefore play the
role of gauge degrees of freedom on these faces, which can be chosen
arbitrarily subject to the constraints in
Eqs.~(\ref{e:GaugeConstraint1})--(\ref{e:GaugeConstraint3}).  While
not unique, one self-consistent way to satisfy these constraints along
the $A\{\alpha\beta\}$ edge is given by
\begin{eqnarray}
  \partial_\beta\delta \bar g_{\alpha\alpha} &=&
  -2 \bar N^{\{\beta\}} \bar K^{\{\beta\}}_{\alpha\alpha},
  \label{e:EdgeCondition1}\\
    \partial_\beta\delta \bar g_{\alpha\beta}&=& 0,\label{e:EdgeCondition2}\\
    \partial_\beta \delta \bar g_{\alpha\gamma}&=&-\bar N^{\{\beta\}}
    \bar K^{\{\beta\}}_{\alpha\gamma}.\label{e:EdgeCondition3}\\
  \partial_\alpha\delta \bar g_{\beta\beta} &=&
  -2 \bar N^{\{\alpha\}} \bar K^{\{\alpha\}}_{\beta\beta},
  \label{e:EdgeCondition4}\\
    \partial_\alpha\delta \bar g_{\beta\alpha}&=& 0,\label{e:EdgeCondition5}\\
    \partial_\alpha \delta \bar g_{\beta\gamma}&=&-\bar N^{\{\alpha\}}
    \bar K^{\{\alpha\}}_{\beta\gamma}.\label{e:EdgeCondition6}
  \end{eqnarray}
We note that the equations for $\delta \bar g_{\beta\beta}$, $\delta
\bar g_{\beta\alpha}$, and $\delta \bar g_{\beta\gamma}$ in
Eqs.~(\ref{e:EdgeCondition4})--(\ref{e:EdgeCondition6}) can be
obtained from those for $\delta \bar g_{\alpha\alpha}$, $\delta
\bar g_{\alpha\beta}$, and $\delta \bar g_{\alpha\gamma}$ in
Eqs.~(\ref{e:EdgeCondition1})--(\ref{e:EdgeCondition3}) simply by
exchanging the $\alpha$ and $\beta$ indices.

The intrinsic components of the metric perturbations $\delta \bar
g_{ab}$ must vanish on the $A\{\alpha\}$ face, $0 = \delta \bar
g_{\beta\beta}= \delta \bar g_{\beta\gamma} = \delta \bar
g_{\gamma\gamma}$. It follows from
Eqs.~(\ref{e:EdgeCondition1})--(\ref{e:EdgeCondition3}) that the full
set of boundary conditions on $\delta \bar g_{ab}$ along the
$A\{\alpha\beta\}$ edge of the $A\{\alpha\}$ face are given by
\begin{eqnarray}
  \delta \bar g_{ab}&=&0,\label{e:EdgeBC1}\\
  \partial_\beta\delta \bar g_{\alpha\alpha} &=&
  -2 \bar N^{\{\beta\}} \bar K^{\{\beta\}}_{\alpha\alpha},
  \label{e:EdgeBC2}\\
    \partial_\beta \delta \bar g_{\alpha\gamma}&=&-\bar N^{\{\beta\}}
    \bar K^{\{\beta\}}_{\alpha\gamma},\label{e:EdgeBC3}\\
    \partial_\beta\delta \bar g_{\alpha\beta}&=&
    \partial_\beta\delta \bar g_{\beta\beta}\,\,=\,\,
    \partial_\beta\delta \bar g_{\beta\gamma}\,\,=\,\,
    \partial_\beta\delta \bar g_{\gamma\gamma}\,\,=\,\,
    0.\label{e:EdgeBC4}
\end{eqnarray}

When the analogs of the conditions in
Eqs.~(\ref{e:EdgeBC1})--(\ref{e:EdgeBC4}) are enforced along all four
edges of the $A\{\alpha\}$ face, they constitute both Dirichlet and
Neumann conditions for the metric perturbations, $\delta \bar
g^{\{\alpha\}}_{ab}$, on this face.  One convenient way to find
$\delta g^{\{\alpha\}}_{ab}$ that satisfy these boundary conditions is
to solve the bi-harmonic equation on this face:
\begin{eqnarray}
\left(\partial_\beta^{\,4} + 2\, \partial_\beta^{\,2}\,\partial_\gamma^{\,2}
+\partial_\gamma^{\,4}\right)\delta \bar g^{\{\alpha\}}_{ab}&=&0.
\label{e:Biharmonic2Da}
\end{eqnarray}
Solutions to the bi-harmonic equation are uniquely determined by
specifying both Dirichlet and Neumann boundary conditions on the
boundary of a compact domain.  For the intrinsic components on this
face, the solutions with these boundary conditions are trivial: $0 =
\delta \bar g^{\{\alpha\}}_{\beta\beta}= \delta \bar
g^{\{\alpha\}}_{\beta\gamma} = \delta \bar
g^{\{\alpha\}}_{\gamma\gamma}$.  For the non-trivial gauge components,
$\delta \bar g^{\{\alpha\}}_{\alpha\alpha}$, $\delta \bar
g^{\{\alpha\}}_{\alpha\beta}$, and $\delta \bar
g^{\{\alpha\}}_{\alpha\gamma}$, we use pseudo-spectral methods to
solve this equation numerically, as described in
\ref{s:BiharmonicMethods}.  We repeat this procedure to determine
$\delta \bar g_{ab}$ satisfying all the edge boundary conditions on
all the faces of each cubic region.

The solutions to Eq.~(\ref{e:Biharmonic2Da}) determine Dirichlet
boundary conditions for $\delta \bar g_{ab}$ on all the faces of cubic
region $\mathcal{B}_A$.  The normal derivatives of $\delta \bar
g_{ab}$ on the $A\{\alpha\}$ face are not prescribed, except the
requirement that they be compatible with the tangential derivatives on
the adjoining $A\{\beta\}$ faces.  The complete interior solutions for
$\delta \bar g_{ab}$ that are compatible with these boundary
conditions can be determined in a variety of ways.  For example the
three-dimensional Laplace equation could be solved for each component
of $\delta \bar g_{ab}$ with the Dirichlet boundary conditions $\delta
\bar g^{\{\alpha\}}_{ab}$ prescribed by the solutions to
Eq.~(\ref{e:Biharmonic2Da}).

This study has adopted the computationally more efficient approach
described in Sec.~\ref{s:Step3.1}.  This approach extrapolates the
values of $\delta \bar g^{\{\alpha\}}_{ab}$ from the $A\{\alpha\}$
face into the interior using expressions of the form $\delta\bar
g_{ab}=h(s^\alpha_{A})\,\delta \bar g^{\{\alpha\}}_{ab}$, where the
smooth function $h(s)$ is defined in
Eq.~(\ref{e:PartitionOfUnityhAlt}), and $s^\alpha_{A}$ is defined in
Eq.~(\ref{e:s_alphaDef}). The values of each component of $\delta \bar
g^{\{\alpha\}}_{ab}$ vanish on each edge of $\mathcal{B}_A$ because of
the boundary conditions, Eq.~(\ref{e:EdgeBC1}), imposed on the
solutions to Eq.~(\ref{e:Biharmonic2Da}).  It follows that the
extrapolations $\delta\bar g_{ab}=h(s^\alpha_{A})\,\delta \bar
g^{\{\alpha\}}_{ab}$ from the $A\{\alpha\}$ face will vanish on all
the other faces of the multicube structure.  These face extrapolations
can therefore be combined to give a complete interior solution for
$\delta \bar g_{ab}$ in region $\mathcal{B}_A$,
\begin{equation}
  \delta\bar g_{ab}=\sum_\alpha h(s^\alpha_{A})\,\delta \bar
  g^{\{\alpha\}}_{ab},\label{e:deltagbar}
\end{equation}
that automatically satisfies all the required Dirichlet boundary
conditions.  Adding the resulting metric perturbation to $\bar g_{ab}$
results in a new metric $\bar{\bar g}_{ab}$:
\begin{eqnarray}
  \bar{\bar g}_{ab} &=& \bar g_{ab} +\delta \bar g_{ab}.
  \label{e:gbarbarDef}
\end{eqnarray}
The boundary conditions imposed on $\delta \bar g_{ab}$ in this
construction ensure that $\bar{\bar g}_{ab}$ satisfies the two
important properties outlined at the beginning of this subsection. In
particular, the intrinsic components of $\bar{\bar g}_{ab}$ are
identical to those of $\bar g_{ab}$ on each cube face, and the
boundary conditions imposed on $\delta \bar g_{ab}$ ensure that the
extrinsic curvatures, $\bar{\bar K}^{\{\alpha\}}_{ab}$, vanish
identically along each edge of each cube.

\subsection{Step 3: Convert $\bar{\bar g}_{ab}$ into $\tilde g_{ab}$.}
\label{s:Step3.3}

Let $\bar{\bar g}_{ab}$ denote the global $C^{\,0}$ metric constructed
in Sec.~\ref{s:Step3.2}.  The goal of this third step is to convert
$\bar{\bar g}_{ab}$ into a metric $\tilde g_{ab}$ having two important
properties: first, the extrinsic curvatures $\tilde
K^{\{\alpha\}}_{ab}$ associated with $\tilde g_{ab}$ must vanish
identically on each face of each cubic region; and second, the metric
$\tilde g_{ab}$ must be identical to $\bar{\bar g}_{ab}$ on each face
of each cubic region.  We note that while $\bar{\bar g}_{ab}$ was
constructed to have intrinsic parts that are only $C^{\,0}$ across the
interface boundaries, within each region $\bar{\bar g}_{ab}$ is
actually smooth.

We define the unit normal vectors $\tilde n^{\{\alpha\}}_a$, the
projection tensors $\tilde P^{\{\alpha\}\,b}_a$, and the extrinsic
curvatures $\tilde K^{\{\alpha\}}_{ab}$ associated with the metric
$\tilde g_{ab}$ using expressions analogous to those given in
Eqs.~(\ref{e:nbarbarDef})--(\ref{e:KbarbarDefC}).  Similarly, we
define the differences between the $\tilde g_{ab}$ and $\bar{\bar
  g}_{ab}$ metrics, $\delta \bar{\bar g}_{ab}$, and the associated
differences between extrinsic curvatures, $\delta \bar{\bar
  K}^{\{\alpha\}}_{ab}$:
\begin{eqnarray}
  \delta \bar{\bar g}_{ab}&=&  \tilde g_{ab}-\bar{\bar g}_{ab},\\
  \delta \bar{\bar K}^{\{\alpha\}}_{ab}&=&\tilde K^{\{\alpha\}}_{ab}
  -\bar{\bar K}^{\{\alpha\}}_{ab}. 
\end{eqnarray}
We note that these differences are not assumed to be small.  To ensure
that $\tilde g_{ab}$ is identical to $\bar{\bar g}_{ab}$ on the cube
faces, we choose $\delta\bar{\bar g}_{ab}$ to be a smooth tensor in
the interior of each cubic region that satisfies
\begin{eqnarray}
\delta \bar{\bar g}_{ab}=0
  \label{e:MetricContinuitybarbar}
\end{eqnarray}
on each cube face $A\{\alpha\}$.  The condition that
the extrinsic curvature $\tilde K^{\{\alpha\}}_{ab}$ vanishes on the
$A\{\alpha\}$ face is equivalent to the following
condition on $\delta\bar{\bar K}^{\{\alpha\}}_{ab}$,
\begin{eqnarray}
    \delta \bar{\bar K}^{\{\alpha\}}_{ab}&=&-\bar{\bar K}^{\{\alpha\}}_{ab}.
    \label{e:DeltaKbarbarCondition}
\end{eqnarray} 

In analogy with Eq.~(\ref{e:KbarbarDefC}), an exact expression for
$\tilde K^{\{\alpha\}}_{ab}$ is given by
\begin{eqnarray}
  \tilde K^{\{\alpha\}}_{ab}
  &=&\tfrac{1}{2} \tilde P^{\{\alpha\}\,c}_a
  \tilde P^{\{\alpha\}\, d}_b \,\tilde n^{\{\alpha\}\,e}
  \left(\partial_e\,\tilde g_{cd}-\partial_c\,\tilde g_{ed}
  -\partial_d\,\tilde g_{ec}\right)\label{e:KtildeDefA}\\
  &=&\tfrac{1}{2} \tilde P^{\{\alpha\}\, c}_a
  \tilde P^{\{\alpha\}\,d}_b \,\tilde n^{\{\alpha\}\,e}
  \left(\partial_e\,\bar{\bar g}_{cd}-\partial_c\,\bar{\bar g}_{ed}
  -\partial_d\,\bar{\bar g}_{ec}\right)\nonumber\\
  &&+\tfrac{1}{2} \tilde P^{\{\alpha\}\, c}_a
  \tilde P^{\{\alpha\}\, d}_b \,\tilde n^{\{\alpha\}\,e}
  \left(\partial_e\,\delta\bar{\bar g}_{cd}-\partial_c\,\delta\bar{\bar g}_{ed}
  -\partial_d\,\delta\bar{\bar g}_{ec}\right).\label{e:KtildeDefB}
\end{eqnarray}
The metric perturbation $\delta \bar{\bar g}_{ab}$ vanishes on each
cube face, Eq.~(\ref{e:MetricContinuitybarbar}), therefore $\tilde
n^{\{\alpha\}\,a}=\bar{\bar n}^{\{\alpha\}\,a}$ and $\tilde
P^{\{\alpha\}\,b}_a=\bar{\bar P}^{\{\alpha\}\,b}_a$ on those faces as
well.  Consequently, Eq.~(\ref{e:KtildeDefB}) can be re-written as an
exact expression for $\delta\bar{\bar K}^{\{\alpha\}}_{ab}$ on those
faces:
\begin{eqnarray}
  \delta \bar{\bar K}^{\{\alpha\}}_{ab}
  &=&
  \tfrac{1}{2} \bar{\bar P}^{\{\alpha\}\, c}_a
  \bar{\bar P}^{\{\alpha\}\, d}_b \,\bar{\bar n}^{\{\alpha\}\,e}
  \left(\partial_e\,\delta\bar{\bar g}_{cd}-\partial_c\,\delta\bar{\bar g}_{ed}
  -\partial_d\,\delta\bar{\bar g}_{ec}\right)\label{e:deltaKbarbarDefA}\\
  &=& 
  \tfrac{1}{2}\left(\bar{\bar N}^{\{\alpha\}}\right)^{-1} 
  \bar{\bar P}^{\{\alpha\}\, c}_a \bar{\bar P}^{\{\alpha\}\, d}_b
  \partial_\alpha\,\delta\bar{\bar g}_{cd}.\label{e:deltaKbarbarDefB}
\end{eqnarray}
The second equality, Eq.~(\ref{e:deltaKbarbarDefB}), follows from the
fact that $\delta\bar{\bar g}_{ab}$ vanishes on the
$A\{\alpha\}$ face. This implies that $\partial_\beta
\delta\bar{\bar g}_{ab} = \partial_\gamma \delta\bar{\bar g}_{ab}=0$,
$\bar{\bar P}^{\{\alpha\}\,b}_{a}\bar{\bar n}^{\{\alpha\}}_b=0$ so
$\bar{\bar P}^{\{\alpha\}\,\alpha}_{a}=0$, and $\bar{\bar
  n}^{\{\alpha\}\,\alpha}= \left(\bar{\bar
  N}^{\{\alpha\}}\right)^{-1}$.

Let $\mathcal{N}^{\{\alpha\}}_{ab}$ denote the boundary conditions on
$\partial_\alpha \delta \bar{\bar g}^{\{\alpha\}}_{ab}$ on the
$A\{\alpha\}$ face. Only the intrinsic components of
$\delta\bar{\bar g}_{ab}$ contribute to the right side of
Eq.~(\ref{e:deltaKbarbarDefB}). Therefore
Eqs.~(\ref{e:DeltaKbarbarCondition}) and (\ref{e:deltaKbarbarDefB})
provide the needed boundary conditions for the intrinsic metric
components:
\begin{eqnarray} 
  \partial_\alpha \delta \bar{\bar g}^{\{\alpha\}}_{\beta\beta}
  &=&\mathcal{N}^{\{\alpha\}}_{\beta\beta}\,\,
    =\,\,-2 \bar{\bar N}^{\{\alpha\}} \bar {\bar K}^{\{\alpha\}}_{\beta\beta},
  \label{e:IntrinsicNeumanbarbarBCA}\\
  \partial_\alpha \delta \bar{\bar g}^{\{\alpha\}}_{\beta\gamma}
  &=&\mathcal{N}^{\{\alpha\}}_{\beta\gamma}\,\,
    =\,\,-2 \bar{\bar N}^{\{\alpha\}} \bar {\bar K}^{\{\alpha\}}_{\beta\gamma},
  \label{e:IntrinsicNeumanbarbarBCB}\\
  \partial_\alpha \delta \bar{\bar g}^{\{\alpha\}}_{\gamma\gamma}
  &=&\mathcal{N}^{\{\alpha\}}_{\gamma\gamma}\,\,
    =\,\,-2 \bar{\bar N}^{\{\alpha\}} \bar {\bar K}^{\{\alpha\}}_{\gamma\gamma}.
  \label{e:IntrinsicNeumanbarbarBCC}
\end{eqnarray}
The normal derivatives specified in
Eqs.~(\ref{e:IntrinsicNeumanbarbarBCA})--(\ref{e:IntrinsicNeumanbarbarBCC})
vanish along the edges of the $A\{\alpha\}$ face,
because $\bar{\bar K}^{\{\alpha\}}_{ab}$ was constructed to vanish
along those edges.  These edge conditions are needed to ensure that
the normal derivatives $\partial_\alpha \delta \bar{\bar
  g}^{\{\alpha\}}_{ab}$ are consistent along the $A\{\alpha\beta\}$
edge with the tangential derivatives $\partial_\alpha \delta \bar{\bar
  g}^{\{\beta\}}_{ab}$ from the adjoining
$A\{\beta\}$ face.

The normal derivative boundary conditions on the intrinsic components
of $\delta \bar{\bar g}^{\{\alpha\}}_{ab}$ in
Eqs.~(\ref{e:IntrinsicNeumanbarbarBCA})--(\ref{e:IntrinsicNeumanbarbarBCC})
are sufficient to guarantee that the entire extrinsic curvature
vanishes, $\tilde K_{ab}=0$, on the $A\{\alpha\}$ face.  The boundary
conditions on the gauge components are not fixed in this way, however.
These gauge boundary conditions can be chosen arbitrarily so long as
they vanish along each cube edge.  One choice is simply to set $0 =
\mathcal{N}^{\{\alpha\}}_{\alpha\alpha} =
\mathcal{N}^{\{\alpha\}}_{\alpha\beta} =
\mathcal{N}^{\{\alpha\}}_{\alpha\gamma}$ everywhere on the
$A\{\alpha\}$ face.  Somewhat better numerical convergence can be
achieved, however, by choosing $\mathcal{N}^{\{\alpha\}}_{ab}$ to make
the second derivatives $\partial_\alpha\partial_\beta\delta \bar{\bar
  g}^{\{\alpha\}}_{ab}$ consistent along the $\{\alpha\beta\}$ edge
with their values on the adjacent $A\{\beta\}$ face.  These conditions
on $\mathcal{N}^{\{\alpha\}}_{ab}$ along the $\{\alpha\beta\}$ edge
require
\begin{eqnarray}
  \partial_\beta \mathcal{N}^{\{\alpha\}}_{\alpha\alpha} &=&
  \partial_\alpha\mathcal{N}^{\{\beta\}}_{\alpha\alpha} \,\,
  =\,\,-2\partial_\alpha\left(\bar{\bar N}^{\{\beta\}} \bar{\bar
    K}^{\{\beta\}}_{\alpha\alpha}\right),
  \label{e:edgebc1}\\
  \partial_\beta \mathcal{N}^{\{\alpha\}}_{\alpha\beta} &=&0
  \label{e:edgebc2},\\
  \partial_\beta \mathcal{N}^{\{\alpha\}}_{\alpha\gamma}
  &=&\partial_\alpha\mathcal{N}^{\{\beta\}}_{\alpha\gamma} \,\,
  =\,\,-2\partial_\alpha\left(\bar{\bar N}^{\{\beta\}} \bar{\bar
    K}^{\{\beta\}}_{\alpha\gamma}\right).
  \label{e:edgebc3}
\end{eqnarray}
Analogous conditions on each edge of the $A\{\alpha\}$ face provide
Neumann boundary conditions for the gauge components of
$\mathcal{N}^{\{\alpha\}}_{ab}$ along the edges of this face.
Together with the Dirichlet conditions
$\mathcal{N}^{\{\alpha\}}_{ab}=0$ along these edges, they provide the
boundary conditions needed to determine
$\mathcal{N}^{\{\alpha\}}_{\alpha\alpha}$,
$\mathcal{N}^{\{\alpha\}}_{\alpha\beta}$, and
$\mathcal{N}^{\{\alpha\}}_{\alpha\gamma}$ everywhere on this face by
solving the two-dimensional bi-harmonic equations
\begin{eqnarray}
\left(\partial_\beta^{\,4} + 2\, \partial_\beta^{\,2}\,\partial_\gamma^{\,2}
+\partial_\gamma^{\,4}\right)\mathcal{N}^{\{\alpha\}}_{\alpha\alpha}&=&0,
\label{e:Biharmonic2DNA}\\
\left(\partial_\beta^{\,4} + 2\, \partial_\beta^{\,2}\,\partial_\gamma^{\,2}
+\partial_\gamma^{\,4}\right)\mathcal{N}^{\{\alpha\}}_{\alpha\beta}&=&0,
\label{e:Biharmonic2DNB}\\
\left(\partial_\beta^{\,4} + 2\, \partial_\beta^{\,2}\,\partial_\gamma^{\,2}
+\partial_\gamma^{\,4}\right)\mathcal{N}^{\{\alpha\}}_{\alpha\gamma}&=&0.
\label{e:Biharmonic2DNC}
\end{eqnarray}
The pseudo-spectral numerical methods used in this study
to solve this equation are described in~\ref{s:BiharmonicMethods}. 

Equation~(\ref{e:MetricContinuitybarbar}) provides Dirichlet boundary
conditions for $\delta\bar{\bar g}_{ab}$, and the
$\mathcal{N}^{\{\alpha\}}_{ab}$ from
Eq.~(\ref{e:IntrinsicNeumanbarbarBCA})--(\ref{e:IntrinsicNeumanbarbarBCC})
together with the solutions to
Eqs.~(\ref{e:Biharmonic2DNA})--(\ref{e:Biharmonic2DNC}) provide
Neumann boundary conditions on each face of each cubic region.  The
perturbations $\delta \bar{\bar g}_{ab}$ can therefore be determined
throughout the region by solving the three-dimensional biharmonic
equation,
\begin{eqnarray}
   \left(\partial_x^{\,4} +\partial_y^{\,4} +\partial_z^{\,4}
  +2\,\partial_x^{\,2}\,\partial_y^{\,2}
  +2\,\partial_x^{\,2}\,\partial_z^{\,2}
  +2\,\partial_y^{\,2}\,\partial_z^{\,2}\right)\delta \bar{\bar
    g}_{ab}=0,
  \label{e:Biharmonicbarbar3D}
\end{eqnarray}
with these boundary conditions.  The pseudo-spectral numerical methods
used here to solve Eq.~(\ref{e:Biharmonicbarbar3D}) for $\delta
\bar{\bar g}_{ab}$ are discussed in \ref{s:BiharmonicMethods}.  Adding
the resulting $\delta \bar{\bar g}_{ab}$ to $\bar{\bar g}_{ab}$
results in the new metric $\tilde g_{ab}$:
\begin{eqnarray}
  \tilde g_{ab} &=& \bar{\bar g}_{ab} +\delta \bar{\bar g}_{ab}.
  \label{e:gtildeDef}
\end{eqnarray}
The boundary conditions imposed on $\delta \bar{\bar g}_{ab}$ ensure
that $\tilde g_{ab}$ satisfies the two important properties outlined
at the beginning of this subsection; namely, the components of $\tilde
g_{ab}$ are identical to those of $\bar{\bar g}_{ab}$ on each cube
face, and the extrinsic curvatures, $\tilde K^{\{\alpha\}}_{ab}$,
vanish identically on each cube face.  It follows that the intrinsic
components of $\tilde g_{ab}$ and the extrinsic curvatures $\tilde
K^{\{\alpha\}}_{ab}$ are continuous across the interface boundaries
between all the cubic regions.  Therefore $\tilde g_{ab}$ satisfies
the Israel~\cite{Israel1966} junction conditions across all the
boundaries of the multicube structure, and is therefore $C^1$
globally.

\section{Numerical Examples}  
\label{s:NumericalExamples} 
 
Multicube structures for a collection of manifolds have been developed
here to test the numerical reference metric construction methods
described in Secs.~\ref{s:C0ReferenceMetrics} and
\ref{s:C1ReferenceMetrics}.  All the multicube structures used in
these examples satisfy the local reflection symmetry property
described in Sec.~\ref{s:C0ReferenceMetrics}.  This condition is
needed to permit the construction of flat metrics in the neighborhood
of each vertex having uniform dihedral angles around each edge of the
cubic regions.  These example manifolds are listed in
Table~\ref{t:manifold_list}, including their Thurston geometrization
classes (see Ref~\cite{Petronio}).  They include representatives from
five of the eight Thurston geometrization classes, missing only the
$S\!L_2$, $Nil$, and $S\!ol$ classes.

\begin{table}[!b] 
\scriptsize \renewcommand{\arraystretch}{1.5}
\begin{center}
\caption{Manifolds used in numerical tests of the $C^1$ metric
    construction methods developed in Secs.~\ref{s:C0ReferenceMetrics}
    and \ref{s:C1ReferenceMetrics}.  First part of the table lists
    this multicube structures constructed by hand, while the second
    part lists those constructed from triangulations by the code
    described in~\ref{s:TriangulationToMultiCubeCode}. Names used for
    the manifolds constructed from triangulations are those used in
    Ref.~\cite{Regina}.  The L(p,q) manifolds are lens spaces, i.e.,
    quotients of the three-sphere $S^3$ with a discrete group
    characterized by parameters (p,q).  The manifolds $S2\times S1$,
    $T\times S1$, $KB/n2\times\!\!\sim\!\!S1$, and
    $SFS[B:(p{}_1,q{}_1)(p{}_2,q{}_2)(p{}_3,q{}_3)]$, are
    Seifert fibered spaces.  $S1$ represents a circle. The $\times$
    operator is the Cartesian product, e.g. $S2 \times S1$, while
    $\times\!\!\sim$ is the twisted product used to undo the
    non-orientability of the base manifold, e.g. in
    $KB/n2\times\!\!\sim\!\!S1$.  The base spaces B include the two-sphere
    $S2$, the real projective plane $RP2/n2$, the Klein bottle $KB/n2$, and
    the two-torus $T$.  The parameters, e.g. $(p{}_1,q{}_1)$, describe
    ``singular'' fibers whose neighborhoods have been replaced by
    fibers twisted by an amount determined by $(p{}_1,q{}_1)$.
    \label{t:manifold_list}}
  \begin{tabular}{c}
    \toprule 
    Three Dimensional Multicube Structures Constructed by Hand\\
  \end{tabular}
\begin{tabular}{l|l||l|l}  
  \midrule Manifold & Geometry Class & Manifold & Geometry Class\\
  \midrule Three-Torus (E1) & $E^3$ & Half-Turn Space (E2) & $E^3$\\
  Quarter-Turn Space (E3) & $E^3$ & Third-Turn Space (E4) & $E^3$\\
  Sixth-Turn Space (E5) & $E^3$ & Hantzsche-Wendt Space (E6) & $E^3$\\
  Three-Sphere (S3) & ${S}^3$ & $S2\times S1$ & $S^2\times S^1$ \\
  $G2\times S1$ & $H^2\times S^1$ & Seifert-Weber Space &
  $H^3$\\ Poincar\'e Dodecahedral Space & $S^3$ \\
  \bottomrule
\end{tabular}
\end{center} 
\begin{center}
  \begin{tabular}{c}
    \toprule
    Three Dimensional Multicube Structures Constructed From
    Triangulations\\
  \end{tabular}
\begin{tabular}{l|l||l|l||l|l}  
  \midrule
  Manifold & Geometry Class & Manifold & Geometry Class
  & Manifold & Geometry Class\\
  \midrule
  $L(5,2)$ & $S^3$ & $L(40,19)$ & $S^3$
  & $SFS[S2:(2,1)(2,1)(7,-6)]$ & $S^3$\\
  $L(8,3)$ & $S^3$ & $L(44,21)$ & $S^3$
  & $SFS[S2:(2,1)(2,1)(8,-7)]$ & $S^3$ \\ 
  $L(10,3)$ & $S^3$ & $T \times S1$ & $E^3$
  & $SFS[S2:(2,1)(2,1)(9,-8)]$ & $S^3$\\
  $L(12,5)$ & $S^3$ & $KB/n2\times\!\!\sim\!\!S1$ & $E^3$
  &  $SFS[S2:(2,1)(2,1)(10,-9)]$ & $S^3$ \\ 
  $L(16,7)$ & $S^3$ & $SFS[RP2/n2:(2,1)(2,-1)]$ & $E^3$
  &  $SFS[S2:(2,1)(2,1)(11,-10)]$ & $S^3$ \\
  $L(20,9)$ & $S^3$ & $SFS[S2:(2,1)(2,1)(2,-1)]$ & $S^3$
  & $SFS[S2:(2,1)(3,1)(5,-4)]$ & $S^3$ \\
  $L(24,11)$ & $S^3$ & $SFS[S2:(2,1)(2,1)(3,-2)]$ & $S^3$
  & $SFS[S2:(2,1)(3,2)(3,-1)]$ & $S^3$ \\
  $L(28,13)$ & $S^3$ & $SFS[S2:(2,1)(2,1)(4,-3)]$ & $S^3$
  & $SFS[S2:(2,1)(4,1)(4,-3)]$ & $S^3$ \\
  $L(32,15)$ & $S^3$ & $SFS[S2:(2,1)(2,1)(5,-4)]$ & $S^3$
  & $SFS[S2:(3,1)(3,1)(3,-2)]$ & $S^3$ \\
  $L(36,17)$ & $S^3$ & $SFS[S2:(2,1)(2,1)(6,-5)]$ & $S^3$ \\
  \bottomrule
\end{tabular}
\end{center}
\end{table}

Some of the multicube structures used in these examples were
constructed by hand, while most were constructed from triangulations
obtained from Ref.~\cite{Regina} using the method developed in
Ref.~\cite{Lindblom2013}.  The multicube structures constructed from
triangulations were done automatically by the code described in
\ref{s:TriangulationToMultiCubeCode}.  Those constructed by hand
include the Three-Torus (E1), S3 and S$\times$S1, as described in
Ref.~\cite{Lindblom2013}.  Manifolds of the form $G_n\times S^1$,
where $G_n$ is the compact orientable two-manifold with genus number
$n$, can be constructed easily by hand from the two-dimensional
multicube structures developed for arbitrary $G_n$ in
Ref.~\cite{Lindblom2015}.  This study includes G${}_2\times$S1 as an
example.  Multicube structures have also been constructed by hand for
several manifolds that can be defined by identifying the faces of
three-dimensional polygonal solids.  The numerical examples presented
here include the Poincar\'e dodecahedral
space~\cite{ThrelfallSeifert1931}, Seifert-Weber
space~\cite{SeifertWeber1933}, and all six compact orientable
three-manifolds that admit flat metrics (sometimes called
E1-E6)~\cite{riazuelo2004, Hitchman2018}, and the Hantzsche-Wendt
space~\cite{Hantzsche1934} (also called E6).
\ref{s:3DMulticubeManifolds} gives the complete descriptions of the
previously unpublished multicube structures constructed by hand for
this study, along with a representative selection of those constructed
automatically from triangulations by the code described in
\ref{s:TriangulationToMultiCubeCode}.

Reference metrics $\tilde g_{ab}$ have been constructed numerically
for each of the manifolds listed in Table~\ref{t:manifold_list}.  The
methods developed in Secs.~\ref{s:C0ReferenceMetrics} and
\ref{s:C1ReferenceMetrics} are designed to make the intrinsic parts of
$\tilde g_{ab}$ continuous across the interface boundaries between the
cubic regions, and also to make the associated extrinsic curvatures,
$\tilde K^{\{\alpha\}}_{ab}$, vanish on each interface boundary.
These conditions satisfy the Israel junction
conditions~\cite{Israel1966} that ensure $\tilde g_{ab}$ is $C^1$
across those interfaces.

\begin{figure}[!htb]
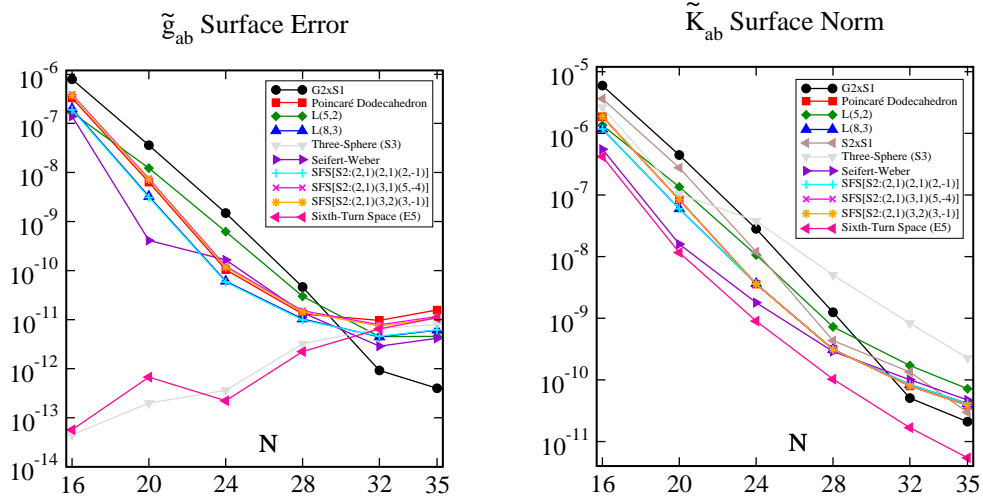
 
\begin{center}
  \subfigure{
    \includegraphics[width=0.43\textwidth]{Fig5a.eps}
    \label{f:gSurfaceError1}
  }\hspace{1.0cm}\vspace{-0.2cm}
  \subfigure{
    \includegraphics[width=0.43\textwidth]{Fig5b.eps}
                    \label{f:KSurfaceError1}
  }
\end{center}
\caption{Graphs representing $L_2$ norms, for different values of the
  spatial resolution $N$, of the intrinsic metric discontinuities of
  $\tilde g_{ab}$ across the multicube interface boundaries in the
  left Fig.~\ref{f:gSurfaceError1}, and $L_2$ norms of the associated
  extrinsic curvatures $\tilde K^{\{\alpha\}}_{ab}$ of those
  boundaries in the right Fig.~\ref{f:KSurfaceError1}. In
  Fig.~\ref{f:gSurfaceError1} the graph for S2$\times$S1 is not shown
  because the errors are at the $10^{-15}$ level.
\label{f:SurfaceError1}}
\end{figure}

\begin{figure}[!htb]
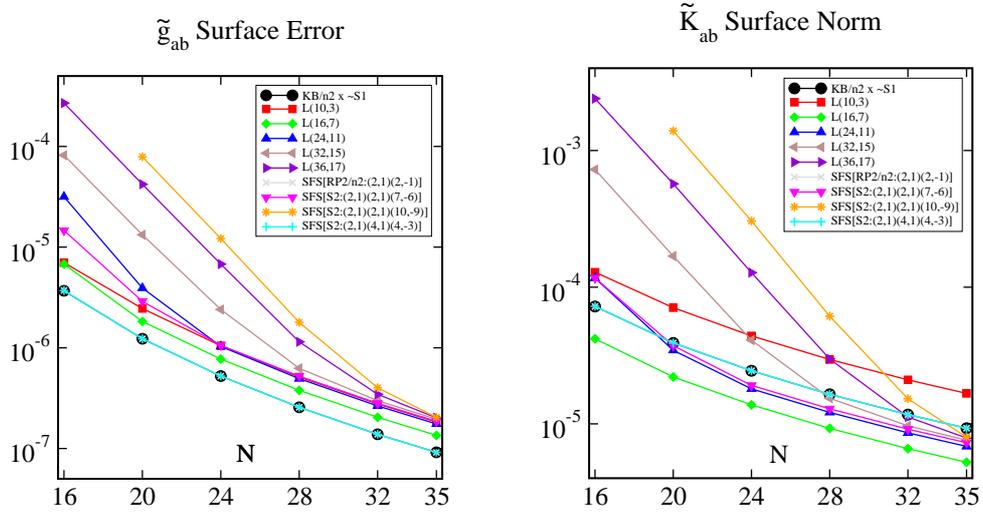
 
\begin{center}
  \subfigure{
    \includegraphics[width=0.43\textwidth]{Fig6a.eps}
    \label{f:gSurfaceError2}
  }\hspace{1.0cm}\vspace{-0.2cm}
  \subfigure{
    \includegraphics[width=0.43\textwidth]{Fig6b.eps}
                    \label{f:KSurfaceError2}
  }
\end{center}
\caption{Graphs representing $L_2$ norms, for different
  values of the spatial resolution $N$, of the intrinsic metric
  discontinuities of $\tilde g_{ab}$ across the multicube interface
  boundaries in the left Fig.~\ref{f:gSurfaceError2}, and $L_2$ norms
  of the associated extrinsic curvatures $\tilde K^{\{\alpha\}}_{ab}$
  of those boundaries in the right Fig.~\ref{f:KSurfaceError2}. 
\label{f:SurfaceError2}}
\end{figure}

The methods introduced in Secs.~\ref{s:C0ReferenceMetrics} and
\ref{s:C1ReferenceMetrics} have been implemented numerically for this
study in the SpEC pseudo-spectral code (developed originally by the
Caltech/Cornell numerical relativity collaboration~\cite{Kidder2000a,
  Scheel2006, Szilagyi:2009qz}).  Figures~\ref{f:SurfaceError1} and
\ref{f:SurfaceError2} show $L_2$ norms of the surface discontinuities
of the intrinsic parts of $\tilde g_{ab}$ and the extrinsic curvatures
$\tilde K^{\{\alpha\}}_{ab}$ as functions of the spatial resolution
parameter $N$ (the number of spectral collocation points used in each
dimension) for most of the manifolds listed in
Table~\ref{t:manifold_list}.  These $L_2$ norms were computed by
averaging the squares of all the intrinsic components of each tensor
over all the grid points on all interface surfaces, and finally taking
the square root of this average.  These results show that the
numerical methods developed and implemented here produce $C^1$
reference metrics having small errors that converge toward satisfying
the Israel junction conditions as the spatial resolution is increased.
The results for the manifolds not included in these graphs are similar
to those shown in Fig.~\ref{f:SurfaceError2} (except for the flat
manifolds E1-E4 and E6 whose $\tilde K_{ab}$ errors are at or below
the $10^{-12}$ level for all $N$).

The results of these numerical tests have been divided into two
groups.  Those represented in Fig.~\ref{f:SurfaceError1} have
significantly smaller errors than those shown in
Fig.~\ref{f:SurfaceError2}.  The reason for these differences appears
to be the amount of distortion caused by the dihedral angles needed to
allow the cubic regions to fit together without introducing conical
edge singularities.  Higher resolutions are needed to represent models
having larger distortions at a particular accuracy level.  All the
manifolds in the larger error group, Fig.~\ref{f:SurfaceError2}, have
some edges with small dihedral angles, $\min (\psi)\leq 2\pi/6$, while
those in the smaller error group, Fig.~\ref{f:SurfaceError1}, have
larger minimum dihedral angles $\min( \psi)\geq 2\pi/5$ (except for
G${}_2$$\times$S1 and Sixth-Turn Space, E5, which have
$\min(\psi)=2\pi/6$).

The surface errors in $\tilde g_{ab}$ and $\tilde K_{ab}$ for the
examples shown in Fig.~\ref{f:SurfaceError1} decrease (approximately)
exponentially with increasing $N$ for $N\leq 28$.  Double precision
roundoff error is probably limiting convergence in these cases for
$N>28$.  Some of the examples in Fig.~\ref{f:SurfaceError2} also show
exponential convergence for $N\leq 28$.  However most of the examples
in Fig.~\ref{f:SurfaceError2} show slower power law convergence in
$N$.  For example the errors in one of the slowest converging cases,
$KB/n2\,\, \times\!\!\sim S^1$, are well fit by the power laws
$N^{-14/3}$ for $\tilde g_{ab}$ and $N^{-8/3}$ for $\tilde
K^{\{\alpha\}}_{ab}$.  There is some indication that the examples in
Fig.~\ref{f:SurfaceError2} with exponential convergence transition to
power law convergence for larger values of $N$.  This transition is
probably caused by errors due to discontinuities in the mixed partial
derivatives of $\delta \bar{\bar g}_{ab}$ at some of the edges.  These
discontinuities are caused by disagreements between the tangential
derivatives of the Neumann boundary data on the faces that intersect
along those edges.  At some resolution these higher-order
discontinuity errors become dominant and power law convergence takes
over.  The examples in Fig.~\ref{f:SurfaceError2} with the largest
errors are also those with the most distorted multicube structures,
some with dihedral angles as small as $\min(\psi)=2\pi/8$.  This
supports the idea that the larger distortions cause the larger errors
at a given resolution $N$.

\section{Discussion}
\label{s:Discussion}

New methods have been presented in Secs.~\ref{s:C0ReferenceMetrics}
and \ref{s:C1ReferenceMetrics} for building three-dimensional
differentiable manifolds numerically.  These methods involve the
construction of $C^{\,1}$ reference metrics that are used to construct
special Jacobians to define the continuity of tensors, and a covariant
derivative to define the differentiability of those tensors, across
the interface boundaries between coordinate charts.  These methods
have been applied in Sec.~\ref{s:NumericalExamples} to a selection of
forty three-dimensional manifolds, including examples from five of the
eight Thurston geometrization classes.  Test results on these examples
show that the methods developed in Secs.~\ref{s:C0ReferenceMetrics}
and \ref{s:C1ReferenceMetrics}, and our implementation of those
methods in the SpEC pseudo-spectral code, are numerically convergent.

The methods developed here are general enough to be applied to a
larger variety of differentiable three-manifolds than has been studied
previously using existing numerical methods.  However, the methods
presented here make very restrictive assumptions about the multicube
structures to which they can be applied.  Perhaps the most obvious
limitation is the assumption in Sec.~\ref{s:C0ReferenceMetrics} that
the multicube structure exhibit a particular local reflection
symmetry.  A diverse collection of manifolds that satisfy this
restriction has been constructed, however, this assumption is not
satisfied by most multicube structures.  We do not think that this
assumption is essential.  It was made here because it was easy to
implement numerically in our code.  We think it will be possible to
relax this assumption.  We plan to investigate ways to do that in a
future study.

Another obvious limitation of the results presented in
Sec.~\ref{s:NumericalExamples} is the relatively slow numerical
convergence of the reference metrics constructed on manifolds having
highly distorted multicube structures.  One significant part of this
problem is probably caused by the discontinuities in the derivatives
of the Neumann boundary data used to determine the $C^{\,1}$ reference
metrics in Sec.~\ref{s:Step3.3} (at cube edges where some intrinsic
metric component is present on both faces, e.g.  the $\tilde
g_{\gamma\gamma}$ component along the $ A{\{\alpha\beta\}}$ edge).  We
think this particular problem can be ameliorated by enforcing somewhat
different boundary conditions on the gauge components of the metric in
Sec.~\ref{s:Step3.2}.  We plan to investigate this and other
approaches to improving the numerical convergence of these methods in
a future study.

Most of the differential equations used in the physical sciences,
e.g. systems of symmetric hyperbolic evolution equations, or systems
of second-order elliptic equations, require specifying some
combination of the values of fields and their derivatives at the
boundaries of computational domains.  The $C^{\,1}$ reference metrics
developed in this paper are sufficient to provide the needed
transformations of these data at the interface boundaries between
coordinate patches.  We showed in Ref.~\cite{Lindblom2015} that the
differentiable structures produced by different $C^{\,1}$ reference
metrics are equivalent.  The needed continuity of the boundary data at
the interfaces between computational domains can therefore be done
correctly and exactly using the $C^{\,1}$ reference metrics
constructed here.  There is no fundamental need to refine these
reference metrics by increasing their global differentiability.

For various reasons it may be desirable, however, to transform these
metrics further to produce metrics that are smoother at the interface
boundaries, or perhaps that have more uniform spatial structures which
can be resolved numerically at lower resolutions.  In
Ref.~\cite{Lindblom2015} we used numerical Ricci flow to evolve the
$C^{\,1}$ reference metrics developed there for two-dimensional
manifolds.  Ricci flow is a system of parabolic evolution equations
that transform $C^{\,1}$ initial data into $C^\infty$ solutions at
later times \cite{Chow2004, Hamilton1988, Chow1991, Chen2006,
  DeTurck1983}.  The initial metrics for Ricci flow are required to
have bounded curvatures~\cite{Bemelmans1984, Bando1987} to ensure that
even very short evolutions become real analytic.  The Israel junction
conditions~\cite{Israel1966} guarantee that while our $C^1$ reference
metrics may have curvature discontinuities across the interfaces, they
will not be unbounded there. Ricci flow in two dimensions also evolves
all initial data into constant curvature geometries.  We plan to use
numerical Ricci flow to evolve the three-dimensional $C^{\,1}$
reference metrics produced here in a future study.  In three
dimensions, Ricci flow may form singularities before the manifold
attains constant curvature, even for manifolds like the Three-Sphere
(S3) having very simple topologies~\cite{Garfinkle2008}.  While there
is no guarantee that the Ricci flow of our $C^{\,1}$ metrics will
necessarily produce more uniform geometries, it will be interesting to
see what happens.  If singularities occur then it will be interesting
to explore the nature of those singularities.  If these evolutions
proceed to uniform curvature solutions, then it will be interesting to
determine and to verify that the resulting geometries satisfy the
appropriate properties associated with their Thurston geometrization
classes.

Finally, we plan to use the reference metrics developed here in a
future study to solve Einstein's equations numerically on a diverse
collection of manifolds.  Solving Einstein's equations involves
finding solutions to an elliptic system to obtain acceptable initial
data, and then to evolve those data using a system of hyperbolic
equations that determine the structure of the resulting spacetime.
The appropriate representation of Einstein's equations to use in
spacetimes with non-trivial topologies was developed in
Ref.~\cite{Lindblom2014}.  We plan to use those methods to explore
solutions representing cosmological models evolved from initial data
on a diverse collection of compact three-manifolds.  It might also be
interesting to explore solutions to Einstein's equations on manifolds
with non-trivial topologies and asymptotically flat initial data.
These geometries are expected to evolve into black-hole
spacetimes~\cite{Friedman1995, Eichmair2013, Galloway2017}.  with any
non-trivial topological structures hidden behind event horizons.  By
studying these evolutions, it will be interesting to explore whether
observers outside the black holes could identify the presence of these
topological structures in some indirect way.

\section*{Acknowledgments}

This research was supported in part by NSF grant 2012857 to the
University of California at San Diego.

\appendix
\counterwithin*{figure}{section}
\counterwithin*{table}{section}
\section{Converting Three-Manifold Triangulations to Multicube Structures}
\label{s:TriangulationToMultiCubeCode}

This appendix describes the method used by our code to convert a
three-dimensional triangulation into a multicube structure.  A
three-dimensional triangulation consists of a set of tetrahedra and
the identification maps that identify each tetrahedron face with the
appropriate face of its neighbor.  These face identifications are
determined by specifying which vertices of one tetrahedron are
identified with which vertices of its neighbor.  Large numbers of
triangulations specified in this way are published in the Regina
catalog~\cite{Regina}.  Our code is designed to read the triangulation
structures exported into files by the Regina software.

\begin{figure}[ht]
\subfigure{
  \includegraphics[width=0.23\textwidth]{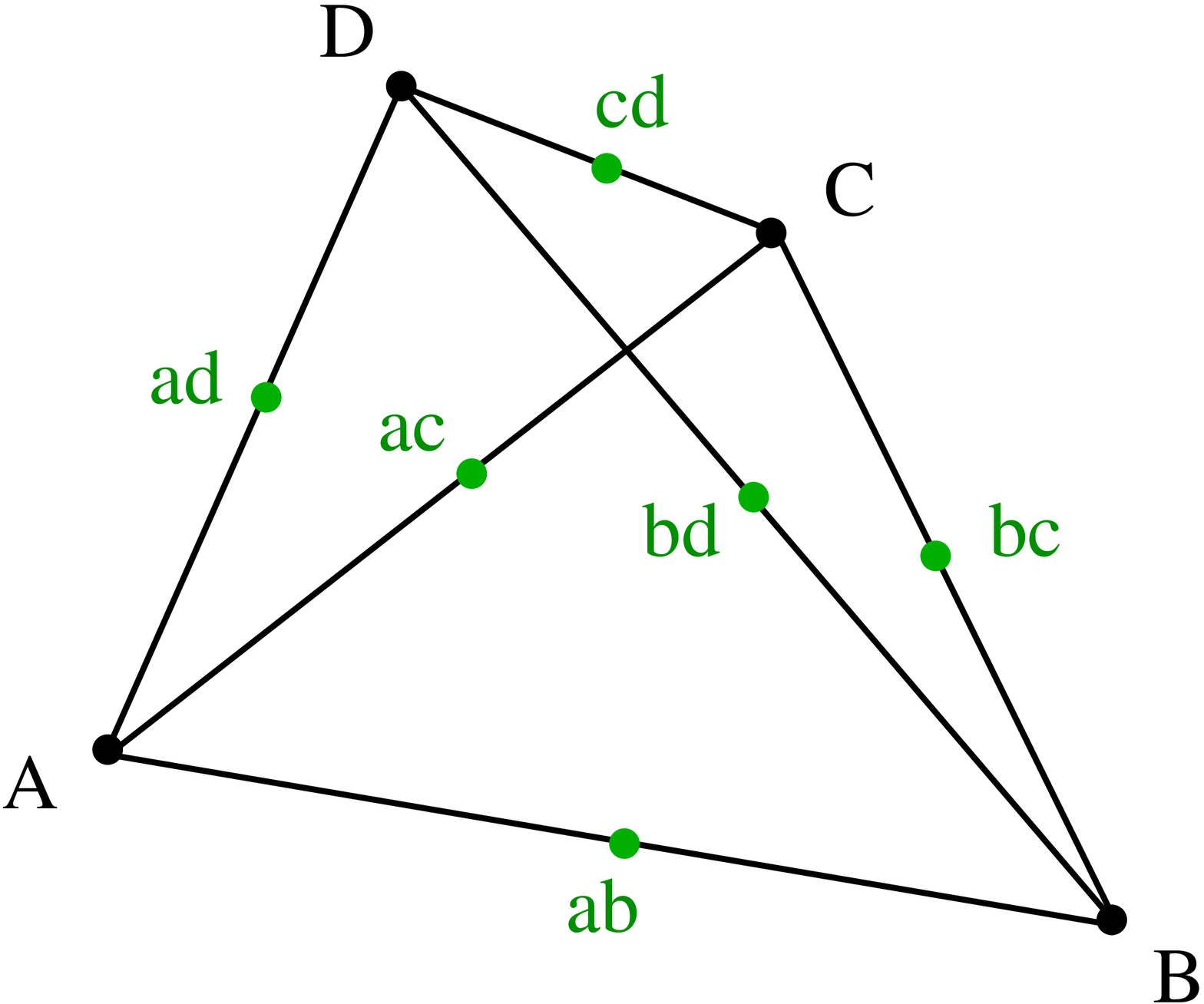}
\label{f:CubedTetrahedron0}}
\subfigure{
  \includegraphics[width=0.23\textwidth]{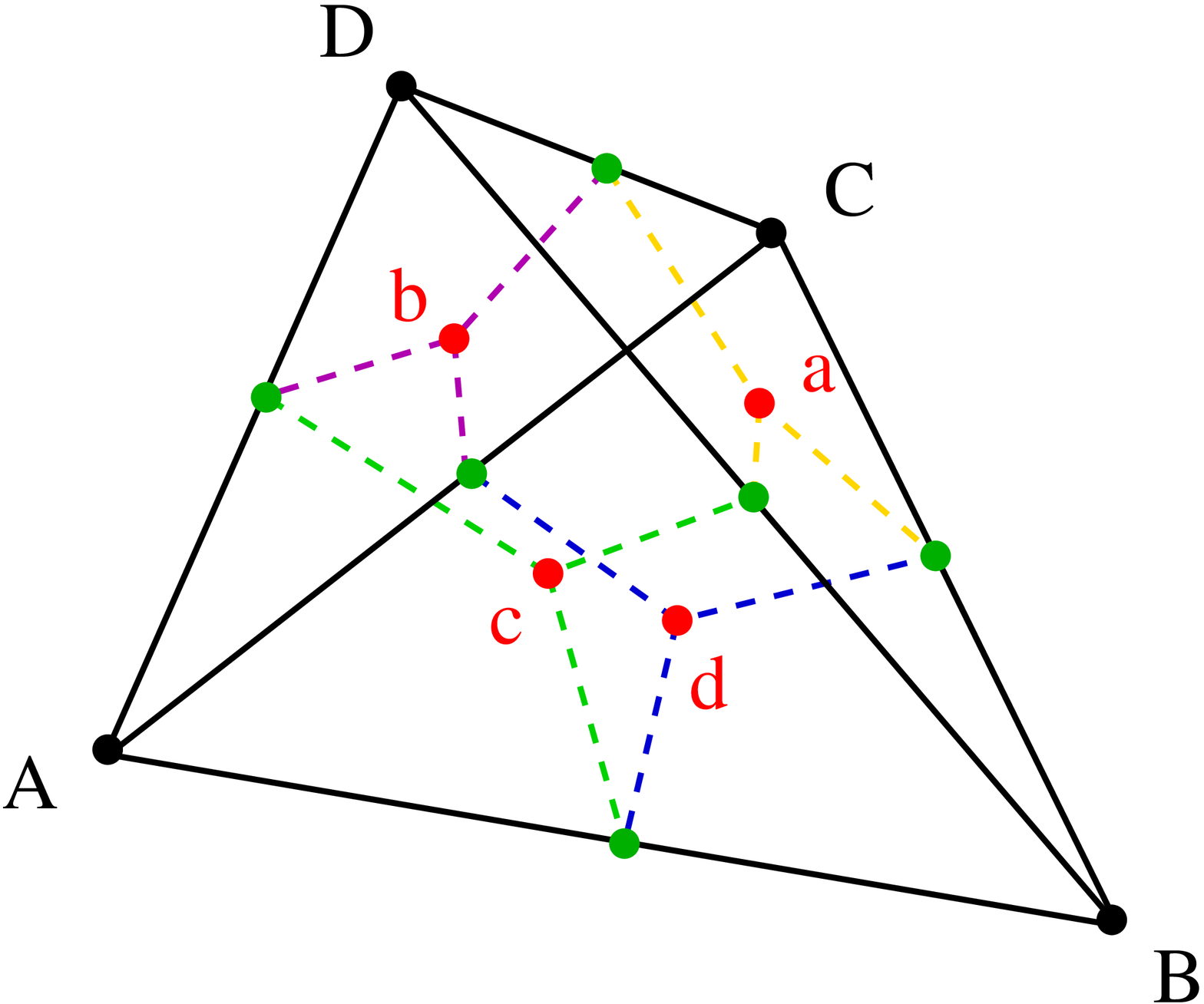}
\label{f:CubedTetrahedron1}}  
\subfigure{
  \includegraphics[width=0.23\textwidth]{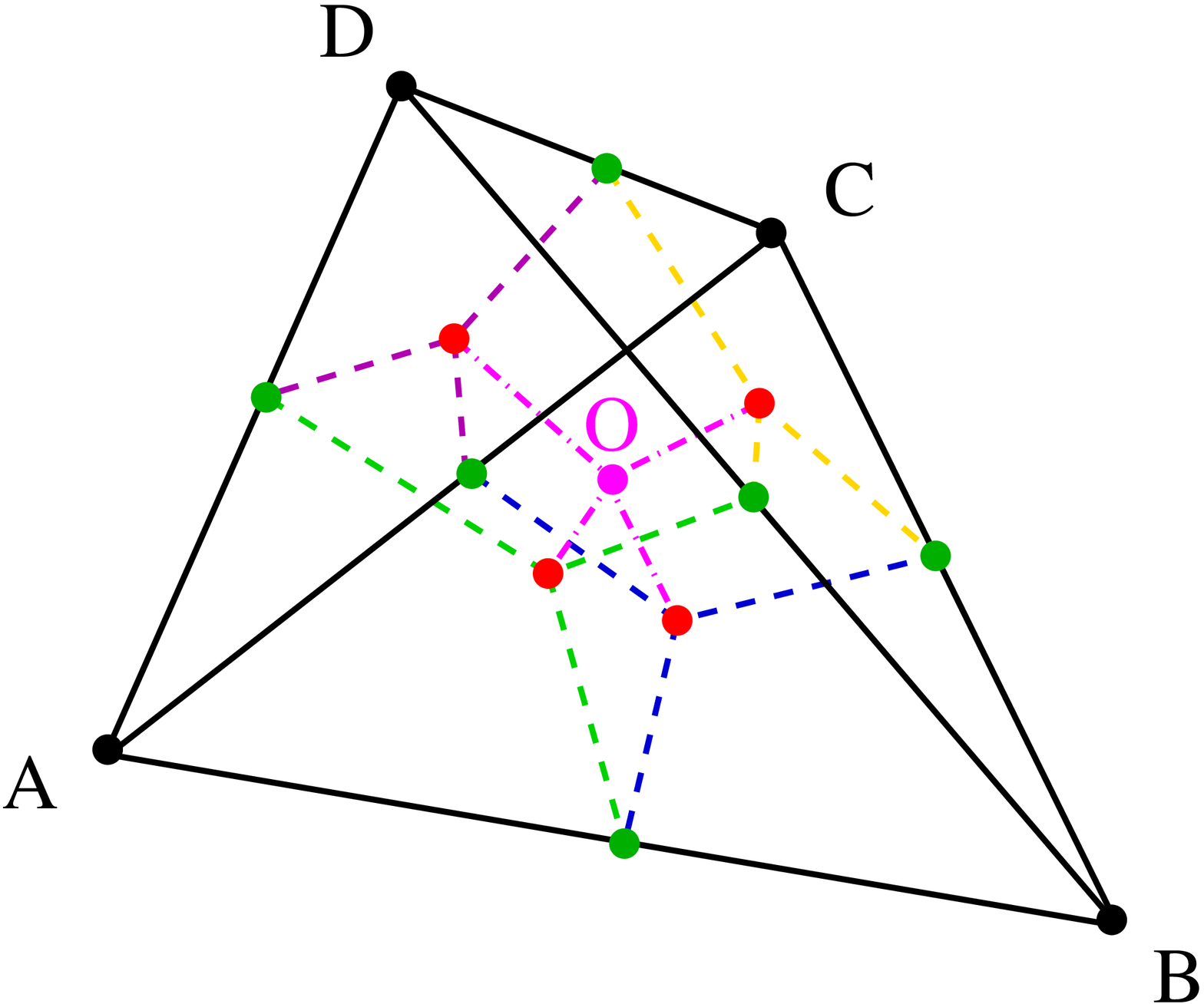}
\label{f:CubedTetrahedron2}}
\subfigure{
  \includegraphics[width=0.23\textwidth]{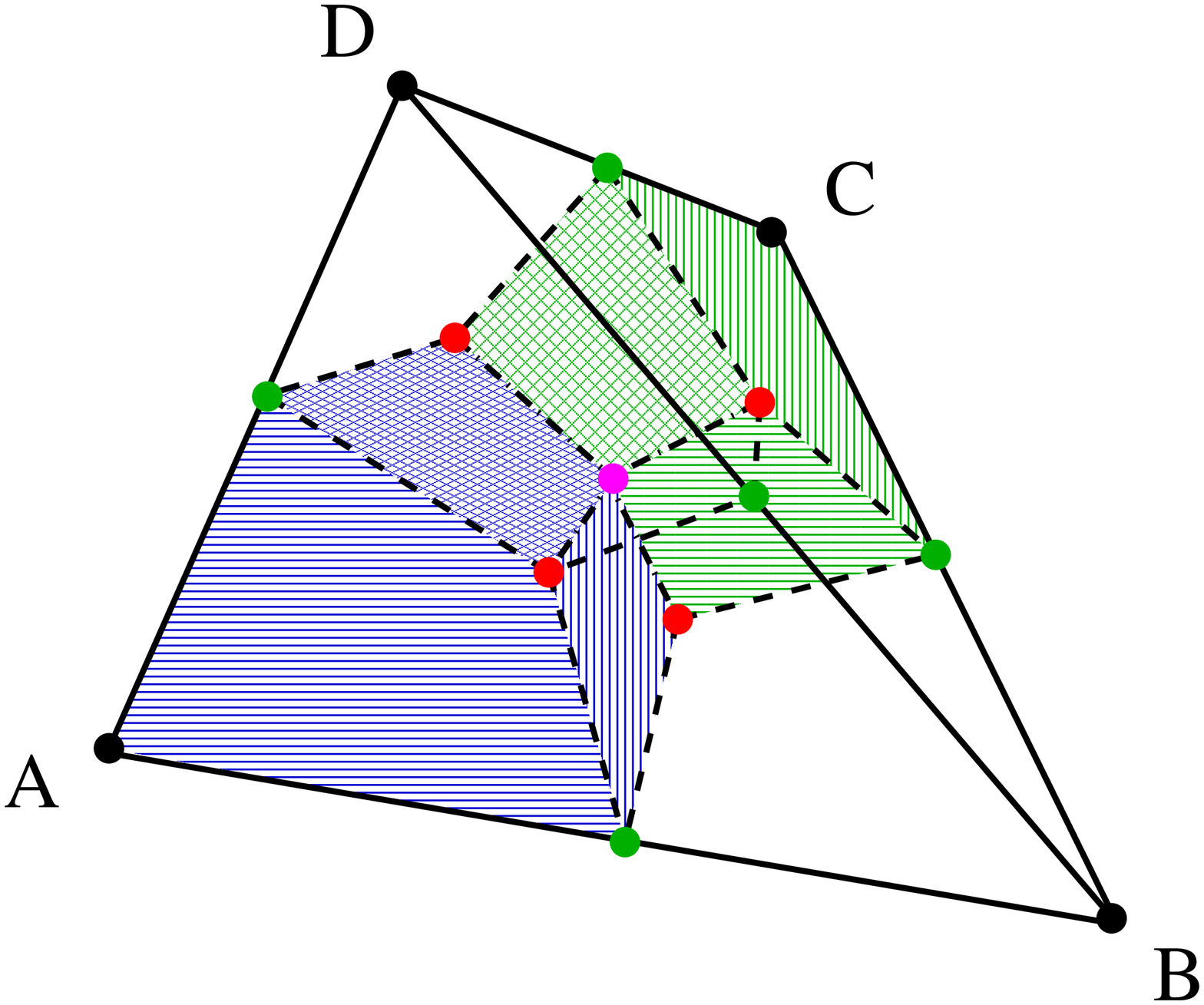}
\label{f:CubedTetrahedron3}} 
     
\caption{ In Fig.~\ref{f:CubedTetrahedron0}
  label the vertexes of the tetrahedron ``A'', ``B'', ``C'' and ``D'', 
  and add vertexes at the midpoints of each edge.  In
  Fig.~\ref{f:CubedTetrahedron1} add additional vertexes at the
  centroid of each face of the tetrahedron, labeled ``a'' for the
  centroid of face ``BCD'', ``b'' for face ``ACD'', etc.  Then add
  additional edges (shown as dashed line segments) connecting each
  centroid to the midpoint of each adjoining edge.  In
  Fig.~\ref{f:CubedTetrahedron2} add one additional vertex, labeled
  ``O'' at the centroid of the tetrahedron.  Add additional edges
  (shown as dash-dot line segments) that connect ``O'' to the
  centroids of each face, and add six additional faces that include ``O''
  as a vertex.  In Fig.~\ref{f:CubedTetrahedron3} the ``distorted''
  cubes that make up the tetrahedron are illustrated.  The two cubes
  adjacent to vertexes ``A'' and ``C'' are shown with opaque shaded
  faces, while the faces of the cubes adjacent to ``B'' and ``D'' are
  transparent.
  \label{f:CubedTetrahedron}} 
\end{figure}

\begin{figure}
  \begin{center} 
    \vspace{-0.2cm}
    \includegraphics[height=0.35\textwidth]{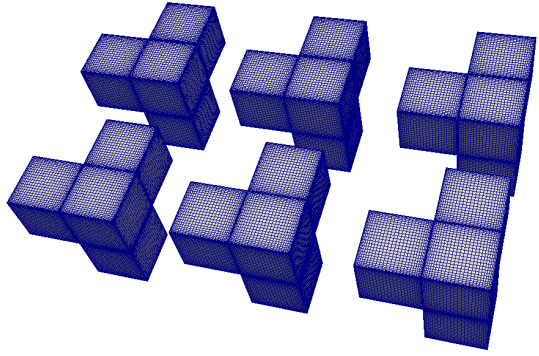}
    \vspace{-0.2cm}
  \end{center}
  \caption{Locations of the cubic regions in $\Rn$ assigned by our
    code for the manifold SFS[RP2/n2:(2,1)(2,-1)] based on the
    triangulation given in the Regina~\cite{Regina} catalog.  Each
    tetrahedron is divided into four cubes, which are placed in groups
    with some of the identified internal faces
    overlapping.\label{f:CubeLocations}}
\end{figure}

Given a three-dimensional triangulation, it is straightforward to
convert it to a multicube structure following the method described in
Ref.~\cite{Lindblom2013}.  The idea is to cut each tetrahedron into
four cubes by adding vertices and edges as illustrated in
Fig.~\ref{f:CubedTetrahedron}, and described in some detail in the
caption. Our code creates a list of cubic regions from the list of
tetrahedrons, then it assigns unique locations in $\Rn$ to each cube.
These locations are chosen so the four cubes associated with each
tetrahedron are grouped together, and these tetrahedron based groups
are arranged in a 2D lattice for convenience of 3D visualization.
Figure~\ref{f:CubeLocations} illustrates the locations of the cubes
assigned by our code for the multicube structure constructed for the
SFS[RP2/n2:(2,1)(2,-1)] manifold from the triangulation given in the
Regina catalog.

Finally our code constructs the appropriate maps in $\Rn$ between cube
faces using Eq.~(\ref{e:CoordinateMap}), following the prescription
given in Ref.~\cite{Lindblom2014}.  In addition to the locations of
each cube, these maps depend on knowing the appropriate
rotation/reflection matrix, ${\mathbf C}_{B\beta}^{A\alpha}$, that
aligns the faces $A\{\alpha\}$ and $B\{\beta\}$ in the appropriate
way.  Each cube has six faces, three of which correspond to internal
connections between the four cubes that make up a tetrahedron.  The
rotation/reflection matrices needed for these internal face
transformations are the same for every tetrahedron group of cubes.  So
they are easily included in the code.  The three additional faces of
each cube are parts of the faces of the tetrahedra.  The appropriate
rotation/reflection matrices for those faces depend on the face
mappings of the triangulations.  There are, however, a reasonably
small number of ways the faces of any two tetrahedra can be
identified.  Our code includes the appropriate matrices for all the
possible cube face matchings (which we determined by systematically
reproducing each possibility with a collection of paper models).  Once
a triangulation with its face mappings is read into our code, it
automatically determines the appropriate cube mappings from its table
of all possibilities.  Our code can then output the complete multicube
structure in any desired format.  For example Tables~ \ref{t:L(5,2)},
\ref{t:SFS[S2:(2,1)(2,1)(2,-1)]}, \ref{t:KB} and \ref{t:SFSRP2} in
\ref{s:3DMulticubeManifolds} are output from our code in \LaTeX
format.  Our code also generates the appropriately formatted input
files used by the SpEC code to compute the reference metrics $\tilde
g_{ab}$ described in Sec.~\ref{s:NumericalExamples}.

Our code can be used to construct a multicube structure from any
three-dimensional triangulation.  However, the methods for
constructing reference metrics presented in
Sec.~\ref{s:C0ReferenceMetrics} only work for special multicube
structures that allow uniform dihedral angles around each edge.  The
code therefore tests several identities to determine when this is
possible.

Once the multicube structure has been constructed by the code, it
determines the dihedral angles $\psi_{A\{\alpha\beta\}}$ around each
edge using the uniform dihedral angle assumption given in
Eq.~(\ref{e:DihedralAngleDef}) of Sec.~\ref{s:Step2}.  The most
important identity that must be satisfied by these
$\psi_{A\{\alpha\beta\}}$ involves the associated angles
$\theta_{A\{\alpha\}}$ between the axes that define the edges of the
$A\{\alpha\}$ face.  These angles must agree with the angles
$\theta_{B\{\alpha\}}$ between the axes of the $B\{\alpha\}$ face
identified with it in the multicube structure.  Without this condition
the intrinsic metric of region $\mathcal{B}_A$ would not be continuous
across that face with the intrinsic metric of region $\mathcal{B}_B$.
The edge angle $\theta_{A\{\alpha\}}$ is related to the dihedral
angles $\psi_{A\{\alpha\beta\}}$ using the spherical law of cosines,
\begin{equation}
  \cos \theta_{A\{\alpha\}} = \frac{\cos\psi_{A\{\beta\gamma\}} +
    \cos\psi_{A\{\alpha\beta\}}\cos\psi_{A\{\alpha\gamma\}}}{\sin\psi_{A\{\alpha\beta\}}
    \sin\psi_{A\{\alpha\gamma\}}}.
\end{equation}
Our code evaluates the $\theta_{A\{\alpha\}}$ for each vertex of each
cube face and determines whether it agrees with the corresponding
angles $\theta_{B\{\alpha\}}$ at those vertices.  Multicube
structures that do not satisfy this condition could not be used in the
present study.

Our code also checks two other less restrictive identities.  One
ensures that the determinant of the flat inverse metric
$e^{ab}_{A\{\alpha\beta\gamma\}}$ defined in
Eq.~(\ref{e:GlobalFlatMetricS2}) is positive in each cubic region:
\begin{equation}
  \det e^{ab}_{A\{\alpha\beta\gamma\}} = 1
  + 2 \cos\psi_{A\{\alpha\beta\}}\cos\psi_{A\{\alpha\gamma\}}\cos\psi_{A\{\beta\gamma\}}
  - \cos\psi_{A\{\alpha\beta\}}^2 - \cos\psi_{A\{\alpha\gamma\}}^2
  - \cos\psi_{A\{\beta\gamma\}}^2 > 0. 
\end{equation}
A second identity ensures that the areas of the spherical triangles
created by the intersection of each cubic region with small spheres
located at their vertices (see Fig.~\ref{f:ThreeDWedgeB}) are positive.
This requires
\begin{equation}
  \psi_{A\{\alpha\beta\}} + \psi_{A\{\alpha\gamma\}} + \psi_{A\{\beta\gamma\}} > \pi.
\end{equation}

We have run this code on all the triangulations consisting of up to
eleven tetrahedra listed in the catalogs of all closed prime orientable
three-manifolds in Refs.~\cite{Matveev2005, Martelli2001,
  Martelli2006, Burton2011, Regina}.  We find that of these only the
$29$ manifolds listed in Table~\ref{t:manifold_list} satisfy all these
constraints.

\section{Solving the Biharmonic Equation Using Pseudo-Spectral Methods}
\label{s:BiharmonicMethods}

The biharmonic equations in two and three dimensions are given by
\begin{eqnarray}
0&=&\left(\partial_x^{\,4} + 2\, \partial_x^{\,2}\,\partial_y^{\,2}
+\partial_y^{\,4}\right)\,U,
\label{e:Biharmonic2Du}\\
0&=&   \left(\partial_x^{\,4} +\partial_y^{\,4} +\partial_z^{\,4}
  +2\,\partial_x^{\,2}\,\partial_y^{\,2}
  +2\,\partial_x^{\,2}\,\partial_z^{\,2} 
  +2\,\partial_y^{\,2}\,\partial_z^{\,2}\right)\,U. 
  \label{e:Biharmonicbarbar3Du}
\end{eqnarray}
The solutions to these equations on compact domains are determined
uniquely by the values of $U$ and its normal derivative $\mathrm{d}U/
\mathrm{d}n$ (the Dirichlet and Neumann conditions respectively) on
the boundaries of the domain~\cite{Barton2014}.

In this study these equations are solved using pseudo-spectral
numerical methods.  A function $U$ is specified in this approach by
its values on a special mesh.  The mesh points used here are located
at the Gauss-Lobatto collocation points~\cite{Boyd2000}.  This choice
makes it possible to transform easily and exactly back and forth
between the mesh representation and a Chebyshev polynomial based
spectral representation of $U$.  The value of $U$ at a particular mesh
point is written here as $U_{\{ij\kern 0.07em\}}$ in two dimensions or
$U_{\{ijk\}}$ in three.  Partial derivatives of a function, which 
are exact for this spectral representation, can be written as special
linear combinations of its values on these mesh points,
\begin{eqnarray}
  \partial_xU_{\{ij\kern 0.07em\}} &=&
  \mathcal{D}^{(x)}{}_i{}^s\,U_{\{sj\kern 0.07em\}},\\
  \partial_yU_{\{ij\kern 0.07em\}} &=&
  \mathcal{D}^{(y)}{}_j{}^t\,U_{\{it\}},
\end{eqnarray}
where the repeated indices $s$ or $t$ are summed over all the mesh
points in the particular direction.  The discrete pseudo-spectral
representation of the two-dimensional biharmonic equation can
therefore be written as
\begin{eqnarray}
  0&=& \mathcal{D}^{(x)}{}_i{}^s\mathcal{D}^{(x)}{}_s{}^t
  \mathcal{D}^{(x)}{}_t{}^u\mathcal{D}^{(x)}{}_u{}^v\,U_{\{vj\kern 0.07em\}}
  +2\,\mathcal{D}^{(x)}{}_i{}^s\mathcal{D}^{(x)}{}_s{}^t
  \mathcal{D}^{(y)}{}_j{}^u\mathcal{D}^{(y)}{}_u{}^v\,U_{\{tv\}}\nonumber\\
  &&+\mathcal{D}^{(y)}{}_j{}^s\mathcal{D}^{(y)}{}_s{}^t
  \mathcal{D}^{(y)}{}_t{}^u\mathcal{D}^{(y)}{}_u{}^v\,U_{\{iv\}}.
  \label{e:DiscreteBiharmonic2D}
\end{eqnarray}
An analogous expression is used for the discrete representation of the
biharmonic equation in three dimensions.

\begin{figure}[!htb]
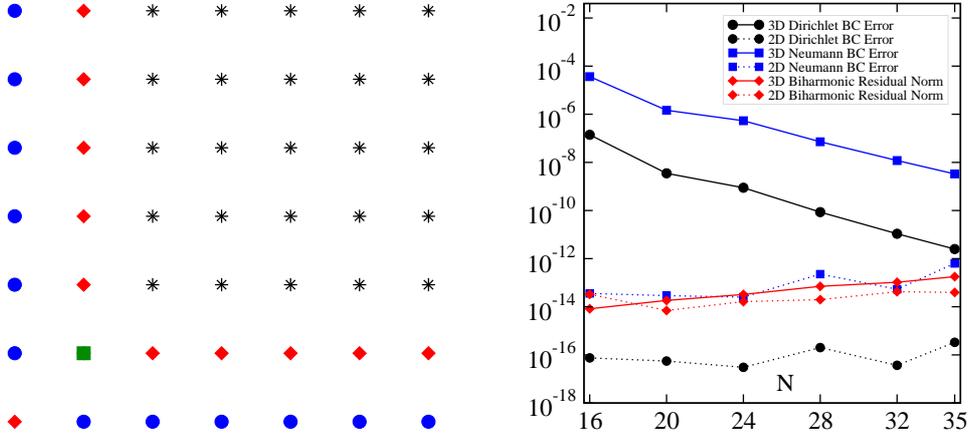
 
  \begin{center}
  \subfigure{
    \includegraphics[height=0.42\textwidth]{FigB1a.eps}
    \label{f:grid}
  }\hspace{1.0cm}\vspace{-0.3cm}
  \subfigure{
    \includegraphics[width=0.43\textwidth]
                    {FigB1b.eps}
                    \label{f:residual}
  }
\end{center}
\caption{Figure on the left represents one corner of a two-dimensional
  mesh used to solve the biharmonic equation.  Discrete
  representations of the boundary conditions for points along the
  boundaries replace the biharmonic equation at the points marked with
  (red) diamonds for Dirichlet and (blue) circles for Neumann
  conditions, respectively.  The average of the Dirichlet conditions
  from both nearby boundary points are imposed at the point marked
  with a (green) square. Discrete representations of the biharmonic
  equation are imposed at the remaining interior mesh points marked
  with (black) stars.  Figure on the right illustrates the average
  (rms) boundary errors in the Dirichlet and Neumann boundary
  conditions for examples of the numerical 2D and 3D biharmonic
  solutions used in this study.}
\end{figure}

Boundary conditions are imposed by replacing the discrete biharmonic
equations along the outer layer of mesh points with discrete versions
of the Neumann boundary conditions at those points.  Dirichlet
boundary conditions are also needed along the boundaries, and those
are imposed by replacing the discrete biharmonic equation on the mesh
points at the next layer of points adjacent to the boundary with the
Dirichlet condition evaluated at the nearest boundary points.
Figure~\ref{f:grid} illustrates where these boundary conditions are
imposed for the case of a two-dimensional mesh.  The three-dimensional
case is analogous, but more difficult to illustrate in two-dimensional
figures.

The functions $U_{\{ij\kern 0.07em\}}$ in two dimensions (or
$U_{\{ijk\}}$ in three) can be thought of as vectors $U_\mathcal{A}$
on a space where the super-index $\mathcal{\tiny{A}}$ ranges over all
the mesh points, i.e. $\mathcal{\tiny{A}}=\{ij\kern 0.07em\}$ in two
dimensions (or $\mathcal{\tiny{A}}=\{ijk\}$ in three).  The discrete
biharmonic equation can be thought of as a linear matrix equation on
this space:
\begin{eqnarray}
  \sum_\mathcal{B}\mathcal{O_A{}^B} \,U_\mathcal{B}=h_\mathcal{\tiny{A}},
  \label{e:MeshBiharEq}
\end{eqnarray} 
where $\mathcal{O_A{}^B}$ is defined in two dimensions by the
expression in Eq.~(\ref{e:DiscreteBiharmonic2D}) at the interior mesh
points. The discrete versions of the Dirichlet and Neumann conditions
are imposed on the components of this equation representing the
surface layers of the mesh. The vector $h_\mathcal{\tiny{A}}$ holds
the boundary data for those conditions, in addition (if any) to the
inhomogeneous source for the equation at the interior points.  The
expression used here for $\mathcal{O_A{}^B}$ in three dimensions is
completely analogous.

Our primary interest is finding smooth functions $U$ that satisfy the
boundary conditions as accurately as possible. The components of the
matrices $\mathcal{D}^{(x)}{}_i{}^j$, etc., which provide discrete
representations of the derivative operators, have average magnitudes
that scale like $N$, where $N$ is the number of mesh points used in
each direction.  Therefore the components of the matrix
$\mathcal{O_A{}^B}$ representing the biharmonic operator on interior
mesh points will scale like $N^{4}$, and for large $N$ will therefore
dominate the boundary condition terms.  These interior components have
therefore been scaled in this study by $N^{-4}$ to emphasize the
relative importance of the boundary conditions.  A similar scaling
would also be applied to any source terms in $h_\mathcal{\tiny{A}}$,
however, no additional scaling is needed for the homogeneous equations
considered here.

The linear equations given in Eq.~(\ref{e:MeshBiharEq}) can be solved
numerically using a variety of iterative techniques, e.g. using
solvers such as GMRES~\cite{SAAD1986} or
BI-CGSTAB~\cite{VanDerVorst1992}.  Numerical experiments using
pseudo-spectral methods described above for this problem showed that
faster and more accurate results could be obtained using more direct
non-iterative methods, because the meshes used here are relatively
small (in comparison with those used by standard finite difference or
finite element methods).  The matrix $\mathcal{O_A{}^B}$ has size
$N^2\times N^2$ for the two-dimensional problem and $N^3\times N^3$
for three, where $N$ is the number of mesh points used in each
dimension.  The largest meshes used in this study have $N=35$, so the
largest matrix has size $1,225\times 1,225$ for the two-dimensional
meshes, and $42,875\times 42,875$ for three. For matrices of this
size, it is possible to construct the LU decomposition of
$\mathcal{O_A{}^B}$ directly using modest computing resources.  Very
fast direct algorithms then exist for solving such linear systems
exactly, see e.g. Ref.~\cite{Press1992}.  The construction of the LU
decomposition requires a lot of memory and computer time.  The highest
resolution that could be run on the computing facility available to us
is $N=35$ due to memory limitations.  Constructing the LU
decomposition at this resolution required about 152 hours on a single
processor.  But once computed for each needed resolution $N$, these LU
decompositions can be stored on disk and quickly read in whenever they
are needed.  A very accurate solution of the linear equations in LU
form can then be obtained very quickly and efficiently.  Pre-computing
the LU decompositions in this way reduces the $N=35$ problem of
solving one 3D biharmonic equation (plus six 2D biharmonic equations
to set the boundary conditions) from about 152 hours to about 75
seconds.

\begin{figure}[!htb] 
  \begin{center}
    \includegraphics[height=0.42\textwidth]{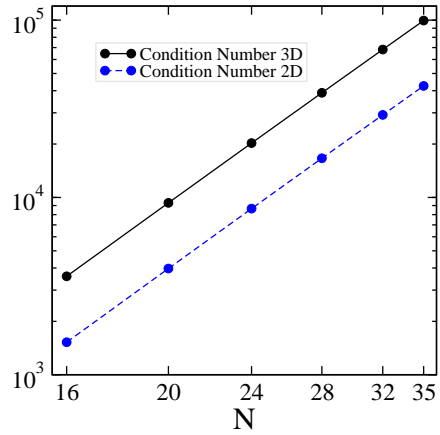}
  \end{center}
  \caption{Figure illustrates the spatial resolution, $N$, dependence
    of the condition number $\kappa_\infty$ (using the estimate given
    in Eq.~\ref{e:ConditionNumber}) for the pseudo-spectral matrix
    representations of the two- and three-dimensional biharmonic
    operators used in this study.\label{f:CN} }
\end{figure}

The condition number $\kappa$ of a matrix operator $\mathcal{O_A{}^B}$
is a measure of how accurately linear equations like
Eq.~(\ref{e:MeshBiharEq}) can be solved numerically~\cite{Cline1979}.
Fractional errors in the solutions $U_\mathcal{B}$ can be as large as
$\kappa$ multiplied by the fractional errors in the matrix
$\mathcal{O_A{}^B}$ and the source $h_\mathcal{A}$. We have estimated
$\kappa$ for the two and three-dimensional representations of the
biharmonic matrices used in our study.  These approximations were
obtained using the simple approximate expression:
\begin{equation}
  \kappa_\infty = ||\mathcal{O}||_\infty \, ||\mathcal{O}^{-1}||_\infty
  \geq \frac{\max_\mathcal{A}\sum_\mathcal{B}|\mathcal{O_A{}^B}|}
       {\min_\mathcal{A}|\mathcal{U_A{}^A}|},\label{e:ConditionNumber}
\end{equation}
where $\mathcal{U_A{}^B}$ is the upper diagonal part of the LU
decomposition: $\sum_\mathcal{C}\mathcal{P}_A{}^C \mathcal{O_C{}^B}
=\sum_\mathcal{C}\mathcal{L_A{}^CU_C{}^B}$, where $\mathcal{P}_A{}^C$
is a permutation matrix.  These estimates of $\kappa_\infty$ using
Eq.~\ref{e:ConditionNumber} are illustrated in Fig.~\ref{f:CN}.  As
the figures shows, these estimates for the condition number scale with
spatial resolution $N$ as a power law: $\sim\! N^{\,4.2}$.  While
condition numbers as large as the $10^5$ seen in this figure might
seem large, they mean that fractional erors in the matrix
$\mathcal{O_A{}^B}$ and the source $h_\mathcal{A}$ at the double
precision roundoff level, $10^{-16}$, could produce fractional errors
only as large as $10^{-11}$ in the solution $U_\mathcal{B}$.  We also
note that condition numbers of this size have not significantly
influenced the boundary condition errors (our primary interest in
these solutions), as illustrated in Figs.~\ref{f:residual},
\ref{f:SurfaceError1} and \ref{f:SurfaceError2}.

In addition to being very quick and efficient, the direct LU solver
method used here provides solutions having better accuracy for our
purposes than those obtained with the iterative solvers that were
tested.  Solutions to the two- and three-dimensional biharmonic
equation are used here at various stages in the construction of a
reference metric.  The important requirement is the need to have those
solutions satisfy the Dirichlet and Neumann boundary conditions as
accurately as possible.  The interior details are not of primary
importance, so long as they are smooth.  Figure~\ref{f:residual}
illustrates the convergence with resolution $N$ of the errors in the
Dirichlet and Neumann boundary conditions satisfied by numerical
examples of two- and three-dimensional biharmonic solutions obtained
with this direct LU solver method.  The two-dimensional results are at
the double-precision roundoff levels for all values of $N$ tested,
while the three-dimensional results show the exponential convergence
that is expected for pseudo-spectral methods.  The average interior
bulk residual errors are also roughly at double-precision roundoff
levels.  The boundary condition accuracies achieved using this direct
LU solver method were much better than anything obtained with the
iterative solvers tested here.

\section{Proof of the Identity $\bar N^{\{\beta\}}\,
  \bar K^{\{\beta\}}_{\alpha\gamma} =-\bar
N^{\{\alpha\}}\, \bar K^{\{\alpha\}}_{\beta\gamma}$}
\label{s:IdentityProof}

The following simple argument shows that this condition is satisfied
along the $A\{\alpha\beta\}$ edge by the $C^0$ metric $\bar g_{ab}$
constructed in Sec.~\ref{s:Step3.1}.  Start with the identity
\begin{equation}
  \gamma^a\,\bar n^{\{\alpha\}\,b}\bar K^{\{\beta\}}_{ab}
  = \gamma^a\,\bar n^{\{\alpha\}\,b}\,\bar \nabla_a \bar n^{\{\beta\}}_b
  = \gamma^a\,\bar \nabla_a\left(\bar n^{\{\alpha\}\,b}\,\bar n^{\{\beta\}}_b\right)
     -\gamma^a\,\bar n^{\{\beta\}\,b}\bar \nabla_a \,\bar n^{\{\alpha\}}_b
     =-\gamma^a\,\bar n^{\{\beta\}\,b}\bar K^{\{\alpha\}}_{ab},
     \label{e:KbarTensorIdentity}
\end{equation}
where the $\gamma^a$ are the components of the vector $\partial_\gamma
= \gamma^a\partial_a$ that is tangent to both the $A\{\alpha\}$ and
$A\{\beta\}$ faces.  The last equality follows from the fact that
$\partial_\gamma(\bar n^{\{\alpha\}\,b}\,\bar n^{\{\beta\}}_b)=0$
because the dihedral angle is constant along the $A\{\alpha\beta\}$
edge.  The additional simple identities $ \bar
n^{\{\alpha\}\,a}\,\gamma^b\,\bar K^{\{\alpha\}}_{ab}= \bar
n^{\{\beta\}\,a}\,\gamma^b\,\bar K^{\{\beta\}}_{ab}=0$ and $\bar
K^{\{\alpha\}}_{\gamma\gamma}=\bar K^{\{\beta\}}_{\gamma\gamma}=0$ can
be used to transform the tensor identity in
Eq.~(\ref{e:KbarTensorIdentity}) into the coordinate identity given in
Eq.~(\ref{e:KbarCoordinateIdentity}).  First obtain the coordinate
representations of the simple identities:
\begin{eqnarray}
  0&=&  \bar n^{\{\alpha\}\,a}\,\gamma^b\, \bar K^{\{\alpha\}}_{ab}
  = \bar N^{\{\alpha\}}\left( \bar g^{\alpha\alpha} \bar K^{\{\alpha\}}_{\alpha \gamma}
  +  \bar g^{\alpha\beta} \bar K^{\{\alpha\}}_{\beta \gamma}\right),
  \label{e:Kcords1}\\
  0&=&  \bar n^{\{\beta\}\,a}\,\gamma^b\, \bar K^{\{\beta\}}_{ab}
  = \bar N^{\{\beta\}}\left( \bar g^{\beta\alpha} \bar K^{\{\beta\}}_{\alpha \gamma}
  +  \bar g^{\beta\beta} \bar K^{\{\beta\}}_{\beta \gamma}\right).\label{e:Kcords2}
\end{eqnarray}
Coordinate representations of $\bar n^{\{\beta\}\,a}\kern 0.05em
\gamma^b\bar K^{\{\alpha\}}_{ab}$ and $\bar n^{\{\alpha\}\,a}\kern
0.05em \gamma^b\bar K^{\{\beta\}}_{ab}$ can be written as
\begin{eqnarray}
    \bar n^{\{\beta\}\,a}\kern 0.05em \gamma^b\bar K^{\{\alpha\}}_{ab} &=& 
 \bar N^{\{\beta\}}\left( \bar g^{\beta\alpha} \bar K^{\{\alpha\}}_{\alpha \gamma}
 +  \bar g^{\beta\beta} \bar K^{\{\alpha\}}_{\beta \gamma}\right),
 \label{e:Kcords3}\\
 \bar n^{\{\alpha\}\,a}\kern 0.05em \gamma^b\bar K^{\{\beta\}}_{ab} &=&
   \bar N^{\{\alpha\}}\left( \bar g^{\alpha\alpha} \bar K^{\{\beta\}}_{\alpha \gamma}
   +  \bar g^{\alpha\beta} \bar K^{\{\beta\}}_{\beta \gamma}\right).
   \label{e:Kcords4}
 \end{eqnarray}
These expressions can be simplified by using Eq.~(\ref{e:Kcords1}) to
express $\bar K^{\{\alpha\}}_{\alpha\gamma}$ in terms of $\bar
K^{\{\alpha\}}_{\beta\gamma}$, and similarly Eq.~(\ref{e:Kcords2}) to
express $\bar K^{\{\beta\}}_{\beta\gamma}$ in terms of $\bar
K^{\{\beta\}}_{\alpha\gamma}$.  Making these
substitutions in Eqs.~(\ref{e:Kcords3}) and (\ref{e:Kcords4}) gives
\begin{eqnarray}
        \bar n^{\{\beta\}\,a}\kern 0.05em \gamma^b\bar K^{\{\alpha\}}_{ab} &=&
        \bar N^{\{\beta\}}\left(\bar N^{\{\alpha\}}\right)^2
        \left[\bar g^{\alpha\alpha}\bar g^{\beta\beta} -(\bar
  g^{\alpha\beta})^2\right]\bar K^{\{\alpha\}}_{\beta\gamma},\\
    \bar n^{\{\alpha\}\,a}\kern 0.05em \gamma^b\bar K^{\{\beta\}}_{ab} &=&
    \bar N^{\{\alpha\}}\left(\bar N^{\{\beta\}}\right)^2
    \left[\bar g^{\alpha\alpha}\bar g^{\beta\beta} -(\bar
      g^{\alpha\beta})^2\right]\bar K^{\{\beta\}}_{\alpha\gamma}.
\end{eqnarray}
It follows that the identity $\bar n^{\{\alpha\}\,a}\kern 0.05em
\gamma^b\bar K^{\{\beta\}}_{ab}= -\bar n^{\{\beta\}\,a}\kern 0.05em
\gamma^b\bar K^{\{\alpha\}}_{ab}$ from
Eq.~(\ref{e:KbarTensorIdentity}) implies the identity $\bar
N^{\{\beta\}}\,\bar K^{\{\beta\}}_{\alpha\gamma} =-\bar
N^{\{\alpha\}}\,\bar K^{\{\alpha\}}_{\beta\gamma}$ given in
Eq.~(\ref{e:KbarCoordinateIdentity}).  
  
\section{Example Three-Dimensional Multicube Manifolds}
\label{s:3DMulticubeManifolds}
This appendix gives detailed descriptions of the multicube structures
of several manifolds used in this study.  New structures are presented
here for all the examples constructed by hand (except the trivial flat
examples, E2 and E3): the Poincar\'e dodecahedral
space~\cite{ThrelfallSeifert1931}, Seifert-Weber
space~\cite{SeifertWeber1933}, G2$\times$S1, and the three non-trivial
compact orientable three-manifolds that admit flat
metrics~\cite{riazuelo2004, Hitchman2018}, E4, E5 and E6
(Hantzsche-Wendt space~\cite{Hantzsche1934}).  In addition, a
selection of the multicube structures constructed automatically from
triangulations using the code described in
\ref{s:TriangulationToMultiCubeCode} are presented here:
KB/n2$\times\!\!\sim$S1, L(5,2), SFS[RP2/n2:(2,1)(2,-1)], and
SFS[S2:(2,1)(2,1)(2,-1)].

The notation used to describe these multicube structures is based on
that introduced in Ref.~\cite{Lindblom2013}.  Each multicube
structure consists of a set of non-overlapping cubes, ${\cal B}_A$,
that cover the manifold, and a set of maps $\Psi^{A\alpha}_{B\beta}$
that identify the faces of neighboring cubes.  The interface boundary
maps used here (written in terms of the global Cartesian coordinates
used for the multicube structure) take points, $x^i_B$, on the
interface boundary $B\{\beta\}$ (or equivalently
$\partial_\beta\mathcal{B}_B$) of region $\mathcal{B}_B$ to the
corresponding points, $x^i_A$, in the boundary $A\{\alpha\}$ (or
equivalently $\partial_\alpha\mathcal{B}_A$) of region $\mathcal{B}_A$
in the following way,
\begin{eqnarray}
x^i_A = c^i_A +f^i_{\alpha} + C_{B\beta\,j}^{A\alpha\, i}(x^j_B
 - c^j_B-f^j_{\beta}).
\label{e:CoordinateMap}
\end{eqnarray}
The vectors $\vec c_A+\vec f_{\alpha}$ and $\vec c_B+\vec f_{\beta}$
are the locations of the centers of the $A\{\alpha\}$ and $B\{\beta\}$
faces respectively, and ${\mathbf C}_{B\beta}^{A\alpha}$ is the
combined rotation/reflection matrix needed to orient the faces
properly.

The following tables include lists of the cubic regions,
$\mathcal{B}_A$, used to cover the manifold in each structure, the
vectors $\vec c_A$ that define the locations of the centers of these
regions in $\mathbb{R}^3$, and the rotation/reflection matrices
${\mathbf C}_{B\beta}^{A\alpha}$ needed to transform each cube face
into the face of its neighbor.\footnote{The vectors $\vec f_\alpha$
are the relative positions of the center of the $A\{\alpha\}$ cube
face with the center of region $\mathcal{B}_A$.  These vectors are the
same for all the cubic regions, and are given explicitly in
Ref.~\cite{Lindblom2013} so they are not repeated here.}  The
identification of the $B\{\beta\}$ face with the $A\{\alpha\}$ face is
indicated in the tables by $\{\alpha A\}\leftrightarrow \{\beta B\}$.
The notation $\mathbf{I}$ in these tables indicates the identity
matrix, while $\mathbf{R}_{\alpha}$ indicates the $+\pi/2$ rotation
about the outward directed normal to the $\{\alpha\}$ cube face.



\begin{table}[!htb]
  \scriptsize
\renewcommand{\arraystretch}{1.5}
\begin{center}
  \caption{ Multi-Cube representation of Third-Turn
    space~\cite{riazuelo2004, Hitchman2018} (E4, one of the six
    compact orientable three-manifolds that admits a global flat
    metric), can be constructed by identifying opposite rectangular
    faces of a hexagonal cylinder, and identifying the two hexagonal
    faces after twisting by $2\pi/3$.  Multicube structure: region
    center locations $\vec c_A$, region face identifications,
    $\{\alpha \,A\} \leftrightarrow \{\beta\, B\}$ , and the rotation
    matrices for the associated interface maps, ${\bf
      C}_{A\alpha}^{B\beta}$.}
    \label{t:ThirdTurnSpaceE4}
\begin{tabular}{c|c|c|c|c|c|c|c}
  \toprule && $\alpha=-x$ & $\alpha=+x$ & $\alpha=-y$ & $\alpha=+y$ &
  $\alpha=-z$ & $\alpha=+z$ \\ A & $\vec c_A$
  &$B\,\,\beta\,\,\,{\mathbf C}_{A\alpha}^{B\beta}$
  &$B\,\,\beta\,\,\,{\mathbf C}_{A\alpha}^{B\beta}$
  &$B\,\,\beta\,\,\,{\mathbf C}_{A\alpha}^{B\beta}$
  &$B\,\,\beta\,\,\,{\mathbf C}_{A\alpha}^{B\beta}$
  &$B\,\,\beta\,\,\,{\mathbf C}_{A\alpha}^{B\beta}$
  &$B\,\,\beta\,\,\,{\mathbf C}_{A\alpha}^{B\beta}$\\ \midrule $1$ &
  $( 0, 0, 0)$ & $2+x\,\,\,\mathbf{I}$ & $2-x\,\,\,\mathbf{I}$ &
  $3+y\,\,\,\mathbf{I}$ & $3-y\,\,\,\mathbf{I}$ &
  $2+z\,\,\,\mathbf{R}_{+z}$ & $3-z\,\,\,\mathbf{R}_{-z}$ \\ $2$ & $(
  1, 0, 0)$ & $1+x\,\,\,\mathbf{I}$ & $1-x\,\,\,\mathbf{I}$ &
  $3-x\,\,\,\mathbf{R}_{+z}$ & $3+x\,\,\,\mathbf{R}_{+z}$ &
  $3+z\,\,\,\mathbf{R}_{+z}^2$ & $1-z\,\,\,\mathbf{R}_{-z}$ \\ $3$ &
  $( 0, 1, 0)$ & $2-y\,\,\,\mathbf{R}_{-z}$ &
  $2+y\,\,\,\mathbf{R}_{-z}$ & $1+y\,\,\,\mathbf{I}$ &
  $1-y\,\,\,\mathbf{I}$ & $1+z\,\,\,\mathbf{R}_{+z}$ &
  $2-z\,\,\,\mathbf{R}_{-z}^2$ \\ \bottomrule
 \end{tabular}
\end{center}
\end{table}


\begin{table}[!htb]
  \scriptsize
\renewcommand{\arraystretch}{1.5}
\begin{center}
  \caption{
      Multi-Cube representation of Sixth-Turn
    space~\cite{riazuelo2004, Hitchman2018} (E5, one of the
    six compact orientable three-manifolds that admits a global flat
    metric), can be constructed by identifying opposite rectangular
    faces of a hexagonal cylinder, and identifying the two hexagonal
    faces after twisting by $2\pi/6$.  Multicube structure: region
    center locations $\vec c_A$, region face identifications,
    $\{\alpha \,A\} \leftrightarrow \{\beta\, B\}$ , and the rotation
    matrices for the associated interface maps, ${\bf
      C}_{A\alpha}^{B\beta}$.}
    \label{t:TableSixthTurnSpaceE5}
\begin{tabular}{c|c|c|c|c|c|c|c}
  \toprule && $\alpha=-x$ & $\alpha=+x$ & $\alpha=-y$ & $\alpha=+y$ &
  $\alpha=-z$ & $\alpha=+z$ \\ A & $\vec c_A$
  &$B\,\,\beta\,\,\,{\mathbf C}_{A\alpha}^{B\beta}$
  &$B\,\,\beta\,\,\,{\mathbf C}_{A\alpha}^{B\beta}$
  &$B\,\,\beta\,\,\,{\mathbf C}_{A\alpha}^{B\beta}$
  &$B\,\,\beta\,\,\,{\mathbf C}_{A\alpha}^{B\beta}$
  &$B\,\,\beta\,\,\,{\mathbf C}_{A\alpha}^{B\beta}$
  &$B\,\,\beta\,\,\,{\mathbf C}_{A\alpha}^{B\beta}$\\ \midrule $1$ &
  $( 0, 1, 0)$ & $5-x\,\,\,\mathbf{R}^2_{+z}$ &
  $6+y\,\,\,\mathbf{R}_{-z}$ & $2+y\,\,\,\mathbf{I}$ &
  $3+x\,\,\,\mathbf{R}_{+z}$ & $2+z\,\,\,\mathbf{R}_{+z}$ &
  $6-z\,\,\,\mathbf{R}^2_{+z}$ \\ $2$ & $( 0, 0, 0)$ &
  $4-x\,\,\,\mathbf{R}^2_{+z}$ & $3-x\,\,\,\mathbf{I}$ &
  $6-y\,\,\,\mathbf{R}^2_{+z}$ & $1-y\,\,\,\mathbf{I}$ &
  $3+z\,\,\,\mathbf{R}_{+z}$ & $1-z\,\,\,\mathbf{R}_{-z}$ \\ $3$ & $(
  1, 0, 0)$ & $2+x\,\,\,\mathbf{I}$ & $1+y\,\,\,\mathbf{R}_{-z}$ &
  $5-y\,\,\,\mathbf{R}^2_{+z}$ & $4+x\,\,\,\mathbf{R}_{+z}$ &
  $4+z\,\,\,\mathbf{R}^2_{+z}$ & $2-z\,\,\,\mathbf{R}_{-z}$ \\ $4$ &
  $( 3, 1, 0)$ & $2-x\,\,\,\mathbf{R}^2_{+z}$ &
  $3+y\,\,\,\mathbf{R}_{-z}$ & $5+y\,\,\,\mathbf{I}$ &
  $6+x\,\,\,\mathbf{R}_{+z}$ & $5+z\,\,\,\mathbf{R}_{+z}$ &
  $3-z\,\,\,\mathbf{R}^2_{+z}$ \\ $5$ & $( 3, 0, 0)$ &
  $1-x\,\,\,\mathbf{R}^2_{+z}$ & $6-x\,\,\,\mathbf{I}$ &
  $3-y\,\,\,\mathbf{R}^2_{+z}$ & $4-y\,\,\,\mathbf{I}$ &
  $6+z\,\,\,\mathbf{R}_{+z}$ & $4-z\,\,\,\mathbf{R}_{-z}$ \\ $6$ & $(
  4, 0, 0)$ & $5+x\,\,\,\mathbf{I}$ & $4+y\,\,\,\mathbf{R}_{-z}$ &
  $2-y\,\,\,\mathbf{R}^2_{+z}$ & $1+x\,\,\,\mathbf{R}_{+z}$ &
  $1+z\,\,\,\mathbf{R}^2_{+z}$ & $5-z\,\,\,\mathbf{R}_{-z}$
  \\ \bottomrule
 \end{tabular}
\end{center}
\end{table}


\begin{table}[!htb]
  \scriptsize
\renewcommand{\arraystretch}{1.5}
\begin{center}
  \caption{ Multi-Cube representation of Hantzsche-Wendt
      space~\cite{Hantzsche1934, riazuelo2004, Hitchman2018} (E6,
      one of the six compact orientable three-manifolds that admits a
      global flat metric), can be constructed by identifying faces on
      two cubic regions (see Ref.~\cite{Hitchman2018} example 8.1.7 for
      details) .  Multicube structure: region center locations $\vec
      c_A$, region face identifications, $\{\alpha \,A\} \rightarrow
      \{\beta\, B\}$ , and the rotation matrices for the associated
      interface maps, ${\bf C}_{A\alpha}^{B\beta}$.}
    \label{t:HantzscheWendtSpace}
\begin{tabular}{c|c|c|c|c|c|c|c}
  \toprule
  && $\alpha=-x$ & $\alpha=+x$ & $\alpha=-y$ & $\alpha=+y$ & $\alpha=-z$ &
  $\alpha=+z$ \\
A  & $\vec c_A$
&$B\,\,\beta\,\,\,{\mathbf C}_{A\alpha}^{B\beta}$
&$B\,\,\beta\,\,\,{\mathbf C}_{A\alpha}^{B\beta}$
&$B\,\,\beta\,\,\,{\mathbf C}_{A\alpha}^{B\beta}$
&$B\,\,\beta\,\,\,{\mathbf C}_{A\alpha}^{B\beta}$
&$B\,\,\beta\,\,\,{\mathbf C}_{A\alpha}^{B\beta}$
&$B\,\,\beta\,\,\,{\mathbf C}_{A\alpha}^{B\beta}$\\
\midrule
$0$
& $( 0, 0, 0)$
& $1+x\,\,\,\mathbf{R}_{+x}^2$ & $1-x\,\,\,\mathbf{R}_{+x}^2$
& $1-y\,\,\,\mathbf{R}_{+z}^2$ & $1+y\,\,\,\mathbf{R}_{+z}^2$
& $1+z\,\,\,\mathbf{I}$ & $1-z\,\,\,\mathbf{I}$
\\
$1$
& $( 0, 0, 1)$
& $0+x\,\,\,\mathbf{R}_{+x}^2$ & $0-x\,\,\,\mathbf{R}_{+x}^2$
& $0-y\,\,\,\mathbf{R}_{+z}^2$ & $0+y\,\,\,\mathbf{R}_{+z}^2$
& $0+z\,\,\,\mathbf{I}$ & $0-z\,\,\,\mathbf{I}$
\\
\bottomrule
\end{tabular}
\end{center}
\end{table}

 
\begin{table}[htb]
\scriptsize
  \renewcommand{\arraystretch}{1.5}
\begin{center}
  \caption{Multicube representation of the product space G2$\times$S1
    constructed from the genus number $N_g=2$ two-dimensional compact
    orientable manifold. Multicube Structure: region center locations
    $\vec c_A$, region face identifications, $\{\alpha \,A\}
    \leftrightarrow \{\beta\, B\}$ , and the rotation matrices for the
    associated interface maps, ${\bf C}_{A\alpha}^{B\beta}$.
   \label{t:TableG2xS1}}
\begin{tabular}{c|c|c|c|c|c|c|c}
  \toprule
  && $\alpha=-x$ & $\alpha=+x$ & $\alpha=-y$ & $\alpha=+y$ & $\alpha=-z$ &
  $\alpha=+z$ \\
A  & $\vec c_A$
&$B\,\,\beta\,\,\,{\mathbf C}_{A\alpha}^{B\beta}$
&$B\,\,\beta\,\,\,{\mathbf C}_{A\alpha}^{B\beta}$
&$B\,\,\beta\,\,\,{\mathbf C}_{A\alpha}^{B\beta}$
&$B\,\,\beta\,\,\,{\mathbf C}_{A\alpha}^{B\beta}$
&$B\,\,\beta\,\,\,{\mathbf C}_{A\alpha}^{B\beta}$
&$B\,\,\beta\,\,\,{\mathbf C}_{A\alpha}^{B\beta}$\\
\midrule
$1$
& $(L,2L,0)$
& $8+x\,\,\,\mathbf{I}$ & $10-x\,\,\,\mathbf{I}$
& $2+y\,\,\,\mathbf{I}$ & $4-y\,\,\,\mathbf{I}$
& $1+z\,\,\,\mathbf{I}$ & $1-z\,\,\,\mathbf{I}$
\\
$2$
& $(L,L,0)$
& $7+x\,\,\,\mathbf{I}$ & $4+x\,\,\,\mathbf{R}^2_{+z}$
& $3+y\,\,\,\mathbf{I}$ & $1-y\,\,\,\mathbf{I}$
& $2+z\,\,\,\mathbf{I}$ & $2-z\,\,\,\mathbf{I}$
\\
$3$
& $(L,0,0)$
& $6+x\,\,\,\mathbf{I}$ & $9-x\,\,\,\mathbf{I}$
& $4+y\,\,\,\mathbf{I}$ & $2-y\,\,\,\mathbf{I}$
& $3+z\,\,\,\mathbf{I}$ & $3-z\,\,\,\mathbf{I}$
\\
$4$
& $(L,-L,0)$
& $5+x\,\,\,\mathbf{I}$ & $2+x\,\,\,\mathbf{R}^2_{-z}$
& $1+y\,\,\,\mathbf{I}$ & $3-y\,\,\,\mathbf{I}$
& $4+z\,\,\,\mathbf{I}$ & $4-z\,\,\,\mathbf{I}$
\\
$5$
& $(0,-L,0)$
& $7-x\,\,\,\mathbf{R}^2_{+z}$ & $4-x\,\,\,\mathbf{I}$
& $8+y\,\,\,\mathbf{I}$ & $6-y\,\,\,\mathbf{I}$
& $5+z\,\,\,\mathbf{I}$ & $5-z\,\,\,\mathbf{I}$
\\
$6$
& $( 0, 0, 0)$
& $9+x\,\,\,\mathbf{I}$ & $3-x\,\,\,\mathbf{I}$
& $5+y\,\,\,\mathbf{I}$ & $7-y\,\,\,\mathbf{I}$
& $6+z\,\,\,\mathbf{I}$ & $6-z\,\,\,\mathbf{I}$
\\
$7$
& $(0,L,0)$
& $5-x\,\,\,\mathbf{R}^2_{-z}$ & $2-x\,\,\,\mathbf{I}$
& $6+y\,\,\,\mathbf{I}$ & $8-y\,\,\,\mathbf{I}$
& $7+z\,\,\,\mathbf{I}$ & $7-z\,\,\,\mathbf{I}$
\\
$8$
& $(0,2L,0)$
& $10+x\,\,\,\mathbf{I}$ & $1-x\,\,\,\mathbf{I}$
& $7+y\,\,\,\mathbf{I}$ & $5-y\,\,\,\mathbf{I}$
& $8+z\,\,\,\mathbf{I}$ & $8-z\,\,\,\mathbf{I}$
\\
$9$
& $(-L,0,0)$
& $3+x\,\,\,\mathbf{I}$ & $6-x\,\,\,\mathbf{I}$
& $9+y\,\,\,\mathbf{I}$ & $9-y\,\,\,\mathbf{I}$
& $9+z\,\,\,\mathbf{I}$ & $9-z\,\,\,\mathbf{I}$
\\
$10$
& $(-L,2L,0)$
& $1+x\,\,\,\mathbf{I}$ & $8-x\,\,\,\mathbf{I}$
& $10+y\,\,\,\mathbf{I}$ & $10-y\,\,\,\mathbf{I}$
& $10+z\,\,\,\mathbf{I}$ & $10-z\,\,\,\mathbf{I}$
\\
\bottomrule
 \end{tabular}
\end{center}
\end{table}

 
\begin{table}[htb]
  \scriptsize
\renewcommand{\arraystretch}{1.5}
\begin{center}
\caption{Multicube representation of the Poincar\'e dodecahedral space
  (also called the Poincar\'e homology
  three-sphere)~\cite{ThrelfallSeifert1931}.  This multicube structure
  is based on cutting a dodecahedron into twenty cubes (each vertex of
  the dodecahedron is the vertex of one of the cubes, opposite
  vertices of these cubes all intersect at the center of the
  dodecahedron) and identifying opposite faces of the dodecahedron
  after rotation by $\pi/5$.  Multicube Structure: region center
  locations $\vec c_A$, region face identifications, $\{\alpha \,A\}
  \leftrightarrow \{\beta\, B\}$ , and the rotation matrices for the
  associated interface maps, ${\bf C}_{A\alpha}^{B\beta}$.
   \label{t:TablePoincareS3}}
\begin{tabular}{c|c|c|c|c|c|c|c}
  \toprule
  && $\alpha=-x$ & $\alpha=+x$ & $\alpha=-y$ & $\alpha=+y$ & $\alpha=-z$ &
  $\alpha=+z$ \\
A  & $\vec c_A$
&$B\,\,\beta\,\,\,{\mathbf C}_{A\alpha}^{B\beta}$
&$B\,\,\beta\,\,\,{\mathbf C}_{A\alpha}^{B\beta}$
&$B\,\,\beta\,\,\,{\mathbf C}_{A\alpha}^{B\beta}$
&$B\,\,\beta\,\,\,{\mathbf C}_{A\alpha}^{B\beta}$
&$B\,\,\beta\,\,\,{\mathbf C}_{A\alpha}^{B\beta}$
&$B\,\,\beta\,\,\,{\mathbf C}_{A\alpha}^{B\beta}$\\
\midrule
$1$
& $(2L,3L,0)$
& $12+y\,\,\,\mathbf{R}_{+y}\mathbf{R}_{+z}$ & $15-y\,\,\,\mathbf{R}_{+z}$
& $8-y\,\,\,\mathbf{R}_{-y}\mathbf{R}^2_{+x}$ & $4+x\,\,\,\mathbf{R}_{+x}\mathbf{R}_{+z}$
& $10-y\,\,\,\mathbf{R}_{-y}\mathbf{R}_{+x}$ & $2-z\,\,\,\mathbf{I}$
\\ 
$2$
& $(2L,3L,L)$
& $6+x\,\,\,\mathbf{I}$ & $16-y\,\,\,\mathbf{R}_{+z}$
& $18+y\,\,\,\mathbf{I}$ & $13-z\,\,\,\mathbf{R}_{-z}\mathbf{R}_{+x}$
& $1+z\,\,\,\mathbf{I}$ & $7-z\,\,\,\mathbf{R}_{+z}$
\\
$3$
& $(4L,0,3L)$
& $7+x\,\,\,\mathbf{I}$ & $12+z\,\,\,\mathbf{R}_{+z}\mathbf{R}_{+y}$
& $19+y\,\,\,\mathbf{I}$ & $9-x\,\,\,\mathbf{R}_{-z}$
& $18+z\,\,\,\mathbf{R}_{-z}$ & $4-z\,\,\,\mathbf{I}$
\\
$4$
& $(4L,0,4L)$
& $17-x\,\,\,\mathbf{R}_{+x}\mathbf{R}^2_{+z}$ & $1+y\,\,\,\mathbf{R}_{-z}\mathbf{R}_{-x}$
& $13+x\,\,\,\mathbf{R}_{+x}\mathbf{R}_{-z}$ & $10-x\,\,\,\mathbf{R}_{-z}$
& $3+z\,\,\,\mathbf{I}$ & $15-x\,\,\,\mathbf{R}_{-x}\mathbf{R}_{+y}$
\\
$5$
& $(0,3L,0)$
& $8+y\,\,\,\mathbf{R}_{+y}\mathbf{R}_{+z}$ & $12+x\,\,\,\mathbf{R}_{+x}\mathbf{R}^2_{+y}$
& $16+y\,\,\,\mathbf{R}_{+y}$ & $19+x\,\,\,\mathbf{R}_{+z}$
& $14-y\,\,\,\mathbf{R}_{+x}$ & $6-z\,\,\,\mathbf{I}$
\\
$6$
& $(0,3L,L)$
& $17-z\,\,\,\mathbf{R}_{-z}\mathbf{R}_{+y}$ & $2-x\,\,\,\mathbf{I}$
& $10+y\,\,\,\mathbf{I}$ & $20+x\,\,\,\mathbf{R}_{+z}$
& $5+z\,\,\,\mathbf{I}$ & $11-z\,\,\,\mathbf{R}_{+z}$
\\
$7$
& $(2L,3L,3L)$
& $13-x\,\,\,\mathbf{R}^2_{+z}$ & $3-x\,\,\,\mathbf{I}$
& $11+y\,\,\,\mathbf{I}$ & $16+z\,\,\,\mathbf{R}_{-x}$
& $2+z\,\,\,\mathbf{R}_{-z}$ & $8-z\,\,\,\mathbf{I}$
\\
$8$
& $(2L,3L,4L)$
& $14-x\,\,\,\mathbf{R}^2_{+z}$ & $17+y\,\,\,\mathbf{R}_{+y}\mathbf{R}_{-z}$
& $1-y\,\,\,\mathbf{R}^2_{-x}\mathbf{R}_{+y}$ & $5-x\,\,\,\mathbf{R}_{-z}\mathbf{R}_{-y}$
& $7+z\,\,\,\mathbf{I}$ & $19-y\,\,\,\mathbf{R}_{-y}\mathbf{R}_{-x}$
\\
$9$
& $(0,L,0)$
& $3+y\,\,\,\mathbf{R}_{+z}$ & $20-x\,\,\,\mathbf{R}_{-x}$
& $12-x\,\,\,\mathbf{R}_{-x}\mathbf{R}_{+z}$ & $16+x\,\,\,\mathbf{R}_{-x}\mathbf{R}_{+z}$
& $18+x\,\,\,\mathbf{R}_{+y}$ & $10-z\,\,\,\mathbf{I}$
\\
$10$
& $(0,L,L)$
& $4+y\,\,\,\mathbf{R}_{+z}$ & $14+y\,\,\,\mathbf{R}_{-z}$
& $1-z\,\,\,\mathbf{R}_{-x}\mathbf{R}_{+y}$ & $6-y\,\,\,\mathbf{I}$
& $9+z\,\,\,\mathbf{I}$ & $15-z\,\,\,\mathbf{I}$
\\
$11$
& $(0,3L,3L)$
& $20+z\,\,\,\mathbf{R}_{-y}$ & $15+y\,\,\,\mathbf{R}_{-z}$
& $17-y\,\,\,\mathbf{R}^2_{+z}$ & $7-y\,\,\,\mathbf{I}$
& $6+z\,\,\,\mathbf{R}_{-z}$ & $12-z\,\,\,\mathbf{I}$
\\
$12$
& $(0,3L,4L)$
& $9-y\,\,\,\mathbf{R}_{-z}\mathbf{R}_{+x}$ & $5+x\,\,\,\mathbf{R}^2_{-y}\mathbf{R}_{-x}$
& $18-y\,\,\,\mathbf{R}^2_{+z}$ & $1-x\,\,\,\mathbf{R}_{-z}\mathbf{R}_{-y}$
& $11+z\,\,\,\mathbf{I}$ & $3+x\,\,\,\mathbf{R}_{-y}\mathbf{R}_{-z}$
\\
$13$
& $(2L,0,0)$
& $7-x\,\,\,\mathbf{R}^2_{-z}$ & $4-y\,\,\,\mathbf{R}_{+z}\mathbf{R}_{-x}$
& $16-x\,\,\,\mathbf{R}_{-x}\mathbf{R}_{+z}$ & $20+y\,\,\,\mathbf{R}_{-y}\mathbf{R}^2_{+z}$
& $2+y\,\,\,\mathbf{R}_{-x}\mathbf{R}_{+z}$ & $14-z\,\,\,\mathbf{I}$
\\
$14$
& $(2L,0,L)$
& $8-x\,\,\,\mathbf{R}^2_{-z}$ & $18-x\,\,\,\mathbf{I}$
& $5-z\,\,\,\mathbf{R}_{-x}$ & $10+x\,\,\,\mathbf{R}_{+z}$
& $13+z\,\,\,\mathbf{I}$ & $19-z\,\,\,\mathbf{R}_{+z}$
\\
$15$
& $(0,L,3L)$
& $4+z\,\,\,\mathbf{R}_{-y}\mathbf{R}_{+x}$ & $19-x\,\,\,\mathbf{I}$
& $1+x\,\,\,\mathbf{R}_{-z}$ & $11+x\,\,\,\mathbf{R}_{+z}$
& $10+z\,\,\,\mathbf{I}$ & $16-z\,\,\,\mathbf{I}$
\\
$16$
& $(0,L,4L)$
& $13-y\,\,\,\mathbf{R}_{-z}\mathbf{R}_{+x}$ & $9+y\,\,\,\mathbf{R}_{-z}\mathbf{R}_{+x}$
& $2+x\,\,\,\mathbf{R}_{-z}$ & $5-y\,\,\,\mathbf{R}_{-y}$
& $15+z\,\,\,\mathbf{I}$ & $7+y\,\,\,\mathbf{R}_{+x}$
\\
$17$
& $(4L,0,0)$
& $4-x\,\,\,\mathbf{R}^2_{-z}\mathbf{R}_{-x}$ & $20-y\,\,\,\mathbf{R}_{-y}\mathbf{R}_{+z}$
& $11-y\,\,\,\mathbf{R}^2_{-z}$ & $8+x\,\,\,\mathbf{R}_{+z}\mathbf{R}_{-y}$
& $6-x\,\,\,\mathbf{R}_{-y}\mathbf{R}_{+z}$ & $18-z\,\,\,\mathbf{I}$
\\
$18$
& $(4L,0,L)$
& $14+x\,\,\,\mathbf{I}$ & $9-z\,\,\,\mathbf{R}_{-y}$
& $12-y\,\,\,\mathbf{R}^2_{-z}$ & $2-y\,\,\,\mathbf{I}$
& $17+z\,\,\,\mathbf{I}$ & $3-z\,\,\,\mathbf{R}_{+z}$
\\
$19$
& $(2L,0,3L)$
& $15+x\,\,\,\mathbf{I}$ & $5+y\,\,\,\mathbf{R}_{-z}$
& $8+z\,\,\,\mathbf{R}_{+x}\mathbf{R}_{+y}$ & $3-y\,\,\,\mathbf{I}$
& $14+z\,\,\,\mathbf{R}_{-z}$ & $20-z\,\,\,\mathbf{I}$
\\
$20$
& $(2L,0,4L)$
& $9+x\,\,\,\mathbf{R}_{+x}$ & $6+y\,\,\,\mathbf{R}_{-z}$
& $17+x\,\,\,\mathbf{R}_{-z}\mathbf{R}_{+y}$ & $13+y\,\,\,\mathbf{R}^2_{-z}\mathbf{R}_{+y}$
& $19+z\,\,\,\mathbf{I}$ & $11-x\,\,\,\mathbf{R}_{+y}$
\\
\bottomrule
 \end{tabular}
\end{center}
\end{table}


\begin{table}[t]
  \scriptsize
  \renewcommand{\arraystretch}{1.2}
  \begin{center}
    \caption{Multicube representation of the Regina triangulation of
      the lens space L(5,2). Multicube Structure: region center
      locations $\vec c_A$, region face identifications, $\{\alpha
      \,A\} \leftrightarrow \{\beta\, B\}$ , and the rotation matrices
      for the associated interface maps, ${\bf C}_{A\alpha}^{B\beta}$.
    \label{t:L(5,2)}}
    \begin{tabular}{c|c|c|c|c|c|c|c}
      \toprule
      && $\alpha=-x$ & $\alpha=+x$ & $\alpha=-y$ & $\alpha=+y$ & $\alpha=-z$ &
      $\alpha=+z$ \\
      A  & $\vec c_A$
      &$B\,\,\beta\,\,\,{\mathbf C}_{A\alpha}^{B\beta}$
      &$B\,\,\beta\,\,\,{\mathbf C}_{A\alpha}^{B\beta}$
      &$B\,\,\beta\,\,\,{\mathbf C}_{A\alpha}^{B\beta}$
      &$B\,\,\beta\,\,\,{\mathbf C}_{A\alpha}^{B\beta}$
      &$B\,\,\beta\,\,\,{\mathbf C}_{A\alpha}^{B\beta}$
      &$B\,\,\beta\,\,\,{\mathbf C}_{A\alpha}^{B\beta}$\\
      \midrule
      $0.0$
      & $(0, 0, 0)$
      & $0.3+z\,\,\,\mathbf{R}_{-y} \mathbf{R}_{-x}$
      & $0.1-x\,\,\,\mathbf{I}$
      & $0.2-z\,\,\,\mathbf{R}_{-x} \mathbf{R}_{+y}$
      & $0.2-y\,\,\,\mathbf{I}$
      & $0.1-y\,\,\,\mathbf{R}_{+x} \mathbf{R}_{+z}$
      & $0.3-z\,\,\,\mathbf{I}$
      \\
      $0.1$
      & $(L, 0, 0)$
      & $0.0+x\,\,\,\mathbf{I}$
      & $0.2-x\,\,\,\mathbf{R}_{+x}$
      & $0.0-z\,\,\,\mathbf{R}_{-x} \mathbf{R}_{+y}$
      & $0.2+x\,\,\,\mathbf{R}_{+z}$
      & $0.3-y\,\,\,\mathbf{R}^2_{+z} \mathbf{R}_{-x}$
      & $0.3+x\,\,\,\mathbf{R}_{-y}$
      \\
      $0.2$
      & $(0, L, 0)$
      & $0.1+x\,\,\,\mathbf{R}_{-x}$
      & $0.1+y\,\,\,\mathbf{R}_{-z}$
      & $0.0+y\,\,\,\mathbf{I}$
      & $0.3-x\,\,\,\mathbf{R}^2_{+y} \mathbf{R}_{+z}$
      & $0.0-y\,\,\,\mathbf{R}_{+x} \mathbf{R}_{+z}$
      & $0.3+y\,\,\,\mathbf{R}_{+x}$
      \\
      $0.3$
      & $(0, 0, L)$
      & $0.2+y\,\,\,\mathbf{R}^2_{+x} \mathbf{R}_{-z}$
      & $0.1+z\,\,\,\mathbf{R}_{+y}$
      & $0.1-z\,\,\,\mathbf{R}^2_{+y} \mathbf{R}_{+x}$
      & $0.2+z\,\,\,\mathbf{R}_{-x}$
      & $0.0+z\,\,\,\mathbf{I}$
      & $0.0-x\,\,\,\mathbf{R}_{+y} \mathbf{R}_{+z}$
      \\
      \bottomrule
    \end{tabular}
  \end{center}
\end{table}

 
\begin{table}[htb] 
\scriptsize
  \renewcommand{\arraystretch}{1.5}
\begin{center}
\caption{Multicube representation of Seifert-Weber
  space~\cite{SeifertWeber1933}.  This multicube structure is based on
  cutting a dodecahedron into twenty cubes (each vertex of the
  dodecahedron is the vertex of one of the cubes, opposite vertices of
  these cubes all intersect at the center of the dodecahedron) and
  identifying opposite faces of the dodecahedron after rotation by
  $3\pi/5$.  Multicube Structure: region center locations $\vec c_A$,
  region face identifications, $\{\alpha \,A\} \leftrightarrow
  \{\beta\, B\}$ , and the rotation matrices for the associated
  interface maps, ${\bf C}_{A\alpha}^{B\beta}$.
   \label{t:TableSeifertWeber}}
\begin{tabular}{c|c|c|c|c|c|c|c}
  \toprule
  && $\alpha=-x$ & $\alpha=+x$ & $\alpha=-y$ & $\alpha=+y$ & $\alpha=-z$ &
  $\alpha=+z$ \\
A  & $\vec c_A$
&$B\,\,\beta\,\,\,{\mathbf C}_{A\alpha}^{B\beta}$
&$B\,\,\beta\,\,\,{\mathbf C}_{A\alpha}^{B\beta}$
&$B\,\,\beta\,\,\,{\mathbf C}_{A\alpha}^{B\beta}$
&$B\,\,\beta\,\,\,{\mathbf C}_{A\alpha}^{B\beta}$
&$B\,\,\beta\,\,\,{\mathbf C}_{A\alpha}^{B\beta}$
&$B\,\,\beta\,\,\,{\mathbf C}_{A\alpha}^{B\beta}$\\
\midrule
$1$
& $(2L,3L,0)$
& $12+y\,\,\,\mathbf{R}_{+y}\mathbf{R}_{+z}$ & $19-y\,\,\,\mathbf{R}_{-y}\mathbf{R}_{+z}$
& $8-y\,\,\,\mathbf{R}_{-y}\mathbf{R}^2_{+x}$ & $13-z\,\,\,\mathbf{R}^2_{+z}\mathbf{R}_{+x}$
& $9-y\,\,\,\mathbf{R}^2_{+y}\mathbf{R}_{+x}$ & $2-z\,\,\,\mathbf{I}$
\\ 
$2$
& $(2L,3L,L)$
& $6+x\,\,\,\mathbf{I}$ & $15-y\,\,\,\mathbf{R}_{-y}\mathbf{R}_{+z}$
& $18+y\,\,\,\mathbf{I}$ & $20+x\,\,\,\mathbf{R}_{+x}\mathbf{R}_{+z}$
& $1+z\,\,\,\mathbf{I}$ & $3-z\,\,\,\mathbf{I}$
\\
$3$
& $(2L,3L,2L)$
& $7+x\,\,\,\mathbf{I}$ & $5+y\,\,\,\mathbf{R}_{+y}\mathbf{R}_{-z}$
& $19+y\,\,\,\mathbf{I}$ & $10-x\,\,\,\mathbf{R}_{-x}\mathbf{R}_{-z}$
& $2+z\,\,\,\mathbf{I}$ & $4-z\,\,\,\mathbf{I}$
\\
$4$
& $(2L,3L,3L)$
& $17-x\,\,\,\mathbf{R}_{+x}\mathbf{R}^2_{+z}$ & $12+z\,\,\,\mathbf{R}^2_{+z}\mathbf{R}_{+y}$
& $13+x\,\,\,\mathbf{R}_{+x}\mathbf{R}_{-z}$ & $6-x\,\,\,\mathbf{R}_{-x}\mathbf{R}_{-z}$
& $3+z\,\,\,\mathbf{I}$ & $16-x\,\,\,\mathbf{R}^2_{+x}\mathbf{R}_{+y}$
\\
$5$
& $(0,3L,0)$
& $17-z\,\,\,\mathbf{R}^2_{+z}\mathbf{R}_{+y}$ & $12+x\,\,\,\mathbf{R}_{+x}\mathbf{R}^2_{+y}$
& $16+y\,\,\,\mathbf{R}_{+y}$ & $3+x\,\,\,\mathbf{R}_{+z}\mathbf{R}_{-y}$
& $13-y\,\,\,\mathbf{R}_{-y}\mathbf{R}_{+x}$ & $6-z\,\,\,\mathbf{I}$
\\
$6$
& $(0,3L,L)$
& $4+y\,\,\,\mathbf{R}_{+z}\mathbf{R}_{+x}$ & $2-x\,\,\,\mathbf{I}$
& $10+y\,\,\,\mathbf{I}$ & $19+x\,\,\,\mathbf{R}_{+x}\mathbf{R}_{+z}$
& $5+z\,\,\,\mathbf{I}$ & $7-z\,\,\,\mathbf{I}$
\\
$7$
& $(0,3L,2L)$
& $14-x\,\,\,\mathbf{R}_{+x}\mathbf{R}^2_{+y}$ & $3-x\,\,\,\mathbf{I}$
& $11+y\,\,\,\mathbf{I}$ & $9-x\,\,\,\mathbf{R}_{-x}\mathbf{R}_{-z}$
& $6+z\,\,\,\mathbf{I}$ & $8-z\,\,\,\mathbf{I}$
\\
$8$
& $(0,3L,3L)$
& $10-y\,\,\,\mathbf{R}_{-y}\mathbf{R}_{-z}$ & $17+y\,\,\,\mathbf{R}_{+y}\mathbf{R}_{-z}$
& $1-y\,\,\,\mathbf{R}^2_{-x}\mathbf{R}_{+y}$ & $16+z\,\,\,\mathbf{R}_{+z}\mathbf{R}_{-x}$
& $7+z\,\,\,\mathbf{I}$ & $20-y\,\,\,\mathbf{R}^2_{+y}\mathbf{R}_{-x}$
\\
$9$
& $(0,L,0)$
& $7+y\,\,\,\mathbf{R}_{+z}\mathbf{R}_{+x}$ & $20-x\,\,\,\mathbf{R}_{-x}$
& $1-z\,\,\,\mathbf{R}_{-x}\mathbf{R}^2_{-y}$ & $16+x\,\,\,\mathbf{R}_{-x}\mathbf{R}_{+z}$
& $17+x\,\,\,\mathbf{R}_{+x}\mathbf{R}_{+y}$ & $10-z\,\,\,\mathbf{I}$
\\
$10$
& $(0,L,L)$
& $3+y\,\,\,\mathbf{R}_{+z}\mathbf{R}_{+x}$ & $14+y\,\,\,\mathbf{R}_{-z}$
& $8-x\,\,\,\mathbf{R}_{+z}\mathbf{R}_{+y}$ & $6-y\,\,\,\mathbf{I}$
& $9+z\,\,\,\mathbf{I}$ & $11-z\,\,\,\mathbf{I}$
\\
$11$
& $(0,L,2L)$
& $13-x\,\,\,\mathbf{R}_{-x}\mathbf{R}^2_{+z}$ & $15+y\,\,\,\mathbf{R}_{-z}$
& $18-y\,\,\,\mathbf{R}_{-y}\mathbf{R}^2_{+z}$ & $7-y\,\,\,\mathbf{I}$
& $10+z\,\,\,\mathbf{I}$ & $12-z\,\,\,\mathbf{I}$
\\
$12$
& $(0,L,3L)$
& $20+z\,\,\,\mathbf{R}_{+z}\mathbf{R}_{-y}$ & $5+x\,\,\,\mathbf{R}^2_{-y}\mathbf{R}_{-x}$
& $14-y\,\,\,\mathbf{R}_{-y}\mathbf{R}^2_{+z}$ & $1-x\,\,\,\mathbf{R}_{-z}\mathbf{R}_{-y}$
& $11+z\,\,\,\mathbf{I}$ & $4+x\,\,\,\mathbf{R}_{-y}\mathbf{R}^2_{-z}$
\\
$13$
& $(2L,0,0)$
& $11-x\,\,\,\mathbf{R}^2_{-z}\mathbf{R}_{+x}$ & $4-y\,\,\,\mathbf{R}_{+z}\mathbf{R}_{-x}$
& $5-z\,\,\,\mathbf{R}_{-x}\mathbf{R}_{+y}$ & $20+y\,\,\,\mathbf{R}_{-y}\mathbf{R}^2_{+z}$
& $1+y\,\,\,\mathbf{R}_{-x}\mathbf{R}^2_{-z}$ & $14-z\,\,\,\mathbf{I}$
\\
$14$
& $(2L,0,L)$
& $7-x\,\,\,\mathbf{R}^2_{-y}\mathbf{R}_{-x}$ & $18-x\,\,\,\mathbf{I}$
& $12-y\,\,\,\mathbf{R}^2_{-z}\mathbf{R}_{+y}$ & $10+x\,\,\,\mathbf{R}_{+z}$
& $13+z\,\,\,\mathbf{I}$ & $15-z\,\,\,\mathbf{I}$
\\
$15$
& $(2L,0,2L)$
& $17-y\,\,\,\mathbf{R}_{-y}\mathbf{R}_{-z}$ & $19-x\,\,\,\mathbf{I}$
& $2+x\,\,\,\mathbf{R}_{-z}\mathbf{R}_{+y}$ & $11+x\,\,\,\mathbf{R}_{+z}$
& $14+z\,\,\,\mathbf{I}$ & $16-z\,\,\,\mathbf{I}$
\\
$16$
& $(2L,0,3L)$
& $4+z\,\,\,\mathbf{R}_{-y}\mathbf{R}^2_{-x}$ & $9+y\,\,\,\mathbf{R}_{-z}\mathbf{R}_{+x}$
& $18+x\,\,\,\mathbf{R}_{+x}\mathbf{R}_{-z}$ & $5-y\,\,\,\mathbf{R}_{-y}$
& $15+z\,\,\,\mathbf{I}$ & $8+y\,\,\,\mathbf{R}_{+z}\mathbf{R}_{-x}$
\\
$17$
& $(4L,0,0)$
& $4-x\,\,\,\mathbf{R}^2_{-z}\mathbf{R}_{-x}$ & $9-z\,\,\,\mathbf{R}_{-y}\mathbf{R}_{-x}$
& $15-x\,\,\,\mathbf{R}_{+z}\mathbf{R}_{+y}$ & $8+x\,\,\,\mathbf{R}_{+z}\mathbf{R}_{-y}$
& $5-x\,\,\,\mathbf{R}_{-y}\mathbf{R}^2_{-z}$ & $18-z\,\,\,\mathbf{I}$
\\
$18$
& $(4L,0,L)$
& $14+x\,\,\,\mathbf{I}$ & $16-y\,\,\,\mathbf{R}_{+z}\mathbf{R}_{-x}$
& $11-y\,\,\,\mathbf{R}^2_{-z}\mathbf{R}_{+y}$ & $2-y\,\,\,\mathbf{I}$
& $17+z\,\,\,\mathbf{I}$ & $19-z\,\,\,\mathbf{I}$
\\
$19$
& $(4L,0,2L)$
& $15+x\,\,\,\mathbf{I}$ & $6+y\,\,\,\mathbf{R}_{-z}\mathbf{R}_{-x}$
& $1+x\,\,\,\mathbf{R}_{-z}\mathbf{R}_{+y}$ & $3-y\,\,\,\mathbf{I}$
& $18+z\,\,\,\mathbf{I}$ & $20-z\,\,\,\mathbf{I}$
\\
$20$
& $(4L,0,3L)$
& $9+x\,\,\,\mathbf{R}_{+x}$ & $2+y\,\,\,\mathbf{R}_{-z}\mathbf{R}_{-x}$
& $8+x\,\,\,\mathbf{R}_{+x}\mathbf{R}^2_{-y}$ & $13+y\,\,\,\mathbf{R}^2_{-z}\mathbf{R}_{+y}$
& $19+z\,\,\,\mathbf{I}$ & $12-x\,\,\,\mathbf{R}_{+y}\mathbf{R}_{-z}$
\\
\bottomrule
 \end{tabular}
\end{center}
\end{table}


\begin{table}[t]
  \scriptsize
  \renewcommand{\arraystretch}{1.2}
  \begin{center}
    \caption{Multicube representation of the Regina triangulation of
      the Seifert fiber space SFS[S2:(2,1)(2,1)(2,-1)].  Multicube
      Structure: region center locations $\vec c_A$, region face
      identifications, $\{\alpha \,A\} \leftrightarrow \{\beta\, B\}$
      , and the rotation matrices for the associated interface maps,
      ${\bf C}_{A\alpha}^{B\beta}$.
      \label{t:SFS[S2:(2,1)(2,1)(2,-1)]}}
    \begin{tabular}{c|c|c|c|c|c|c|c}
      \toprule
      && $\alpha=-x$ & $\alpha=+x$ & $\alpha=-y$ & $\alpha=+y$ & $\alpha=-z$ &
      $\alpha=+z$ \\
      A  & $\vec c_A$
      &$B\,\,\beta\,\,\,{\mathbf C}_{A\alpha}^{B\beta}$
      &$B\,\,\beta\,\,\,{\mathbf C}_{A\alpha}^{B\beta}$
      &$B\,\,\beta\,\,\,{\mathbf C}_{A\alpha}^{B\beta}$
      &$B\,\,\beta\,\,\,{\mathbf C}_{A\alpha}^{B\beta}$
      &$B\,\,\beta\,\,\,{\mathbf C}_{A\alpha}^{B\beta}$
      &$B\,\,\beta\,\,\,{\mathbf C}_{A\alpha}^{B\beta}$\\
      \midrule
      $0.0$
      & $(0, 0, 0)$
      & $1.2-x\,\,\,\mathbf{R}^2_{+z}$
      & $0.1-x\,\,\,\mathbf{I}$
      & $1.1-z\,\,\,\mathbf{R}_{-x} \mathbf{R}_{-y}$
      & $0.2-y\,\,\,\mathbf{I}$
      & $1.3-y\,\,\,\mathbf{R}_{+x} \mathbf{R}_{-z}$
      & $0.3-z\,\,\,\mathbf{I}$
      \\
      $0.1$
      & $(L, 0, 0)$
      & $0.0+x\,\,\,\mathbf{I}$
      & $1.3+z\,\,\,\mathbf{R}_{+y}$
      & $1.2-z\,\,\,\mathbf{R}^2_{+y} \mathbf{R}_{+x}$
      & $0.2+x\,\,\,\mathbf{R}_{+z}$
      & $1.0-y\,\,\,\mathbf{R}_{+x} \mathbf{R}_{-z}$
      & $0.3+x\,\,\,\mathbf{R}_{-y}$
      \\
      $0.2$
      & $(0, L, 0)$
      & $1.0-x\,\,\,\mathbf{R}^2_{+z}$
      & $0.1+y\,\,\,\mathbf{R}_{-z}$
      & $0.0+y\,\,\,\mathbf{I}$
      & $1.2+y\,\,\,\mathbf{R}^2_{+x} \mathbf{R}_{+y}$
      & $1.1-y\,\,\,\mathbf{R}^2_{+z} \mathbf{R}_{-x}$
      & $0.3+y\,\,\,\mathbf{R}_{+x}$
      \\
      $0.3$
      & $(0, 0, L)$
      & $1.3-x\,\,\,\mathbf{R}^2_{+z} \mathbf{R}_{-x}$
      & $0.1+z\,\,\,\mathbf{R}_{+y}$
      & $1.0-z\,\,\,\mathbf{R}_{-x} \mathbf{R}_{-y}$
      & $0.2+z\,\,\,\mathbf{R}_{-x}$
      & $0.0+z\,\,\,\mathbf{I}$
      & $1.1+x\,\,\,\mathbf{R}_{-y}$
      \\
      $1.0$
      & $(3L, 0, 0)$
      & $0.2-x\,\,\,\mathbf{R}^2_{+z}$
      & $1.1-x\,\,\,\mathbf{I}$
      & $0.1-z\,\,\,\mathbf{R}_{-x} \mathbf{R}_{-y}$
      & $1.2-y\,\,\,\mathbf{I}$
      & $0.3-y\,\,\,\mathbf{R}_{+x} \mathbf{R}_{-z}$
      & $1.3-z\,\,\,\mathbf{I}$
      \\
      $1.1$
      & $(4L, 0, 0)$
      & $1.0+x\,\,\,\mathbf{I}$
      & $0.3+z\,\,\,\mathbf{R}_{+y}$
      & $0.2-z\,\,\,\mathbf{R}^2_{+y} \mathbf{R}_{+x}$
      & $1.2+x\,\,\,\mathbf{R}_{+z}$
      & $0.0-y\,\,\,\mathbf{R}_{+x} \mathbf{R}_{-z}$
      & $1.3+x\,\,\,\mathbf{R}_{-y}$
      \\
      $1.2$
      & $(3L, L, 0)$
      & $0.0-x\,\,\,\mathbf{R}^2_{+z}$
      & $1.1+y\,\,\,\mathbf{R}_{-z}$
      & $1.0+y\,\,\,\mathbf{I}$
      & $0.2+y\,\,\,\mathbf{R}^2_{+x} \mathbf{R}_{+y}$
      & $0.1-y\,\,\,\mathbf{R}^2_{+z} \mathbf{R}_{-x}$
      & $1.3+y\,\,\,\mathbf{R}_{+x}$
      \\
      $1.3$
      & $(3L, 0, L)$
      & $0.3-x\,\,\,\mathbf{R}^2_{+z} \mathbf{R}_{-x}$
      & $1.1+z\,\,\,\mathbf{R}_{+y}$
      & $0.0-z\,\,\,\mathbf{R}_{-x} \mathbf{R}_{-y}$
      & $1.2+z\,\,\,\mathbf{R}_{-x}$
      & $1.0+z\,\,\,\mathbf{I}$
      & $0.1+x\,\,\,\mathbf{R}_{-y}$
      \\
      \bottomrule
    \end{tabular}
  \end{center}
\end{table}


\begin{table}[t]
  \scriptsize
  \renewcommand{\arraystretch}{1.2}
  \begin{center}
    \caption{Multicube representation of the Regina triangulation of
      the Seifert fiber space KB/n2$\times\sim$S1. Multicube Structure: 
      region center locations $\vec c_A$, region face identifications,
      $\{\alpha \,A\} \leftrightarrow \{\beta\, B\}$ , and the rotation
      matrices for the associated interface maps, ${\bf C}_{A\alpha}^{B\beta}$.
    \label{t:KB}
    }
    \begin{tabular}{c|c|c|c|c|c|c|c}
      \toprule
      && $\alpha=-x$ & $\alpha=+x$ & $\alpha=-y$ & $\alpha=+y$ & $\alpha=-z$ &
      $\alpha=+z$ \\
      A  & $\vec c_A$
      &$B\,\,\beta\,\,\,{\mathbf C}_{A\alpha}^{B\beta}$
      &$B\,\,\beta\,\,\,{\mathbf C}_{A\alpha}^{B\beta}$
      &$B\,\,\beta\,\,\,{\mathbf C}_{A\alpha}^{B\beta}$
      &$B\,\,\beta\,\,\,{\mathbf C}_{A\alpha}^{B\beta}$
      &$B\,\,\beta\,\,\,{\mathbf C}_{A\alpha}^{B\beta}$
      &$B\,\,\beta\,\,\,{\mathbf C}_{A\alpha}^{B\beta}$\\
      \midrule
      $0.0$
      & $(0, 0, 0)$
      & $2.1+x\,\,\,\mathbf{I}$
      & $0.1-x\,\,\,\mathbf{I}$
      & $3.1+x\,\,\,\mathbf{R}_{-z} \mathbf{R}_{-y}$
      & $0.2-y\,\,\,\mathbf{I}$
      & $4.3+z\,\,\,\mathbf{I}$
      & $0.3-z\,\,\,\mathbf{I}$
      \\
      $0.1$
      & $(L, 0, 0)$
      & $0.0+x\,\,\,\mathbf{I}$
      & $1.1+x\,\,\,\mathbf{R}^2_{+z} \mathbf{R}_{+x}$
      & $3.3+z\,\,\,\mathbf{R}^2_{+y} \mathbf{R}_{-x}$
      & $0.2+x\,\,\,\mathbf{R}_{+z}$
      & $4.1+x\,\,\,\mathbf{R}_{+y}$
      & $0.3+x\,\,\,\mathbf{R}_{-y}$
      \\
      $0.2$
      & $(0, L, 0)$
      & $2.2+y\,\,\,\mathbf{R}_{+z}$
      & $0.1+y\,\,\,\mathbf{R}_{-z}$
      & $0.0+y\,\,\,\mathbf{I}$
      & $1.3+z\,\,\,\mathbf{R}_{-x}$
      & $4.2+y\,\,\,\mathbf{R}_{-x}$
      & $0.3+y\,\,\,\mathbf{R}_{+x}$
      \\
      $0.3$
      & $(0, 0, L)$
      & $2.3+z\,\,\,\mathbf{R}_{-y}$
      & $0.1+z\,\,\,\mathbf{R}_{+y}$
      & $3.2+y\,\,\,\mathbf{R}_{-y}$
      & $0.2+z\,\,\,\mathbf{R}_{-x}$
      & $0.0+z\,\,\,\mathbf{I}$
      & $1.2+y\,\,\,\mathbf{R}_{+x}$
      \\
      $1.0$
      & $(0, 3L, 0)$
      & $2.2-x\,\,\,\mathbf{R}^2_{+z}$
      & $1.1-x\,\,\,\mathbf{I}$
      & $5.2+y\,\,\,\mathbf{I}$
      & $1.2-y\,\,\,\mathbf{I}$
      & $3.3-x\,\,\,\mathbf{R}_{-y} \mathbf{R}_{+z}$
      & $1.3-z\,\,\,\mathbf{I}$
      \\
      $1.1$
      & $(L, 3L, 0)$
      & $1.0+x\,\,\,\mathbf{I}$
      & $0.1+x\,\,\,\mathbf{R}^2_{+z} \mathbf{R}_{+x}$
      & $5.1+x\,\,\,\mathbf{R}_{-z}$
      & $1.2+x\,\,\,\mathbf{R}_{+z}$
      & $3.2-x\,\,\,\mathbf{R}^2_{+z} \mathbf{R}_{+y}$
      & $1.3+x\,\,\,\mathbf{R}_{-y}$
      \\
      $1.2$
      & $(0, 4L, 0)$
      & $2.0-x\,\,\,\mathbf{R}^2_{+z}$
      & $1.1+y\,\,\,\mathbf{R}_{-z}$
      & $1.0+y\,\,\,\mathbf{I}$
      & $0.3+z\,\,\,\mathbf{R}_{-x}$
      & $3.0-x\,\,\,\mathbf{R}_{-y} \mathbf{R}_{+z}$
      & $1.3+y\,\,\,\mathbf{R}_{+x}$
      \\
      $1.3$
      & $(0, 3L, L)$
      & $2.3-x\,\,\,\mathbf{R}^2_{+z} \mathbf{R}_{-x}$
      & $1.1+z\,\,\,\mathbf{R}_{+y}$
      & $5.3+z\,\,\,\mathbf{R}_{+x}$
      & $1.2+z\,\,\,\mathbf{R}_{-x}$
      & $1.0+z\,\,\,\mathbf{I}$
      & $0.2+y\,\,\,\mathbf{R}_{+x}$
      \\
      $2.0$
      & $(3L, 0, 0)$
      & $1.2-x\,\,\,\mathbf{R}^2_{+z}$
      & $2.1-x\,\,\,\mathbf{I}$
      & $5.3-x\,\,\,\mathbf{R}_{+z} \mathbf{R}_{+y}$
      & $2.2-y\,\,\,\mathbf{I}$
      & $4.3-x\,\,\,\mathbf{R}_{-y} \mathbf{R}_{+z}$
      & $2.3-z\,\,\,\mathbf{I}$
      \\
      $2.1$
      & $(4L, 0, 0)$
      & $2.0+x\,\,\,\mathbf{I}$
      & $0.0-x\,\,\,\mathbf{I}$
      & $5.0-x\,\,\,\mathbf{R}_{+z} \mathbf{R}_{+y}$
      & $2.2+x\,\,\,\mathbf{R}_{+z}$
      & $4.2-x\,\,\,\mathbf{R}^2_{+z} \mathbf{R}_{+y}$
      & $2.3+x\,\,\,\mathbf{R}_{-y}$
      \\
      $2.2$
      & $(3L, L, 0)$
      & $1.0-x\,\,\,\mathbf{R}^2_{+z}$
      & $2.1+y\,\,\,\mathbf{R}_{-z}$
      & $2.0+y\,\,\,\mathbf{I}$
      & $0.2-x\,\,\,\mathbf{R}_{-z}$
      & $4.0-x\,\,\,\mathbf{R}_{-y} \mathbf{R}_{+z}$
      & $2.3+y\,\,\,\mathbf{R}_{+x}$
      \\
      $2.3$
      & $(3L, 0, L)$
      & $1.3-x\,\,\,\mathbf{R}^2_{+z} \mathbf{R}_{-x}$
      & $2.1+z\,\,\,\mathbf{R}_{+y}$
      & $5.2-x\,\,\,\mathbf{R}^2_{+y} \mathbf{R}_{-z}$
      & $2.2+z\,\,\,\mathbf{R}_{-x}$
      & $2.0+z\,\,\,\mathbf{I}$
      & $0.3-x\,\,\,\mathbf{R}_{+y}$
      \\
      $3.0$
      & $(3L, 3L, 0)$
      & $1.2-z\,\,\,\mathbf{R}_{+y} \mathbf{R}_{+x}$
      & $3.1-x\,\,\,\mathbf{I}$
      & $4.3-y\,\,\,\mathbf{R}^2_{+x}$
      & $3.2-y\,\,\,\mathbf{I}$
      & $5.1-z\,\,\,\mathbf{R}^2_{+y}$
      & $3.3-z\,\,\,\mathbf{I}$
      \\
      $3.1$
      & $(4L, 3L, 0)$
      & $3.0+x\,\,\,\mathbf{I}$
      & $0.0-y\,\,\,\mathbf{R}_{+z} \mathbf{R}_{+x}$
      & $4.1-y\,\,\,\mathbf{R}^2_{+x} \mathbf{R}_{-y}$
      & $3.2+x\,\,\,\mathbf{R}_{+z}$
      & $5.0-z\,\,\,\mathbf{R}^2_{+y}$
      & $3.3+x\,\,\,\mathbf{R}_{-y}$
      \\
      $3.2$
      & $(3L, 4L, 0)$
      & $1.1-z\,\,\,\mathbf{R}^2_{+x} \mathbf{R}_{-y}$
      & $3.1+y\,\,\,\mathbf{R}_{-z}$
      & $3.0+y\,\,\,\mathbf{I}$
      & $0.3-y\,\,\,\mathbf{R}_{+y}$
      & $5.2-z\,\,\,\mathbf{R}^2_{+y} \mathbf{R}_{-z}$
      & $3.3+y\,\,\,\mathbf{R}_{+x}$
      \\
      $3.3$
      & $(3L, 3L, L)$
      & $1.0-z\,\,\,\mathbf{R}_{+y} \mathbf{R}_{+x}$
      & $3.1+z\,\,\,\mathbf{R}_{+y}$
      & $4.0-y\,\,\,\mathbf{R}^2_{+x}$
      & $3.2+z\,\,\,\mathbf{R}_{-x}$
      & $3.0+z\,\,\,\mathbf{I}$
      & $0.1-y\,\,\,\mathbf{R}^2_{+z} \mathbf{R}_{+x}$
      \\
      $4.0$
      & $(6L, 0, 0)$
      & $2.2-z\,\,\,\mathbf{R}_{+y} \mathbf{R}_{+x}$
      & $4.1-x\,\,\,\mathbf{I}$
      & $3.3-y\,\,\,\mathbf{R}^2_{+x}$
      & $4.2-y\,\,\,\mathbf{I}$
      & $5.0-y\,\,\,\mathbf{R}_{+x}$
      & $4.3-z\,\,\,\mathbf{I}$
      \\
      $4.1$
      & $(7L, 0, 0)$
      & $4.0+x\,\,\,\mathbf{I}$
      & $0.1-z\,\,\,\mathbf{R}_{-y}$
      & $3.1-y\,\,\,\mathbf{R}^2_{+x} \mathbf{R}_{-y}$
      & $4.2+x\,\,\,\mathbf{R}_{+z}$
      & $5.1-y\,\,\,\mathbf{R}_{+x}$
      & $4.3+x\,\,\,\mathbf{R}_{-y}$
      \\
      $4.2$
      & $(6L, L, 0)$
      & $2.1-z\,\,\,\mathbf{R}^2_{+x} \mathbf{R}_{-y}$
      & $4.1+y\,\,\,\mathbf{R}_{-z}$
      & $4.0+y\,\,\,\mathbf{I}$
      & $0.2-z\,\,\,\mathbf{R}_{+x}$
      & $5.3-y\,\,\,\mathbf{R}_{+x}$
      & $4.3+y\,\,\,\mathbf{R}_{+x}$
      \\
      $4.3$
      & $(6L, 0, L)$
      & $2.0-z\,\,\,\mathbf{R}_{+y} \mathbf{R}_{+x}$
      & $4.1+z\,\,\,\mathbf{R}_{+y}$
      & $3.0-y\,\,\,\mathbf{R}^2_{+x}$
      & $4.2+z\,\,\,\mathbf{R}_{-x}$
      & $4.0+z\,\,\,\mathbf{I}$
      & $0.0-z\,\,\,\mathbf{I}$
      \\
      $5.0$
      & $(6L, 3L, 0)$
      & $2.1-y\,\,\,\mathbf{R}_{-z} \mathbf{R}_{+x}$
      & $5.1-x\,\,\,\mathbf{I}$
      & $4.0-z\,\,\,\mathbf{R}_{-x}$
      & $5.2-y\,\,\,\mathbf{I}$
      & $3.1-z\,\,\,\mathbf{R}^2_{+y}$
      & $5.3-z\,\,\,\mathbf{I}$
      \\
      $5.1$
      & $(7L, 3L, 0)$
      & $5.0+x\,\,\,\mathbf{I}$
      & $1.1-y\,\,\,\mathbf{R}_{+z}$
      & $4.1-z\,\,\,\mathbf{R}_{-x}$
      & $5.2+x\,\,\,\mathbf{R}_{+z}$
      & $3.0-z\,\,\,\mathbf{R}^2_{+y}$
      & $5.3+x\,\,\,\mathbf{R}_{-y}$
      \\
      $5.2$
      & $(6L, 4L, 0)$
      & $2.3-y\,\,\,\mathbf{R}^2_{+x} \mathbf{R}_{+z}$
      & $5.1+y\,\,\,\mathbf{R}_{-z}$
      & $5.0+y\,\,\,\mathbf{I}$
      & $1.0-y\,\,\,\mathbf{I}$
      & $3.2-z\,\,\,\mathbf{R}^2_{+y} \mathbf{R}_{-z}$
      & $5.3+y\,\,\,\mathbf{R}_{+x}$
      \\
      $5.3$
      & $(6L, 3L, L)$
      & $2.0-y\,\,\,\mathbf{R}_{-z} \mathbf{R}_{+x}$
      & $5.1+z\,\,\,\mathbf{R}_{+y}$
      & $4.2-z\,\,\,\mathbf{R}_{-x}$
      & $5.2+z\,\,\,\mathbf{R}_{-x}$
      & $5.0+z\,\,\,\mathbf{I}$
      & $1.3-y\,\,\,\mathbf{R}_{-x}$
      \\
      \bottomrule
    \end{tabular}
  \end{center}
\end{table}


\begin{table}[t]
  \scriptsize
  \renewcommand{\arraystretch}{1.2}
  \begin{center}
    \caption{Multicube representation of the Regina triangulation of
      the Seifert fiber space SFS[RP2/n2:(2,1)(2,-1)]. Multicube
        Structure: region center locations $\vec c_A$, region face
      identifications, $\{\alpha \,A\} \leftrightarrow \{\beta\, B\}$
      , and the rotation matrices for the associated interface maps,
      ${\bf C}_{A\alpha}^{B\beta}$.
      \label{t:SFSRP2}
    }
    \begin{tabular}{c|c|c|c|c|c|c|c}
      \toprule
      && $\alpha=-x$ & $\alpha=+x$ & $\alpha=-y$ & $\alpha=+y$ & $\alpha=-z$ &
      $\alpha=+z$ \\
      A  & $\vec c_A$
      &$B\,\,\beta\,\,\,{\mathbf C}_{A\alpha}^{B\beta}$
      &$B\,\,\beta\,\,\,{\mathbf C}_{A\alpha}^{B\beta}$
      &$B\,\,\beta\,\,\,{\mathbf C}_{A\alpha}^{B\beta}$
      &$B\,\,\beta\,\,\,{\mathbf C}_{A\alpha}^{B\beta}$
      &$B\,\,\beta\,\,\,{\mathbf C}_{A\alpha}^{B\beta}$
      &$B\,\,\beta\,\,\,{\mathbf C}_{A\alpha}^{B\beta}$\\
      \midrule
      $0.0$
      & $(0, 0, 0)$
      & $2.1+x\,\,\,\mathbf{I}$
      & $0.1-x\,\,\,\mathbf{I}$
      & $3.1+x\,\,\,\mathbf{R}_{-z} \mathbf{R}_{-y}$
      & $0.2-y\,\,\,\mathbf{I}$
      & $4.2+y\,\,\,\mathbf{R}_{-x} \mathbf{R}_{-z}$
      & $0.3-z\,\,\,\mathbf{I}$
      \\
      $0.1$
      & $(L, 0, 0)$
      & $0.0+x\,\,\,\mathbf{I}$
      & $1.1+x\,\,\,\mathbf{R}^2_{+z} \mathbf{R}_{+x}$
      & $3.3+z\,\,\,\mathbf{R}^2_{+y} \mathbf{R}_{-x}$
      & $0.2+x\,\,\,\mathbf{R}_{+z}$
      & $4.3+z\,\,\,\mathbf{R}_{-z}$
      & $0.3+x\,\,\,\mathbf{R}_{-y}$
      \\
      $0.2$
      & $(0, L, 0)$
      & $2.2+y\,\,\,\mathbf{R}_{+z}$
      & $0.1+y\,\,\,\mathbf{R}_{-z}$
      & $0.0+y\,\,\,\mathbf{I}$
      & $1.3+z\,\,\,\mathbf{R}_{-x}$
      & $4.1+x\,\,\,\mathbf{R}^2_{+z} \mathbf{R}_{-y}$
      & $0.3+y\,\,\,\mathbf{R}_{+x}$
      \\
      $0.3$
      & $(0, 0, L)$
      & $2.3+z\,\,\,\mathbf{R}_{-y}$
      & $0.1+z\,\,\,\mathbf{R}_{+y}$
      & $3.2+y\,\,\,\mathbf{R}_{-y}$
      & $0.2+z\,\,\,\mathbf{R}_{-x}$
      & $0.0+z\,\,\,\mathbf{I}$
      & $1.2+y\,\,\,\mathbf{R}_{+x}$
      \\
      $1.0$
      & $(0, 3L, 0)$
      & $2.2-x\,\,\,\mathbf{R}^2_{+z}$
      & $1.1-x\,\,\,\mathbf{I}$
      & $5.2+y\,\,\,\mathbf{I}$
      & $1.2-y\,\,\,\mathbf{I}$
      & $3.3-x\,\,\,\mathbf{R}_{-y} \mathbf{R}_{+z}$
      & $1.3-z\,\,\,\mathbf{I}$
      \\
      $1.1$
      & $(L, 3L, 0)$
      & $1.0+x\,\,\,\mathbf{I}$
      & $0.1+x\,\,\,\mathbf{R}^2_{+z} \mathbf{R}_{+x}$
      & $5.1+x\,\,\,\mathbf{R}_{-z}$
      & $1.2+x\,\,\,\mathbf{R}_{+z}$
      & $3.2-x\,\,\,\mathbf{R}^2_{+z} \mathbf{R}_{+y}$
      & $1.3+x\,\,\,\mathbf{R}_{-y}$
      \\
      $1.2$
      & $(0, 4L, 0)$
      & $2.0-x\,\,\,\mathbf{R}^2_{+z}$
      & $1.1+y\,\,\,\mathbf{R}_{-z}$
      & $1.0+y\,\,\,\mathbf{I}$
      & $0.3+z\,\,\,\mathbf{R}_{-x}$
      & $3.0-x\,\,\,\mathbf{R}_{-y} \mathbf{R}_{+z}$
      & $1.3+y\,\,\,\mathbf{R}_{+x}$
      \\
      $1.3$
      & $(0, 3L, L)$
      & $2.3-x\,\,\,\mathbf{R}^2_{+z} \mathbf{R}_{-x}$
      & $1.1+z\,\,\,\mathbf{R}_{+y}$
      & $5.3+z\,\,\,\mathbf{R}_{+x}$
      & $1.2+z\,\,\,\mathbf{R}_{-x}$
      & $1.0+z\,\,\,\mathbf{I}$
      & $0.2+y\,\,\,\mathbf{R}_{+x}$
      \\
      $2.0$
      & $(3L, 0, 0)$
      & $1.2-x\,\,\,\mathbf{R}^2_{+z}$
      & $2.1-x\,\,\,\mathbf{I}$
      & $5.3-x\,\,\,\mathbf{R}_{+z} \mathbf{R}_{+y}$
      & $2.2-y\,\,\,\mathbf{I}$
      & $4.0-x\,\,\,\mathbf{R}_{-y}$
      & $2.3-z\,\,\,\mathbf{I}$
      \\
      $2.1$
      & $(4L, 0, 0)$
      & $2.0+x\,\,\,\mathbf{I}$
      & $0.0-x\,\,\,\mathbf{I}$
      & $5.0-x\,\,\,\mathbf{R}_{+z} \mathbf{R}_{+y}$
      & $2.2+x\,\,\,\mathbf{R}_{+z}$
      & $4.3-x\,\,\,\mathbf{R}_{-y}$
      & $2.3+x\,\,\,\mathbf{R}_{-y}$
      \\
      $2.2$
      & $(3L, L, 0)$
      & $1.0-x\,\,\,\mathbf{R}^2_{+z}$
      & $2.1+y\,\,\,\mathbf{R}_{-z}$
      & $2.0+y\,\,\,\mathbf{I}$
      & $0.2-x\,\,\,\mathbf{R}_{-z}$
      & $4.2-x\,\,\,\mathbf{R}_{-y}$
      & $2.3+y\,\,\,\mathbf{R}_{+x}$
      \\
      $2.3$
      & $(3L, 0, L)$
      & $1.3-x\,\,\,\mathbf{R}^2_{+z} \mathbf{R}_{-x}$
      & $2.1+z\,\,\,\mathbf{R}_{+y}$
      & $5.2-x\,\,\,\mathbf{R}^2_{+y} \mathbf{R}_{-z}$
      & $2.2+z\,\,\,\mathbf{R}_{-x}$
      & $2.0+z\,\,\,\mathbf{I}$
      & $0.3-x\,\,\,\mathbf{R}_{+y}$
      \\
      $3.0$
      & $(3L, 3L, 0)$
      & $1.2-z\,\,\,\mathbf{R}_{+y} \mathbf{R}_{+x}$
      & $3.1-x\,\,\,\mathbf{I}$
      & $4.0-z\,\,\,\mathbf{R}_{-x}$
      & $3.2-y\,\,\,\mathbf{I}$
      & $5.1-z\,\,\,\mathbf{R}^2_{+y}$
      & $3.3-z\,\,\,\mathbf{I}$
      \\
      $3.1$
      & $(4L, 3L, 0)$
      & $3.0+x\,\,\,\mathbf{I}$
      & $0.0-y\,\,\,\mathbf{R}_{+z} \mathbf{R}_{+x}$
      & $4.1-z\,\,\,\mathbf{R}_{-x}$
      & $3.2+x\,\,\,\mathbf{R}_{+z}$
      & $5.0-z\,\,\,\mathbf{R}^2_{+y}$
      & $3.3+x\,\,\,\mathbf{R}_{-y}$
      \\
      $3.2$
      & $(3L, 4L, 0)$
      & $1.1-z\,\,\,\mathbf{R}^2_{+x} \mathbf{R}_{-y}$
      & $3.1+y\,\,\,\mathbf{R}_{-z}$
      & $3.0+y\,\,\,\mathbf{I}$
      & $0.3-y\,\,\,\mathbf{R}_{+y}$
      & $5.2-z\,\,\,\mathbf{R}^2_{+y} \mathbf{R}_{-z}$
      & $3.3+y\,\,\,\mathbf{R}_{+x}$
      \\
      $3.3$
      & $(3L, 3L, L)$
      & $1.0-z\,\,\,\mathbf{R}_{+y} \mathbf{R}_{+x}$
      & $3.1+z\,\,\,\mathbf{R}_{+y}$
      & $4.2-z\,\,\,\mathbf{R}_{-x}$
      & $3.2+z\,\,\,\mathbf{R}_{-x}$
      & $3.0+z\,\,\,\mathbf{I}$
      & $0.1-y\,\,\,\mathbf{R}^2_{+z} \mathbf{R}_{+x}$
      \\
      $4.0$
      & $(6L, 0, 0)$
      & $2.0-z\,\,\,\mathbf{R}_{+y}$
      & $4.1-x\,\,\,\mathbf{I}$
      & $5.0-y\,\,\,\mathbf{R}^2_{+x} \mathbf{R}_{+y}$
      & $4.2-y\,\,\,\mathbf{I}$
      & $3.0-y\,\,\,\mathbf{R}_{+x}$
      & $4.3-z\,\,\,\mathbf{I}$
      \\
      $4.1$
      & $(7L, 0, 0)$
      & $4.0+x\,\,\,\mathbf{I}$
      & $0.2-z\,\,\,\mathbf{R}^2_{+x} \mathbf{R}_{+y}$
      & $5.3-y\,\,\,\mathbf{R}^2_{+x} \mathbf{R}_{+y}$
      & $4.2+x\,\,\,\mathbf{R}_{+z}$
      & $3.1-y\,\,\,\mathbf{R}_{+x}$
      & $4.3+x\,\,\,\mathbf{R}_{-y}$
      \\
      $4.2$
      & $(6L, L, 0)$
      & $2.2-z\,\,\,\mathbf{R}_{+y}$
      & $4.1+y\,\,\,\mathbf{R}_{-z}$
      & $4.0+y\,\,\,\mathbf{I}$
      & $0.0-z\,\,\,\mathbf{R}_{+x} \mathbf{R}_{+y}$
      & $3.3-y\,\,\,\mathbf{R}_{+x}$
      & $4.3+y\,\,\,\mathbf{R}_{+x}$
      \\
      $4.3$
      & $(6L, 0, L)$
      & $2.1-z\,\,\,\mathbf{R}_{+y}$
      & $4.1+z\,\,\,\mathbf{R}_{+y}$
      & $5.1-y\,\,\,\mathbf{R}^2_{+x} \mathbf{R}_{+y}$
      & $4.2+z\,\,\,\mathbf{R}_{-x}$
      & $4.0+z\,\,\,\mathbf{I}$
      & $0.1-z\,\,\,\mathbf{R}_{+z}$
      \\
      $5.0$
      & $(6L, 3L, 0)$
      & $2.1-y\,\,\,\mathbf{R}_{-z} \mathbf{R}_{+x}$
      & $5.1-x\,\,\,\mathbf{I}$
      & $4.0-y\,\,\,\mathbf{R}^2_{+x} \mathbf{R}_{+y}$
      & $5.2-y\,\,\,\mathbf{I}$
      & $3.1-z\,\,\,\mathbf{R}^2_{+y}$
      & $5.3-z\,\,\,\mathbf{I}$
      \\
      $5.1$
      & $(7L, 3L, 0)$
      & $5.0+x\,\,\,\mathbf{I}$
      & $1.1-y\,\,\,\mathbf{R}_{+z}$
      & $4.3-y\,\,\,\mathbf{R}^2_{+x} \mathbf{R}_{+y}$
      & $5.2+x\,\,\,\mathbf{R}_{+z}$
      & $3.0-z\,\,\,\mathbf{R}^2_{+y}$
      & $5.3+x\,\,\,\mathbf{R}_{-y}$
      \\
      $5.2$
      & $(6L, 4L, 0)$
      & $2.3-y\,\,\,\mathbf{R}^2_{+x} \mathbf{R}_{+z}$
      & $5.1+y\,\,\,\mathbf{R}_{-z}$
      & $5.0+y\,\,\,\mathbf{I}$
      & $1.0-y\,\,\,\mathbf{I}$
      & $3.2-z\,\,\,\mathbf{R}^2_{+y} \mathbf{R}_{-z}$
      & $5.3+y\,\,\,\mathbf{R}_{+x}$
      \\
      $5.3$
      & $(6L, 3L, L)$
      & $2.0-y\,\,\,\mathbf{R}_{-z} \mathbf{R}_{+x}$
      & $5.1+z\,\,\,\mathbf{R}_{+y}$
      & $4.1-y\,\,\,\mathbf{R}^2_{+x} \mathbf{R}_{+y}$
      & $5.2+z\,\,\,\mathbf{R}_{-x}$
      & $5.0+z\,\,\,\mathbf{I}$
      & $1.3-y\,\,\,\mathbf{R}_{-x}$
      \\
      \bottomrule
    \end{tabular}
  \end{center}
\end{table}

\clearpage
\vspace{0.2cm}
\bibliographystyle{model1-num-names}
\bibliography{../../References/References}
\end{document}